\documentclass[mmnp]{edpsmath}
\usepackage{graphicx}
\usepackage{amssymb}
\usepackage{amsmath}
\usepackage{subfigure}
\usepackage{tabularx}
\usepackage[section]{placeins}
\usepackage{float}
\numberwithin{equation}{section}

\begin{document}

\overfullrule=9pt
%%%%%%%%%%%%%%%%%%%%%%%%%%%%%%%%%%%%%%%%%%%%%%%%%%%%%%%%%%%%%%%%%%
%                         THE TOP MATTER                         %
%%%%%%%%%%%%%%%%%%%%%%%%%%%%%%%%%%%%%%%%%%%%%%%%%%%%%%%%%%%%%%%%%%
%
\title{Mathematical and numerical methods for understanding immune cell motion during wound healing}
\thanks{This work was supported by grants H2020-955576, APVV-23-0186, VEGA 1/0249/24.}
\runningtitle{Mathematical and numerical methods for understanding immune cell motion}
\author{Giulia Lupi}
\address{Department of Mathematics, Slovak University of Technology, Radlinsk\'eho 11, 81005 Bratislava, Slovakia;} 
\author{Seol Ah Park}
\sameaddress{1}
\author{Martin Ambroz}
\sameaddress{1}
\author{Resul Ozbilgic}
\address{Laboratory of Pathogens and Host Immunity (LPHI), CNRS/Université de Montpellier, Place Eugène Bataillon, 34095 Montpellier cedex 5, France;} 
\author{Mai Nguyen Chi}
\sameaddress{2}
\author{Georges Lutfalla}
\sameaddress{2}
\author{Karol Mikula}
\sameaddress{1}
\runningauthors{G. Lupi, S.A. Park, M. Ambroz, R. Ozbilgic, M. Nguyen Chi, G. Lutfalla, K. Mikula}
\begin{abstract}
In this paper, we propose a new workflow to analyze macrophage motion during wound healing. These immune cells are attracted to the wound after an injury and they move showing both directional and random motion. Thus, first, we smooth the trajectories and we separate the random from the directional parts of the motion. The smoothing model is based on curve evolution where the curve motion is influenced by the smoothing term and the attracting term. Once we obtain the random sub-trajectories, we analyze them using the mean squared displacement to characterize the type of diffusion. Finally, we compute the velocities on the smoothed trajectories and use them as sparse samples to reconstruct the wound attractant field. To do that, we consider a minimization problem for the vector components and lengths, which leads to solving the Laplace equation with Dirichlet conditions for the sparse samples and zero Neumann boundary conditions on the domain boundary. 
\end{abstract}

\maketitle

\section{Introduction}

The movement of immune cells, particularly macrophages, and neutrophils, has been intensively studied \cite{barros2017live,friedl2008interstitial,li2012live,modes2010matrix} and their movement has been characterized in different tissues \cite{friedl2008interstitial,li2012live} and structures \cite{modes2010matrix}. These leukocytes are key players in the innate immune response; they are among the first cells to respond to tissue damage and infection, participating in the killing of invading pathogens, clearing of cellular debris, and tissue repair. 
 After an injury, danger signals are released from the wound and attract immune cells. Cues that can influence the navigation of macrophages are numerous and include gradients of chemoattractants that cause macrophage chemotaxis, inflammatory mediators, or changes in mechanical force \cite{ZHAO2009674}.
 Finding the wound attractant field is an important biological question, that holds significant clinical importance. Indeed, chronic inflammatory diseases and other pathologies are associated with an inappropriate immune response, and finding how wound chemoattractants propagate and control macrophage movement may help to develop new strategies to prevent an excessive or unbalanced immune response. 
Danger signals are thought to be released from the wound margin and include the so-called DAMPs (damage-associated molecular patterns). These were initially identified \textit{in vitro} and the actual contribution of some of them (such as Reactive oxygen species (ROS) and calcium) has also been studied \textit{in vivo} \cite{Tamara}. However, the understanding of the main contributors is still unclear. To address this problem for neutrophils migration, the authors in \cite{kadirkamanathan2012neutrophil} reconstructed the chemoattractant field that drives chemotaxis of this specific type of immune cells during wound healing using the zebrafish. While macrophages behave differently than neutrophils at the wound \cite{Ellett}, their chemoattractant field has not been determined. 

In this paper, we propose a new workflow to analyze macrophage trajectories and use the information extracted from the trajectories to reconstruct the wound attractant field. We exploit the optical clarity and genetic amenability of the zebrafish larva and, thanks to the macrophage reporter line, we follow the \textit{in vivo} behavior of macrophages in response to tail fin injury. 
The segmentation-based tracking approach developed in \cite{Sora} allowed us to obtain 2D$+$ time macrophage trajectories from microscopy videos. Moving toward a wound, macrophages exhibit both directional and random motion. We assume that the wound attractant field causes the directional part, while the random parts are parts of trajectories where the attractant field does not have a direct influence. Therefore, we first smooth the trajectories and remove the random parts. 
The mathematical model we propose to smooth the trajectories is based on evolving curve approach in the Lagrangian formulation. The smoothing process is influenced by three terms: trajectory \textit{smoothing} term depending on its curvature, term \textit{attracting} trajectories to the original one, and \textit{tangential velocity}. 
%The first term smooths the trajectory, while the attracting term keeps it close to the original trajectory. 
The tangential velocity does not change the shape of the evolving curve, but it redistributes the points on the curve, improving the stability of the numerical computations. It can be defined in different ways, for example by imposing a uniform redistribution of points, by conserving the initial distribution of points, or a redistribution that depends on the curvature \cite{hou1994removing,kimura1997numerical,sevcovic2001evolution,mikula2004direct,vsevvcovivc2011evolution}. In this paper, we use the asymptotically uniform redistribution of points; see \cite{mikula2008simple,mikula2004direct} for more details.

In the proposed trajectory smoothing model, the curve evolution is driven by the following partial differential equation
\begin{equation}
 \frac{\partial \textbf{x}}{\partial t}=-\delta k\textbf{N}+\lambda[(\textbf{x}_0-\textbf{x})\cdot \textbf{N}]\textbf{N}+\alpha \textbf{T}, 
 \label{E: Equation, Def}
\end{equation}
where $\textbf{x}$ is a point on the evolving curve, $-\delta k\textbf{N}$ is the smoothing term depending on the local curvature $k$, $\lambda[(\textbf{x}_0-\textbf{x})\cdot \textbf{N}]\textbf{N}$ is the attracting term, $\alpha$ is the tangential velocity, $\textbf{T}$ the unit tangent vector, $\textbf{N}$ the unit normal vector, and $\lambda$ and $\delta$ are the model parameters. 
We can rewrite equation (\ref{E: Equation, Def}) into the form
\begin{equation}
 \frac{\partial \textbf{x}}{\partial t}=\beta \textbf{N}+\alpha \textbf{T},
 \label{E: General Equation, Def}
\end{equation}
where $\beta$ is the velocity in the normal direction. 

Equation (\ref{E: General Equation, Def}) is a general equation for curve evolution models: the normal velocity $\beta$ can be defined in different ways and drives the evolution of the curve, depending on the purpose of the application. In \cite{mikula2021automated}, the authors developed a segmentation model based on curve evolution in the Lagrangian formulation. They defined the normal velocity taking into account the attracting vector field given by the image gradient. 
In \cite{ambroz2019numerical}, the authors developed a wildland fire propagation model based on surface curve evolution and empirical principles of fire propagation. In particular, they considered the evolution of a 3D surface curve and defined the normal velocity of the projected plane curve $\beta$ taking into account parameters such as wind speed and direction and topographic slope. In \cite{Lupi2022MacrophagesTS}, we proposed a model for the smoothing of macrophage trajectories in the Lagrangian formulation. The normal velocity is composed of the smoothing term and the attracting term so that the result of the smoothing keeps the trajectory close to the original one. To define the attracting term we identified, for each point on the evolving curve, the closest point on the original trajectory. The attracting term was then defined as the normal component of the vector connecting these two points.

In this paper, we propose a new way of defining the attracting term by considering the evolution of the lengths of the original segments. For each point on the evolving curve, we find to which of the original segments it belongs and compute its corresponding point on the original curve. This can not be done directly because the tangential velocity moves the points along the curve to redistribute them asymptotically uniformly. So we follow the evolution of the original points and find where they would have been if there had been no tangential redistribution of points. By this new approach, the evolving curve is kept close to the original in a consistent and more accurate way compared to \cite{Lupi2022MacrophagesTS}.

In \cite{Lupi2022MacrophagesTS}, we used the Hausdorff distance between two evolving curves as a stopping criterion. We stopped the curve evolution steps (also referred to as iterations within the text) when this distance was less than a given tolerance; thus, the stopping criterion considered the numerical steady state of the evolution. In this paper, we present a new way of defining the stopping criterion. Since we want to smooth the curve and keep only the directional parts of the motion, the new approach considers the self-intersections of the curves. In fact, where the curves self-intersect, the cells do not have a clear directionality but are moving back and forth. Therefore, we use the condition of having no further self-intersections as the stopping criterion. To detect the self-intersections of the curves, we developed a self-intersection detection algorithm based on topological changes detection algorithm from \cite{ambroz2019numerical}: the main idea is to construct a background (pixel) grid and check whether two non-consecutive grid points belong to the same pixel. If so, the self-intersection is detected; the details of the algorithm will be explained in Section \ref{SS: Detect self intersections}.

Since the main idea is to smooth the curve while keeping the directional parts close to the original, we considered an adaptive choice of the $\delta$ and $\lambda$ parameters based on the curve self-intersections. Equation (\ref{E: Equation, Def}) then becomes
\begin{equation}
 \frac{\partial \textbf{x}}{\partial t}=-\delta(\textbf{x},t) k\textbf{N}+\lambda(\textbf{x},t)[(\textbf{x}_0-\textbf{x})\cdot \textbf{N}]\textbf{N}+\alpha \textbf{T}, 
 \label{E: Equation adaptive choice of parameters, Def}
\end{equation}
where $\delta(\textbf{x},t)$ and $\lambda(\textbf{x},t)$ now depend on the spatial and temporal coordinates.

The ability to distinguish between directional and random motion is crucial to accurately describe the different modes of motion within a trajectory. A particular method for separating a trajectory into subregions is described. e.g., in \cite{Arcizet2008Temporal}. The authors developed a rolling window algorithm to separate active and passive parts of cell motion, taking into account the local mean squared displacement (MSD) and directional persistence in a local window. The mean squared displacement is a common tool for analyzing random motion; it distinguishes the type of diffusive motion followed by the walker. The motion can be diffusive, such as Brownian motion, superdiffusive, such as Levy-flights (random walks where small steps are interspersed with long jumps), and subdiffusive. Subdiffusive motion is usually observed in crowded environments where the walker is not free to move in any direction but is constrained by the environment \cite{kubala2021diffusion, viswanathan2011physics, sokolov2012models}.
Several studies have analyzed cell movement by examining the MSD \cite{daumas2003confined,de2007cellular,dieterich2008anomalous,harris2012generalized,huda2018levy,kusumi1993confined,lee1991direct}. In \cite{huda2018levy}, the authors analyzed the behavior of metastatic and non-metastatic cancer cells migrating on linear microtracks. The large-scale statistical analysis using the averaged MSD showed that metastatic cancer cells follow a different movement strategy than non-metastatic ones. In \cite{harris2012generalized}, a similar analysis showed that a generalized Levy walk describes T-cell motility in the brains of infected mice.

In our approach, the smoothing model allows us to separate the directional parts from the random parts. We define two different ways of finding the random parts. The first takes into account the evolution of the lengths of the original segments. During the evolution, a segment may disappear. This usually happens in regions where the curvature is high and the segments are short, which appears in random regions. If the segment disappears, we consider it to be a random segment; consecutive disappeared segments form a random sub-trajectory. The other way to define the random parts is to consider the self-intersections. In fact, the random parts of the trajectories are characterized by back-and-forth motion without a clear directionality, which causes the curve self-intersection. A random sub-trajectory is then defined as the segments belonging to a self-intersecting part. Once we obtained the random parts, to characterize the type of diffusion followed by the cells in regions where the vector field does not have a direct influence on the motion, we used the MSD.

Finally, we computed the velocities on the smoothed curves and used them as sparse samples to reconstruct the wound attractant field. 
Many different computational contexts require the reconstruction of vector fields from sparse samples: the applications of such a process include fluid dynamics visualization, texture synthesis, non-photorealistic rendering, optical flow fields, and map registration \cite{DesignTangentVF, VectorfieldReconstruction, article}. In \cite{VectorfieldReconstruction}, the authors proposed local piecewise polynomial approximations using least squares methods. In \cite{Mussa-Ivaldi}, a different approach was taken, where linearly independent vector fields, termed 'basis fields,' were combined to approximate the vector field, maintaining invariance under coordinate transformations. In \cite{DesignTangentVF}, the authors introduced a model for reconstructing smooth tangent vector fields. Using sparse samples as user-specified constraints, they applied a coordinate-free approach and employed the discrete Laplace operator to obtain a smooth field. Laplace interpolation has also been applied to image processing tasks like image inpainting and compression \cite{doi:10.1137/080716396, Schnlieb2015PartialDE}. In \cite{article}, the Locally Affine Globally Laplace (LAGL) transformation model was developed, involving the design of affine transformations between maps, followed by solving the Laplace equation to smoothly reconstruct intermediate points.

 In this paper, we solve a minimization problem for the two vector components and the vector length considering Dirichlet conditions for the sparse samples and zero Neumann boundary conditions at the domain boundary. Solving the minimization problem is equivalent to solving the Laplace equation with prescribed boundary conditions. To prove the existence and uniqueness of a weak solution for the considered problem, we use the theorem proved in \cite{LupiVF}. For the discretization, we place the unknowns at the vertices of the pixels and use them as grid nodes for the finite difference scheme. This results in a smooth vector field where the information from the sparse samples is interpolated in regions where the vectors influence each other and extrapolated in regions where they do not. Considering only the minimization problem for the two vector components leads to a smooth vector field in which the lengths of the vectors decrease with distance from the Dirichlet conditions. This is due to the property of the Laplace equation of averaging neighboring values. Since we did not want to add any new information, but to reconstruct the vector field from the information given by the sparse samples, we also considered the minimization problem for the vector lengths. This leads to a smooth vector field, where the direction of the vectors changes smoothly from one Dirichlet condition to the other, and where the vector lengths are interpolated/extrapolated from the Dirichlet conditions.
 
 The paper is organized as follows. In Section \ref{S: Mathematical model}, we will introduce the mathematical model for smoothing trajectories, define the attracting term, and consider the discretization of the model. Then we will present the self-intersections detection algorithm and illustrate how to choose the parameters adaptively. We will conclude Section \ref{S: Mathematical model} by presenting results on real macrophage trajectories. In Section \ref{S: Analysis of random parts}, we will illustrate how to detect the random parts and analyze them using the mean squared displacement. Finally, in Section \ref{S: Vector Field}, we will define the model we propose for vector field reconstruction from sparse samples and illustrate its discretization. The results will be shown for 3 different datasets of macrophages moving during wound healing.
\subsection{Data acquisition}
\label{SS: Data acquisition}
We considered $3$ different datasets; in all of them, zebrafish larvae were utilized to study macrophage motion during wound healing. 
At 3 days post-fertilization, transgenic reporter larvae, Tg(mfap4:mCherry-F), were employed. These larvae express a farnesylated mCherry fluorescent protein under the control of the macrophage-specific promoter mfap4, enabling macrophages to display mCherry at their membranes. Macrophage cells are followed in the red channel. All larvae were wounded at the caudal fin fold region and subsequently imaged using a Spinning Disk Confocal microscope.

Migrating macrophages were imaged starting 30 minutes to 1-hour post-amputation, over a duration of 13 hours, with a time interval of 2.5 minutes between frames and a z-step of 1 µm. For all datasets, the pixel size was 0.319489 µm. For numerical experiments, we used 2D+time projection images, where the three-dimensional (3D) microscopy data were flattened into a 2D plane by applying a maximum intensity projection. In this method, the maximum intensity value along the z-axis was selected for each pixel.

\section{Trajectories smoothing}
\label{S: Mathematical model}

Let $\Gamma$ be an open plane curve with fixed endpoints
\begin{equation} \label{E: Gamma, Def}  
\begin{split}
    \Gamma:[0& ,1]\rightarrow\mathbb{R}^2,\\
    & u\mapsto\textbf{x}(u),
\end{split}  
\end{equation}
where $\textbf{x}(u)=(x_1(u),x_2(u))$ is a point of the curve $\Gamma$. We consider the time evolution of a point $\textbf{x}(u,t)\in\Gamma_t$ of the curve $\Gamma$ driven by the equation (\ref{E: Equation adaptive choice of parameters, Def}), where $t$ is time, $k$ is the curvature, $\textbf{x}_0$ is a point on the original curve, $\textbf{T}$ is the unit tangent vector, $\textbf{N}$ is the unit normal vector, $\alpha$ is the tangential velocity and $\lambda(\textbf{x},t)$ and $\delta(\textbf{x},t)$ are positive parameters.\\
From now on we will indicate by
\begin{equation}
    w=(\textbf{x}_0-\textbf{x})\cdot \textbf{N}
    \label{E:w,Def}
\end{equation}
the attracting term, and by
\begin{equation}
    \beta=-\delta k+\lambda w,
\end{equation}
the velocity in the normal direction.\\
The tangential velocity is defined in such a way that the points are asymptotically uniformly distributed: the formula is obtained by considering the ratio between the local and the global length and letting it tend to $1$; see \cite{mikula2008simple,mikula2004direct} for details on the derivation and discretization. On one hand, the redistribution of points along the curve improves the stability of the numerical computations; on the other hand, it makes the definition of attracting term more complicated. In fact, the points also move in the tangential direction. If we then define the \textit{corresponding point} of the point $\textbf{x}(u)$ simply as $\textbf{x}_0(u)$, we can obtain that $\textbf{x}(u)$ is attracted to a completely different part of the trajectory. Therefore, to define the attracting term independently of the tangential velocity, we follow the evolution of the original segments as if there was no tangential redistribution. Two points are defined as corresponding if they belong to the same evolving segment.

\subsection{Attracting term}
\label{SS: Lagrangian approach to attracting term}
In order to understand the main idea of the model described below, we must first note how trajectories are found in videos. The original trajectories are piecewise linear segments where the time between two linear segment endpoints is $1$ time unit, and a point of the trajectory corresponds to the position of a cell in a time frame. First, for the numerical discretization, we add grid points to the linear segments: refining the grid is a necessary step for the numerical computations to improve the accuracy. To do this, we set a value $\bar{h}$ and add grid points inside the linear segments in a distance of $\bar{h}$, creating an initial, nearly uniform distribution of grid points.
Then we distinguish between a \textit{segment}, which connects two of the original points of the trajectory, and an \textit{element}, which connects two grid points of the evolving curve. Our interest is in the \textit{segments}: the model we propose computes the new length of the segments on the evolving curve as if there was no tangential velocity, and it finds the endpoints of the segment on the evolving curve. The model to find the new lengths of the segments has already been proposed in \cite{Lupi2022MacrophagesTS}, but was not used in the definition of attracting term; we recall here the derivation and later its discretization for completeness. Once we have the segments on the evolving curve and the original ones, we define the attracting term as the normal component of the vector connecting the points in the corresponding segments.\\
Let us consider the unit arc length parameterization $s$ of the curve $\Gamma$ so that it holds $ds=\lvert\textbf{x}_u\rvert du$. The local length $g$ is defined as
\begin{equation}\label{E: g, Def}
    g=\lvert\textbf{x}_u\rvert,
\end{equation}
while the global length is
\begin{equation}
    L=\int_ {\Gamma} ds=\int_{0}^{1}g du=\int_{0}^{1} \lvert\textbf{x}_u\rvert du.
\end{equation}
We consider the time derivative $L_t$ of the global length; using the Frenet-Serret formulas we obtain
\begin{equation}
    g_t= \lvert\textbf{x}_u\rvert_t=gk\beta+g\alpha_s=gk\beta+\alpha_u.
\end{equation}
Then
\begin{equation}
\label{E: L derivative, Def}
    L_t=\int_{0}^{1}g_t~du=\int_{0}^{1}(gk\beta+\alpha_u)~ du=\int_{0}^{1}gk\beta~ du+\alpha(1)-\alpha(0)=\int_{\Gamma}k\beta~ ds,
\end{equation}
where $\alpha(1)=\alpha(0)=0$ because the endpoints of the curve are fixed.\\
We can consider the same equation (\ref{E: L derivative, Def}) for the evolution of the lengths of the segments of the original trajectory. For the $j-$th segment, we obtain
\begin{equation}
    (L_t)_j=\int_{u_{j-1}}^{u_j}(gk\beta+\alpha_u)~ du=\int_{\textbf{x}(u_{j-1})}^{\textbf{x}(u_j)}k\beta~ ds +\alpha(u_{j})-\alpha(u_{j-1}),
\end{equation}
In this case, in general, it holds $\alpha(u_j)\neq 0$, $\alpha(u_{j-1})\neq 0$. But, for our purpose, it is reasonable to consider the model for the evolution of the length of the segments with $\alpha=0$. We obtain
\begin{equation}\label{E: L_j derivative final formula, Def}
    (L_t)_j=\int_{\textbf{x}(u_{j-1})}^{\textbf{x}(u_j)}k\beta~ ds,
\end{equation}
Note that the tangential velocity does not influence the shape of the evolving curve, it only redistributes the points on the evolving curve. It follows that, in a continuous setting, the shape of the evolved curve will be the same with or without the tangential velocity; therefore, considering $\alpha=0$ allows us to find what is the real length of the segments after smoothing. The method consists of the following steps.
\begin{enumerate}
    \item Find the new length of the segment using equation (\ref{E: L_j derivative final formula, Def}).
    \item Find the endpoints of the segments on the smoothed curve.
    \item Define $N$ the number of grid points inside the new segment and consider a uniform distribution of $N$ points inside the original segment.
    \item Define the attracting term as
    \begin{equation}
        (\textbf{x}^0-\textbf{x})(u)=\textbf{x}^0(u)-\textbf{x}(u),
    \end{equation}
    where $\textbf{x}^0(u)$ and $\textbf{x}(u)$ are grid points in the same segment on the original and smoothed curve, respectively.
\end{enumerate}
 A segment can disappear during the evolution; if this happens, we consider that the segment degenerates to a point, i.e. the two endpoints merge. In this case, the point is attracted to the center of the original segment. If the same happens for consecutive segments, we generalize the formula and consider the center of mass of the segments and the vector 
\begin{equation}
    (\textbf{x}^0-\textbf{x})(u)=\textbf{x}^{0,M}-\textbf{x}(u),
\end{equation}
where $\textbf{x}^{0,M}$ is the center of mass. Since the tangential velocity does not affect the shape of the evolving curve, we consider only the normal component of the attracting vector. The details of the numerical discretization are explained in the next section.

\subsection{Numerical algorithm}
\label{S: Numerical algorithm}
The discretization of Equation (\ref{E: Equation adaptive choice of parameters, Def}) is done as in \cite{Lupi2022MacrophagesTS}. We recall here the main steps for completeness. \\ 
Let $g=\lvert\textbf{x}_u\rvert=\sqrt{(\frac{\partial x_1}{\partial u})^2+(\frac{\partial x_2}{\partial u})^2}$. If we denote by $s$ the unit arc length parametrization of the curve $\Gamma$, then $ds=\lvert\textbf{x}_u\rvert du=g du$ and $du=\frac{1}{g}ds$. The unit tangent vector $\textbf{T}$ is defined as $\textbf{T}=\frac{\partial \textbf{x}}{\partial s}=\textbf{x}_s$ while the unit normal vector $\textbf{N}$ is $\textbf{N}=\textbf{x}_{s}^{\perp}$ such that $\textbf{T}\wedge\textbf{N}=-1$. If $\textbf{T}=(\frac{x_1}{\partial s},\frac{x_2}{\partial s})$ then $\textbf{N}=(\frac{x_2}{\partial s},-\frac{x_1}{\partial s})$. From the Frenet-Serret formulas, we have $\textbf{T}_s=-k\textbf{N}$ and $\textbf{N}_s=k\textbf{T}$, where $k$ is the curvature. It follows that $-k\textbf{N}=\textbf{T}_s=(\textbf{x}_s)_s=\textbf{x}_{ss}$.\\
We can rewrite (\ref{E: Equation adaptive choice of parameters, Def}) into the form of the so-called intrinsic partial differential equation
\begin{equation}
 \textbf{x}_t=\delta\textbf{x}_{ss}+\alpha\textbf{x}_s+\lambda w\textbf{x}_{s}^{\perp}
 \label{E: Intrinsic form of PDE, Def}
\end{equation}
which is suitable for numerical discretization. Since $\textbf{x}=(x_1,x_2)$, (\ref{E: Intrinsic form of PDE, Def}) represents a system of two partial differential equations for the two components of the position vector $\textbf{x}$. 

Let us consider the following form of the intrinsic PDE (\ref{E: Intrinsic form of PDE, Def})
\begin{equation}
 \textbf{x}_t-\alpha\textbf{x}_s=\delta\textbf{x}_{ss}+\lambda w\textbf{x}_{s}^{\perp},
\end{equation}
where $w$ is given by (\ref{E:w,Def})  and $\alpha$ is defined in \cite{Lupi2022MacrophagesTS}.\\
Consider an open curve discretized into $n+2$ grid points, $i=0,...,n+1$. We fixed the two points $ \textbf{x}_0$ and $ \textbf{x}_{n+1}$ at the endpoints of the curve. Consequently, for these two points, the velocity is zero, i.e. $\alpha_0=\beta_0=\alpha_{n+1}=\beta_{n+1}=0$. We define the $i$-th element $\textbf{h}_i=\textbf{x}_i-\textbf{x}_{i-1}$ with length $h_i=\lvert \textbf{x}_i-\textbf{x}_{i-1} \rvert$.\\ 
For the spatial discretization, we use the flowing finite volume method \cite{mikula2008simple}. Consider the finite volume $\textbf{p}_i=[ \textbf{x}_{i-\frac{1}{2}}, \textbf{x}_{i+\frac{1}{2}}]$ where $ \textbf{x}_{i+\frac{1}{2}}$ represents the middle point between $ \textbf{x}_i$ and $ \textbf{x}_{i+1}$. 
Integrating over the finite volume $\textbf{p}_i$, using the Newton-Leibniz formula and considering $\alpha_i$, $w_i$, $\delta_i$, and $\lambda_i$ constant over the finite volume we obtain 
\begin{equation}
\begin{split}
 &\frac{h_i+h_{i+1}}{2}(\textbf{x}_i)_t+\frac{\alpha_i}{2}(\textbf{x}_i-\textbf{x}_{i+1})-\frac{\alpha_i}{2}(\textbf{x}_i-\textbf{x}_{i-1})=\\
 &=\delta_i(\frac{\textbf{x}_{i+1}-\textbf{x}_i}{h_{i+1}}-\frac{\textbf{x}_{i}-\textbf{x}_{i-1}}{h_{i}})+\lambda_i w_i(\frac{\textbf{x}_{i+1}-\textbf{x}_{i-1}}{2})^{\perp}.
 \end{split}
\end{equation}
For the time discretization, we use the semi-implicit time scheme, which is implicit in the intrinsic diffusion term and uses the inflow-implicit/outflow-explicit strategy for the intrinsic advection term \cite{mikula2011inflow}. Let $m$ be the time step index and $\tau$ the length of the discrete time step. The basic idea of the inflow-implicit/outflow-explicit strategy is that we treat outflow from a cell explicitly while inflow implicitly \cite{mikula2011inflow}. Therefore, we define
\begin{equation}
\begin{split}
 b_{i-\frac{1}{2}}^{in}=max(-\alpha_{i}^{{m}},0),b_{i-\frac{1}{2}}^{out}=min(-\alpha_{i}^{{m}},0),\\
 b_{i+\frac{1}{2}}^{in}=max(\alpha_{i}^{{m}},0),b_{i+\frac{1}{2}}^{out}=min(\alpha_{i}^{{m}},0).
 \end{split}
\end{equation}
Approximating the time derivative by finite difference, we obtain the fully discrete scheme
\begin{equation}
 \begin{split}
  &\textbf{x}_{i-1}^{m+1}(-\frac{\delta_i}{h_{i}^{m}}-\frac{b_{i-\frac{1}{2}}^{in}}{2})+\textbf{x}_{i+1}^{m+1}(-\frac{\delta_i}{h_{i+1}^{m}}-\frac{b_{i+\frac{1}{2}}^{in}}{2})+\\
  &\textbf{x}_{i}^{m+1}(\frac{h_{i+1}^{m}+h_{i}^{m}}{2\tau}+\frac{\delta_i}{h_{i}^{m}}+\frac{\delta_i}{h_{i+1}^{m}}+\frac{b_{i-\frac{1}{2}}^{in}}{2}+\frac{b_{i+\frac{1}{2}}^{in}}{2})=\textbf{x}_{i}^{m}\frac{h_{i+1}^{m}+h_{i}^{m}}{2\tau}\\
  &-\frac{b_{i+\frac{1}{2}}^{out}}{2}(\textbf{x}_{i}^{m}-\textbf{x}_{i+1}^{m})-\frac{b_{i-\frac{1}{2}}^{out}}{2}(\textbf{x}_{i}^{m}-\textbf{x}_{i-1}^{m})+\lambda_i w_{i}^{m}(\frac{\textbf{x}_{i+1}^{m}-\textbf{x}_{i-1}^{m}}{2})^{\perp},
 \end{split}
 \label{E: Discr final, Discr}
\end{equation}
for $i=1,...,n$ where $n$ is the number of unknown grid points.\\ 
The system (\ref{E: Discr final, Discr}) is represented by a strictly diagonally dominant matrix, then it is always solvable by the classical Thomas algorithm without any restriction on the time step $\tau$.\\
For complex trajectories, the evolving curve may have singularities. For those parts of the trajectories, we then replace the second-order inflow-implicit/ outflow-explicit strategy to discretize the advection term with the first-order implicit upwind scheme \cite{ambroz2019numerical}. We consider the singularity when the angle between two consecutive elements is less than $120^\circ$, i.e. $\measuredangle \textbf{h}_i\textbf{h}_{i+1}<120^\circ$. In the upwind scheme we consider only $b_{i-\frac{1}{2}}^{in}$ and $b_{i+\frac{1}{2}}^{in}$, therefore we obtain the fully discrete scheme
\begin{equation}
 \begin{split}
  &\textbf{x}_{i-1}^{m+1}(-\frac{\delta_i}{h_{i}^{m}}-b_{i-\frac{1}{2}}^{in})+\textbf{x}_{i+1}^{m+1}(-\frac{\delta_i}{h_{i+1}^{m}}-b_{i+\frac{1}{2}}^{in})+\\
  &\textbf{x}_{i}^{m+1}(\frac{h_{i+1}^{m}+h_{i}^{m}}{2\tau}+\frac{\delta_i}{h_{i}^{m}}+\frac{\delta_i}{h_{i+1}^{m}}+b_{i-\frac{1}{2}}^{in}+b_{i+\frac{1}{2}}^{in})=\textbf{x}_{i}^{m}\frac{h_{i+1}^{m}+h_{i}^{m}}{2\tau}\\&+\lambda_i w_{i}^{m}(\frac{\textbf{x}_{i+1}^{m}-\textbf{x}_{i-1}^{m}}{2})^{\perp},
 \end{split}
 \label{E: Upwind}
\end{equation}
which is also represented by a strictly diagonally dominant matrix, and therefore solvable by the Thomas algorithm. The use of (\ref{E: Upwind}) improves the robustness of the scheme with respect to singularities. In the next Section, we will describe how to discretize the attracting term $w_i$.

\subsubsection{Discretization of attracting term}
\label{SS:DiscretizationAT}
\begin{figure}[htbp]
  %\centering
     \includegraphics[scale=0.52]{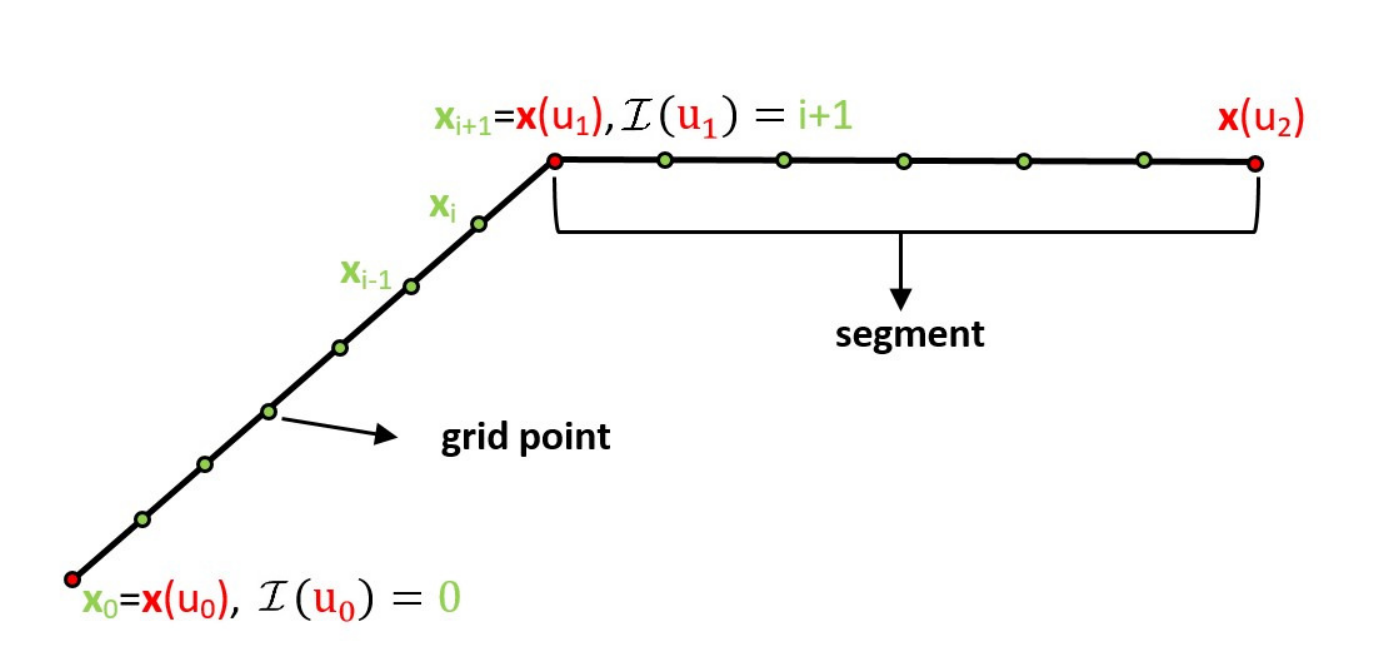}
  \caption{Visualization of curve discretization. The green dots represent the grid points and the red dots represent the endpoints of the segments.}
  \label{I_expl}
\end{figure}

From (\ref{E: g, Def}) we have
\begin{equation}
    g=\lvert\textbf{x}_u\rvert=\sqrt{\left(\frac{\partial x_1}{\partial u}\right)^2+\left(\frac{\partial x_2}{\partial u}\right)^2}
\end{equation}
To discretize (\ref{E: L_j derivative final formula, Def}), we approximate the time derivative by the finite difference and consider the formulas for $k,~\beta$ given by \cite{ambroz2020semi,mikula2021automated,mikula2008simple}
\begin{equation}
 \begin{split}
  &k_{i}^{m}=sgn(\textbf{h}_{i-1}^{m}\wedge \textbf{h}_{i+1}^{m})\frac{1}{2h_{i}^{m}}arccos(\frac{\textbf{h}_{i-1}^{m}\cdot \textbf{h}_{i+1}^{m}}{h_{i-1}^{m} h_{i+1}^{m}}),\\
  &\beta_{i}^{m}=-\delta k_{i}^{{m}}+\lambda  w_{i}^{{m}},\\
  &L^m=\sum_{l=1}^{n+1}h_{l}^{m},
  \end{split}
\end{equation}
where $\textbf{h}_{i}^{m}=\textbf{x}_{i}^{m}-\textbf{x}_{i-1}^{m}$, $\lvert\textbf{h}_{i}^{m}\rvert=h_{i}^{m}$ and $\textbf{h}_{i-1}^{m}\wedge \textbf{h}_{i+1}^{m}$ is the wedge product, i.e. the determinant of the matrix with columns $\textbf{h}_{i-1}^{m}$ and $\textbf{h}_{i+1}^{m}$. For the first and last elements, since we don't have the values of $h_{i-1}^m$, respectively $h_{i+1}^m$, we set $k_{1}^{m}=k_{2}^{m}$ and $k_{n+1}^{m}=k_{n}^{m}$.
Note that, for $m=0$, $w_{i}^{0}=0$ for every $i=1,...,n$, therefore $\beta_{i}^{0}=-\delta k_{i}^{0}$. This is natural since in the first iteration the evolving curve and the original curve are the same: it follows that the attracting vector is $\textbf{0}$.\\
To relate the index $j$ of the segments to the index $i$ of the grid points we will indicate by $\mathcal{I}(u_{j-1})$, $\mathcal{I}(u_j)$ the corresponding $i$ index of the endpoint $\textbf{x}(u_{j-1})$, respectively $\textbf{x}(u_j)$, see Fig. \ref{I_expl}.
We obtain the fully discrete approximation of equation (\ref{E: L_j derivative final formula, Def}) for the new length of the $j$-th segment
\begin{equation}
L_{j}^{m+1}=L_{j}^{m}+\tau \sum_{i=\mathcal{I}(u_{j-1})+1}^{\mathcal{I}(u_j)}h_{i}^{m}k_{i}^{m}\beta_{i}^{m}.
\label{E: L_j, Discr}
\end{equation}
Experimentally, it holds $L_{j}^{m+1} \leq L_{j}^{m}$. From (\ref{E: L_j, Discr}) it is clear that if the segment is in a region with high curvature (like the random parts), its length will decrease faster and it may also disappear, $L_{j}^{m+1}=0$.\\
After obtaining the new length $L_{j}^{m+1}$ of the segment applying the formula in (\ref{E: L_j, Discr}), we need to compute the new $i$ index for $\mathcal{I}(u_{j-1})$ and $\mathcal{I}(u_{j})$. The straightforward approach would be, starting from the first segment $L_{1}^{m+1}$, to sum up the lengths of the elements $h_{i}^{m+1}$ until the sum is equal to $L_{j}^{m+1}$ and consider $\mathcal{I}(u_{j})$ as the last endpoint of the corresponding element $\textbf{h}_{i}^{m+1}$.\\
Consider $\sum_{j=1}^{M} L_{j}^{m+1}$, where $M$ is the number of original segments.
 A problem arises since in the numerical computations it holds
\begin{equation}
    \sum_{i=1}^{n+1} h_{i}^{m+1}\neq \sum_{j=1}^{M} L_{j}^{m+1}.
\end{equation}
It follows that we need to define a way to relate the two global lengths. Consider the ratio
\begin{equation}
 r_{j}^{m+1}=\frac{L_{j}^{m+1}}{\sum_{j=1}^{M} L_{j}^{m+1}},
\end{equation}
and define the discrete length of the segment as
\begin{equation}
 L^{d}_{j}=r_{j}^{m+1}\sum_{i=1}^{n+1} h_{i}^{m+1},
 \label{E: L^d_j}
\end{equation}
so that
\begin{equation}
    \sum_{j=1}^{M}L^{d}_{j}=\sum_{i=1}^{n+1} h_{i}^{m+1}.
\end{equation}
Having obtained the discrete lengths $L_{j}^{d}$, we can apply the straightforward approach. Let us illustrate the approach for $L_{1}^{d}$; the process is then repeated for all others $L_{j}^{d}$, $j=2,...,M$. We sum the lengths of the elements $h_i$ until $\sum_{i=1}^{k} h_i \geq L_{j}^{d}$, for some $k\in \{1,...,n+1\}$. If
\begin{equation}
    \sum_{i=1}^{k} h_i = L_{j}^{d},
    \label{E: sum equal L_1}
\end{equation}
then $\mathcal{I}(u_1)=k$. Otherwise, we find the point $\textbf{x}^*$ on the element $\textbf{h}_i$ such that equality (\ref{E: sum equal L_1}) holds. Then, we find which point between $\textbf{x}_k$ and $\textbf{x}_{k-1}$ is the closest point to $\textbf{x}^*$ and move it to $\textbf{x}^*$. Finally, we set either  $\mathcal{I}(u_1)=k$ or $\mathcal{I}(u_1)=k-1$ depending on which point was moved to $\textbf{x}^*$ and $L_{j}^{m+1}=L_{j}^{d}$. \\
Notice that we focused only on the index of the last endpoint of each segment $\mathcal{I}(u_j)$: that is because the last endpoint of the segment $j$ is the first endpoint of the segment $j+1$.\\

Let us now focus on the disappearing segments. We add the following conditions: 
 if $L_{j}^{m+1}<h_{mean}$, then $L_{j}^{m+1}=0$, where
\[
 h_{mean}=\frac{1}{n+1}\sum_{i=1}^{n+1} h_i.
\]
If we obtain $L_{j}^{m+1}=0$, it means the first and the last endpoints of that segment at time step index $m+1$ became the same point, then $\mathcal{I}(u_{j})=\mathcal{I}(u_{j-1})$. So we modify the algorithm for finding the new indexes $i$ for $\mathcal{I}(u_j)$ adding the following condition
\begin{itemize}
 \item if $r_{j}^{m+1}=0$ then $\mathcal{I}(u_j)=\mathcal{I}(u_{j-1})$.
\end{itemize}
This condition will cause that when we apply formula (\ref{E: L_j, Discr}) for every time step index $m+k$, with $k>1$, we will obtain
\[
 L_{j}^{m+k}=L_{j}^{m+1}=0.
\]
Once we have obtained the new length of the segments, from our algorithm we are also able to calculate how many grid points are inside a segment: since the algorithm gives us the first and the last end points of a segment, it is sufficient to consider $N_j=\mathcal{I}(u_{j})-\mathcal{I}(u_{j-1})+1$. Let us consider for every segment on the original curve, indexed by $j$, a uniform distribution of $N_j$ grid points on the segment. Therefore, we obtain a new parametrization of the original curve. Let us indicate the grid points in the new parametrization on the original curve with the same index $i$, i.e. $\textbf{x}^{0}_{i}$, $i=0,...,n+1$. Then, the attracting vector for the grid points is defined as 
\begin{equation}
    (\textbf{x}^0-\textbf{x})^{m}_{i}=\textbf{x}^{0}_i-\textbf{x}_{i}^{m}.
\end{equation}
If a segment disappears, let's say $L^{m}_{j}=0$, we consider it degenerated to a point, so $\mathcal{I}(u_{j-1})=\mathcal{I}(u_{j})$. In this case, we consider the middle point of the segment
\begin{equation}
    \textbf{x}^{0, M}=\frac{\textbf{x}^{0}_{\mathcal{I}(u_{j-1})}+\textbf{x}^{0}_{\mathcal{I}(u_{j})}}{2}.
\end{equation}
If more subsequent segments disappear, we generalize the formula calculating the center of mass for those segments. For the point where one or more segments degenerated to, we define the attracting vector as
\begin{equation}
     (\textbf{x}^0-\textbf{x})^{m}_{i}=\textbf{x}^{0, M}-\textbf{x}_{i}^{m}
\end{equation}
and $w_{i}^{m}$ is defined as 
\begin{equation}
 w_{i}^{m}=(\textbf{x}^{0}-\textbf{x})^{m}_{i}\cdot \textbf{N}_{i}^{m},
\end{equation}
where $\textbf{N}_{i}^{m}=(\frac{\textbf{x}_{i+1}^{m}-\textbf{x}_{i-1}^{m}}{h_{i}^{m}+h_{i+1}^{m}})^{\perp}$.

\subsection{Numerical experiments}
\label{S: Results}
\begin{figure}[htbp]
\centering
   \includegraphics[width=10cm]{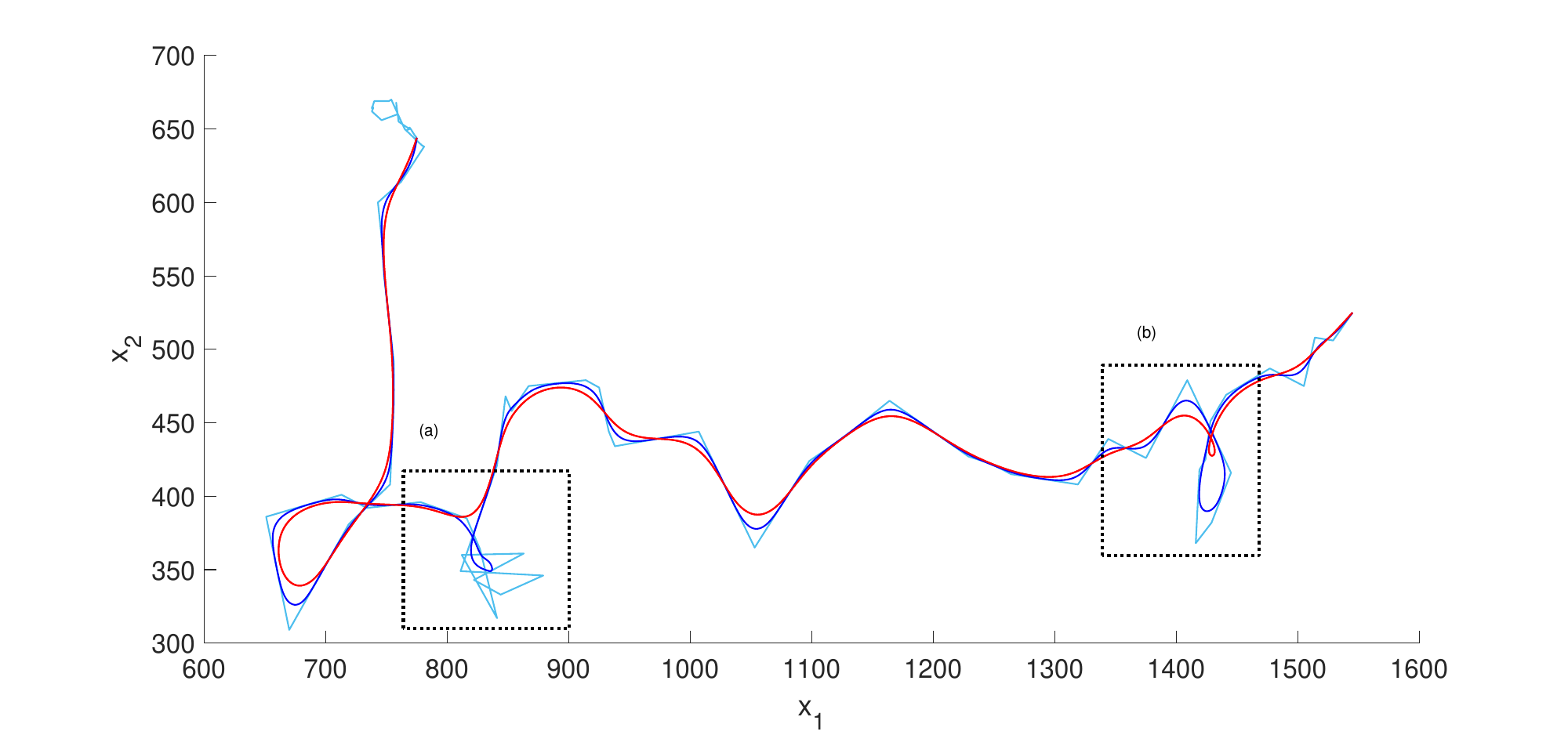} 
\subfigure[]{\includegraphics[width=10cm]{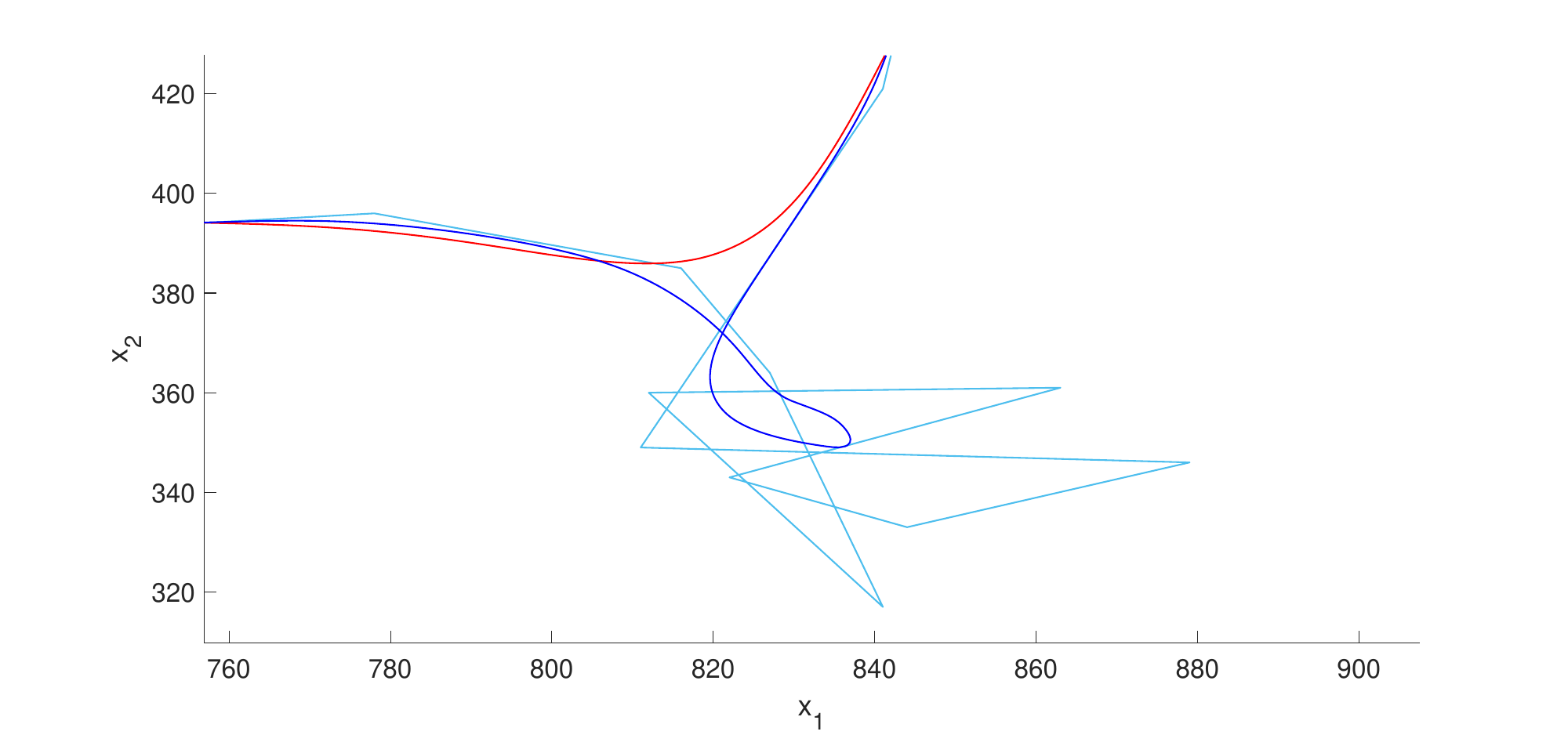}}
\subfigure[]{\includegraphics[width=10cm]{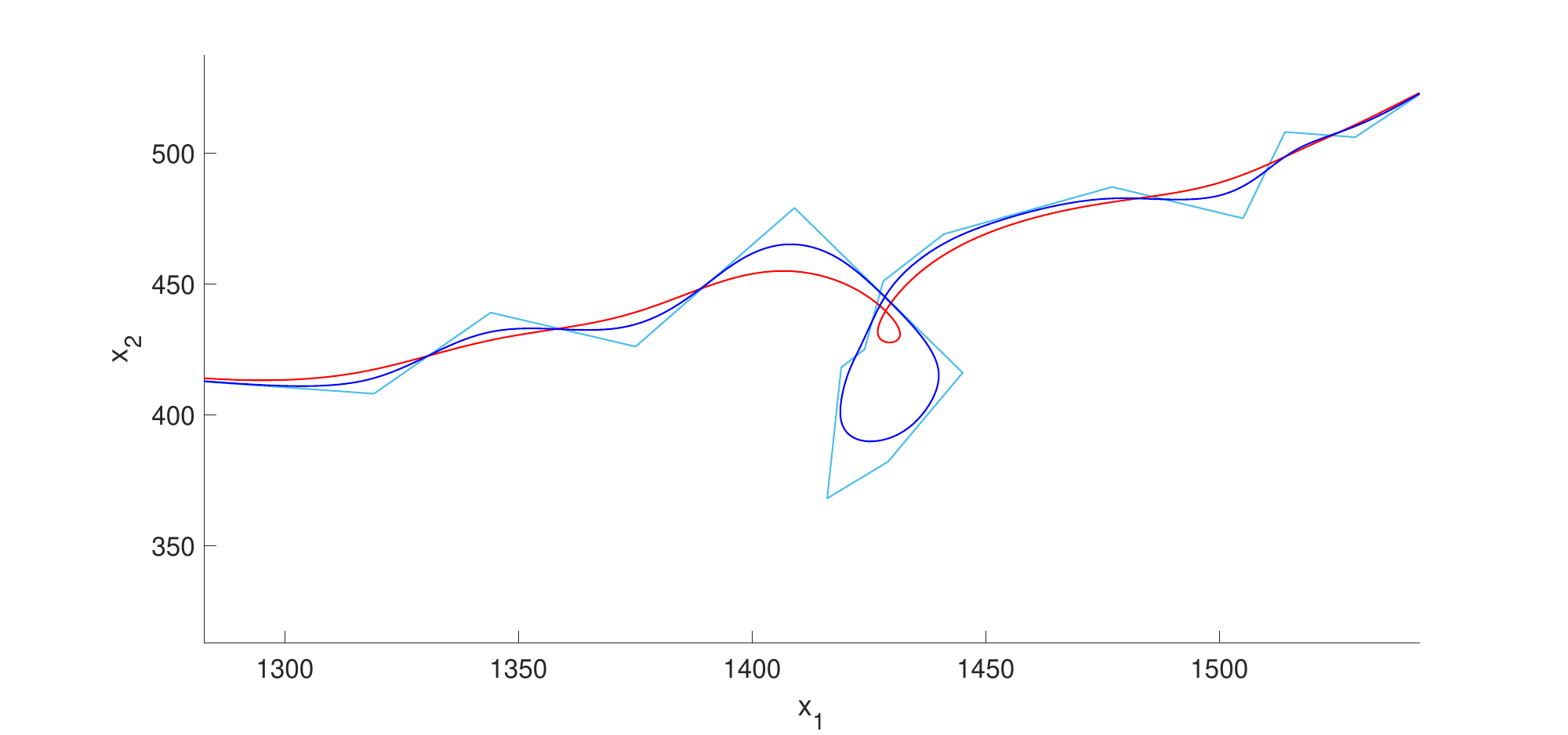}}
 \caption{Original trajectory (light blue line) and results of smoothing for $\lambda=10$ (red line) and $\lambda=50$ (blue line). The black boxes indicate the regions chosen for the zooms. Middle and bottom: zooms of the smoothing process.}
 \label{F:LambdaComparisonT2DS3}
\end{figure}

First, to see how the different parts of the normal velocity work, we studied the results for different values of the parameters $\lambda$ and $\delta$. In \cite{Lupi2022MacrophagesTS} we studied the results for simple curves to see the influence of the attracting term. The steady state keeps the curve closer to the original one than in the case where there is only the curvature term. Moreover, varying the values of $\lambda$ also changes the steady state; the higher the value of $\lambda$, the more the curve is kept closer to the original. Here we consider the new definition of the attracting term defined in Sections \ref{SS: Lagrangian approach to attracting term} and \ref{SS:DiscretizationAT} and different values of $\lambda$, namely $\lambda=10$ and $\lambda=50$. The other parameters were chosen as follows: $\delta=0.003$, $\tau=0.00001$, and $\omega=0.8$. Here, $\omega$ is the parameter that controls the speed of convergence to uniform redistribution of points in the tangential velocity $\alpha$; see \cite{mikula2004direct,mikula2008simple} for details. Fig. \ref{F:LambdaComparisonT2DS3} shows the results for different values of $\lambda$. The stopping criterion in this case considers the Hausdorff distance between two evolving curves and stops the iterations when the distance is less than the tolerance $\epsilon=0.00065$ as in \cite{Lupi2022MacrophagesTS}. As one can see, in both cases the curve is smoothed everywhere, both in the random parts and in the directional parts of the motion, but the loops are not completely removed. 
We noticed that where the curve intersects itself, the cell doesn't have a clear directionality, but moves back and forth. Since we only want to keep the directional parts of the motion, we developed a new stopping criterion based on the curve self-intersections. The algorithm from \cite{ambroz2019numerical} for fast ($O(n)$) detection of topological changes was adjusted to detect self-intersections; we will recall it in the next section. 

\subsubsection{Detection of self-intersections}
\label{SS: Detect self intersections}
The main idea of self-intersection detection is to construct a background grid and check which points of the trajectory belong to the same pixel. If the points are not consecutive, then they are considered to be points where the curve is self-intersecting. To do this, we consider the initial distance $\bar{h}$ of the grid points and construct 2 background grids with pixel size $2\bar{h}$ shifted by $\bar{h}$. Recall that in the beginning, we add grid points to the original segments at intervals of $\bar{h}$, obtaining an almost uniform distribution of points (see Section \ref{SS: Lagrangian approach to attracting term}). In the first part of the algorithm, we go through the points of the curve, marking the pixels in the background grids that correspond to these points as $-1$. Then we go through each grid point and check whether its pixel has been marked as $-1$. If so, we assign that pixel the index $i$ of the current grid point. If we are in a self-intersecting region, a grid point with index $j$ will be positioned in a pixel already marked with index $i$, where $j > i$. We define a self-intersection if $j - i > 4$, i.e. there are at least five other grid points between the two intersecting points. \\
\begin{figure}[htbp]
\centering
   \includegraphics[width=0.6\linewidth]{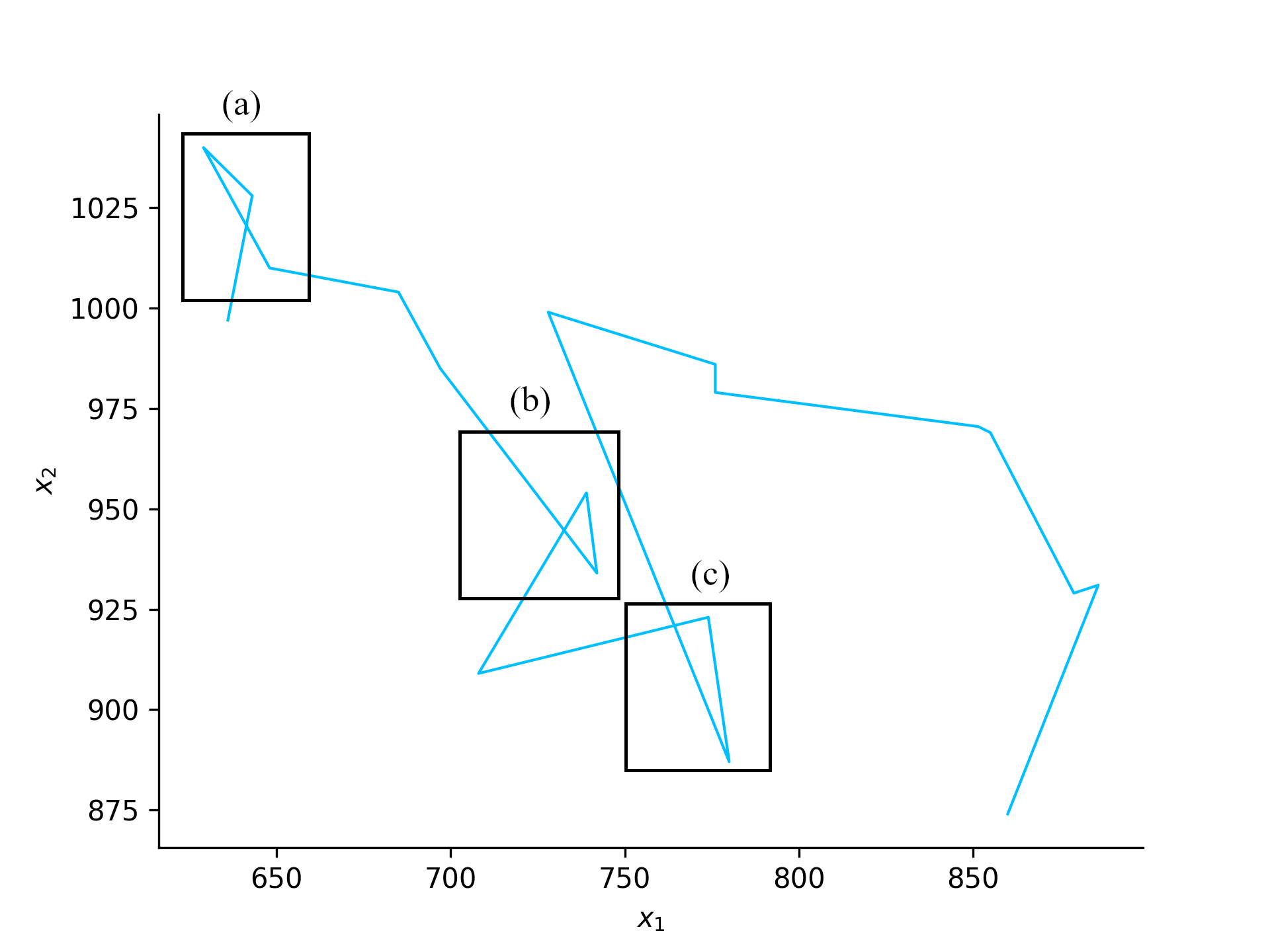} 
   \newline
\subfigure[]{\includegraphics[width=0.3\linewidth]{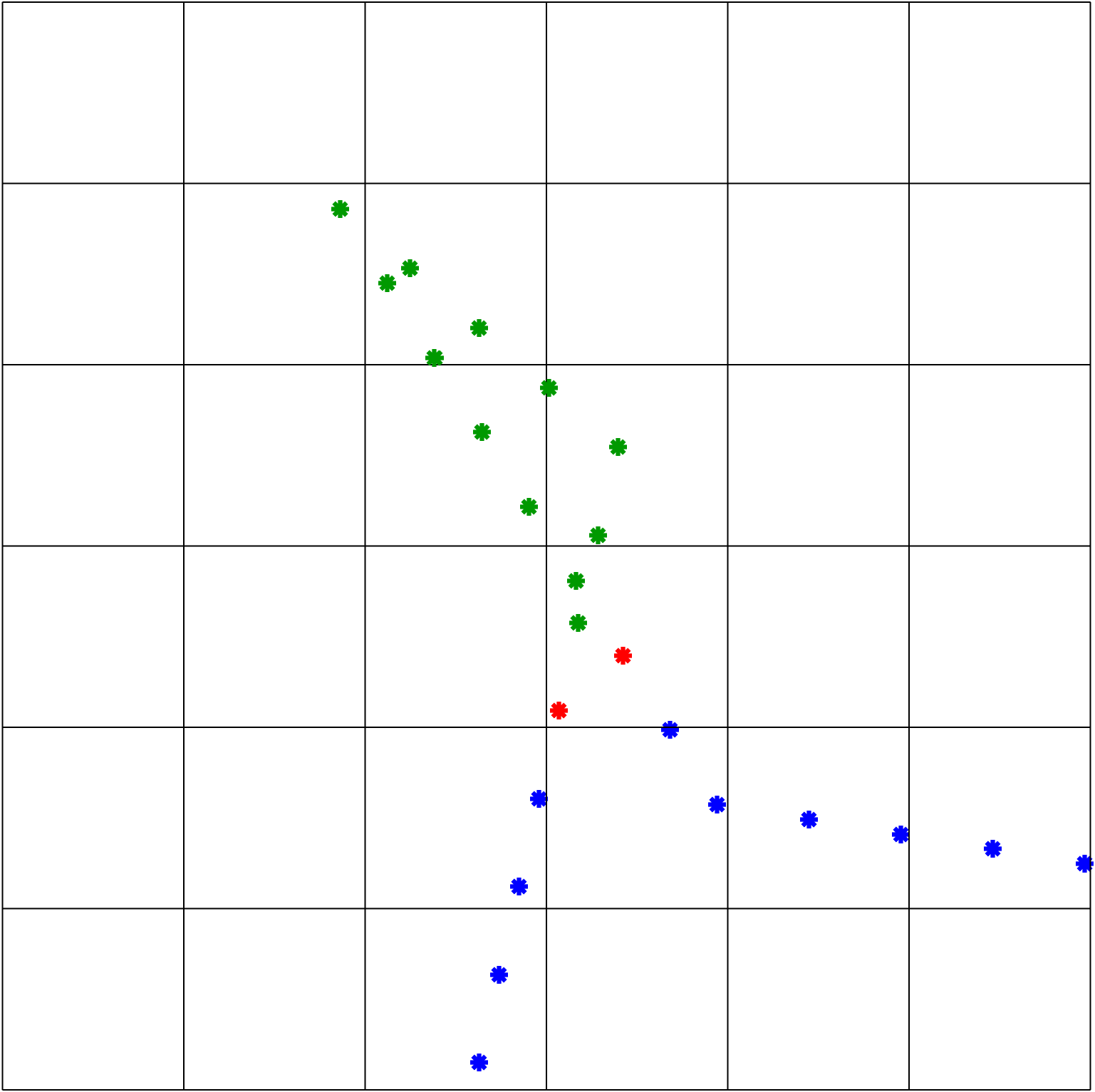}
\label{F: sub1}}
\subfigure[]{\includegraphics[width=0.3\linewidth]{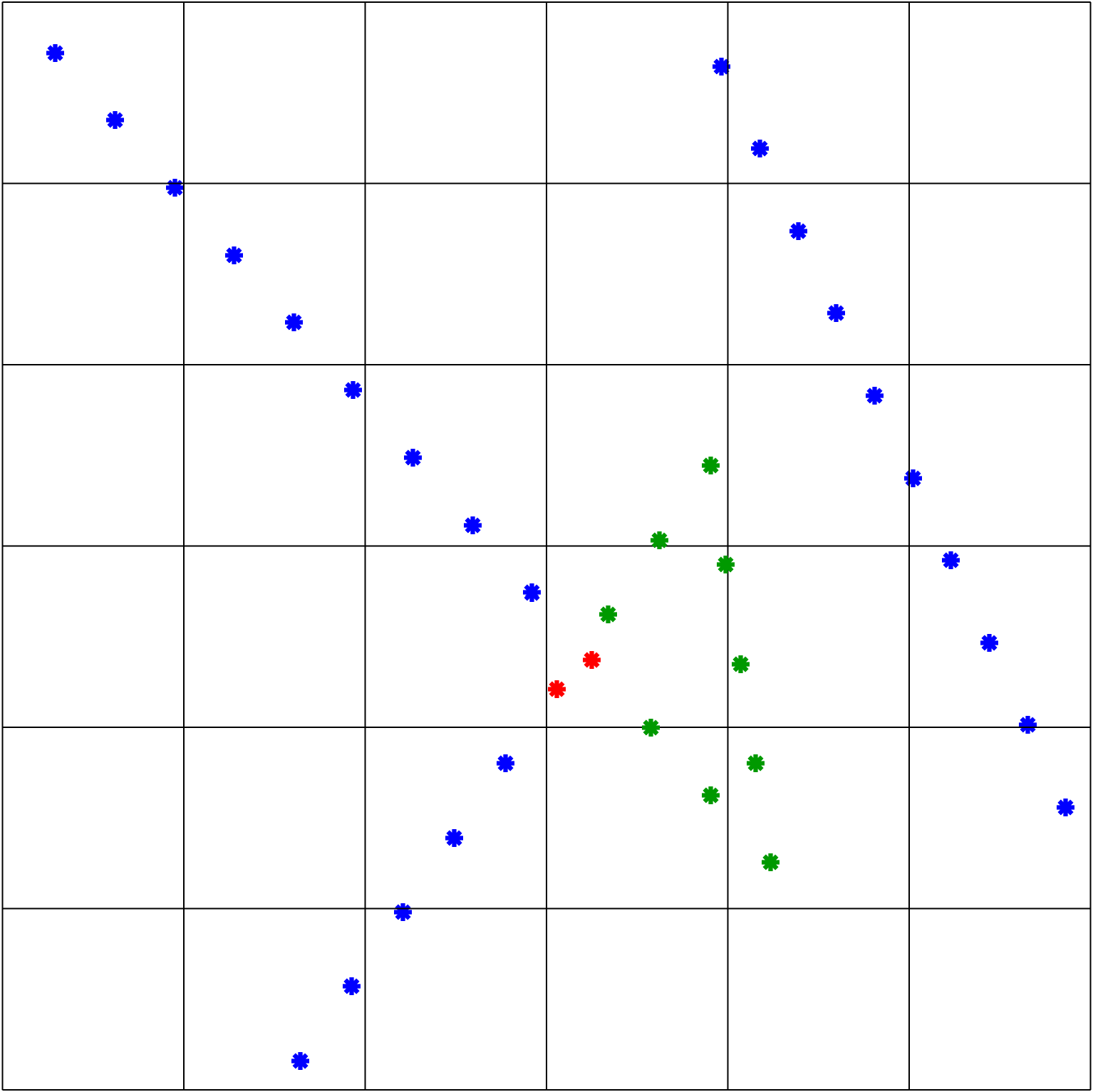}
\label{F: sub2}}
\subfigure[]{\includegraphics[width=0.3\linewidth]{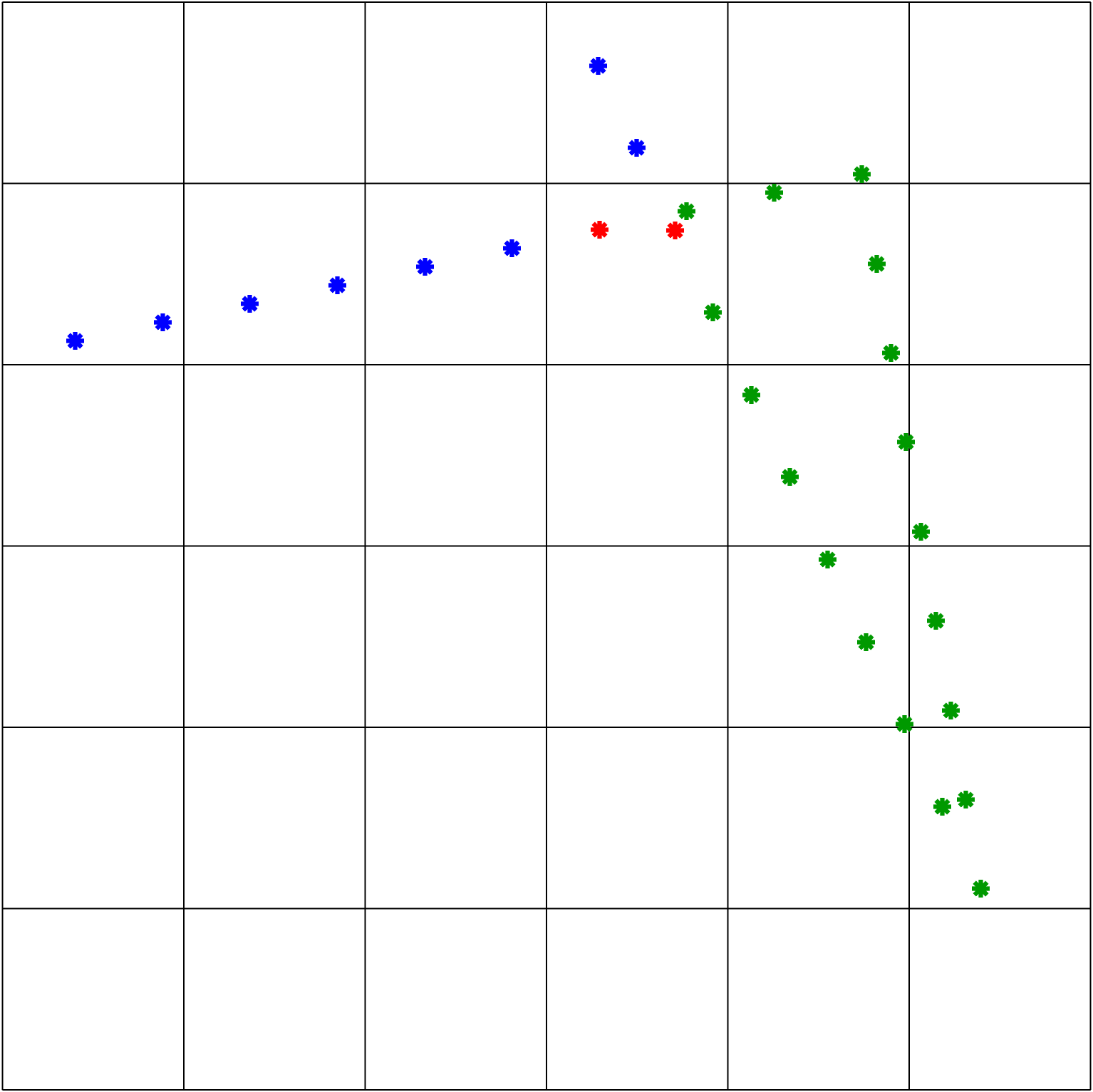}
\label{F: sub3}}
 \caption{Top: original trajectory (light blue curve) and regions of self-intersections (black rectangles). Bottom: zooms of the rectangles. Grid points of the trajectory (blue asterisks), grid points where the self-intersection is detected (red asterisks), and grid points belonging to a self-intersecting part (green asterisks).}
 \label{F: Self intersections}
\end{figure}
Fig. \ref{F: Self intersections} shows the process for a curve; only one of the background grids is shown for clarity. The light blue curve is the original trajectory and the rectangles show the regions of self-intersection. Figs. \ref{F: sub1}, \ref{F: sub2}, \ref{F: sub3} show the zooms of such rectangles. 
In fact, they are not consecutive grid points detected in the same pixel. The blue asterisks are the grid points of the trajectory, the red ones are the grid points where the self-intersection is detected and the green ones are the grid points belonging to the selected self-intersecting part. Note that, as it happens in Fig. \ref{F: sub1}, for a self-intersecting part there may be more pairs of grid points detected as self-intersections. In such a case, we select the start and end points so that the selected self-intersecting part contains all the others. 
The background grid takes into account the initial distance between the points, which is fixed at the beginning. During the evolution, the points move in the normal direction and also along the curve due to tangential redistribution. As a result, the average distance between the points changes. For this reason, before detecting the self-intersections, we add another step to the algorithm and remove every second point if the mean distance is less than $\bar{h}/2$. In the next iteration, the tangential velocity redistributes the points and the mean distance is kept closer to $\bar{h}$. The self-intersection detection is done at each time step and requires $2n$ operations, so the computational cost at each time step is $O(n)$, where $n$ is the number of grid points, and thus it increases the necessary computational cost optimally. \\

 \begin{figure}[htbp]
\centering
   \includegraphics[width=10cm]{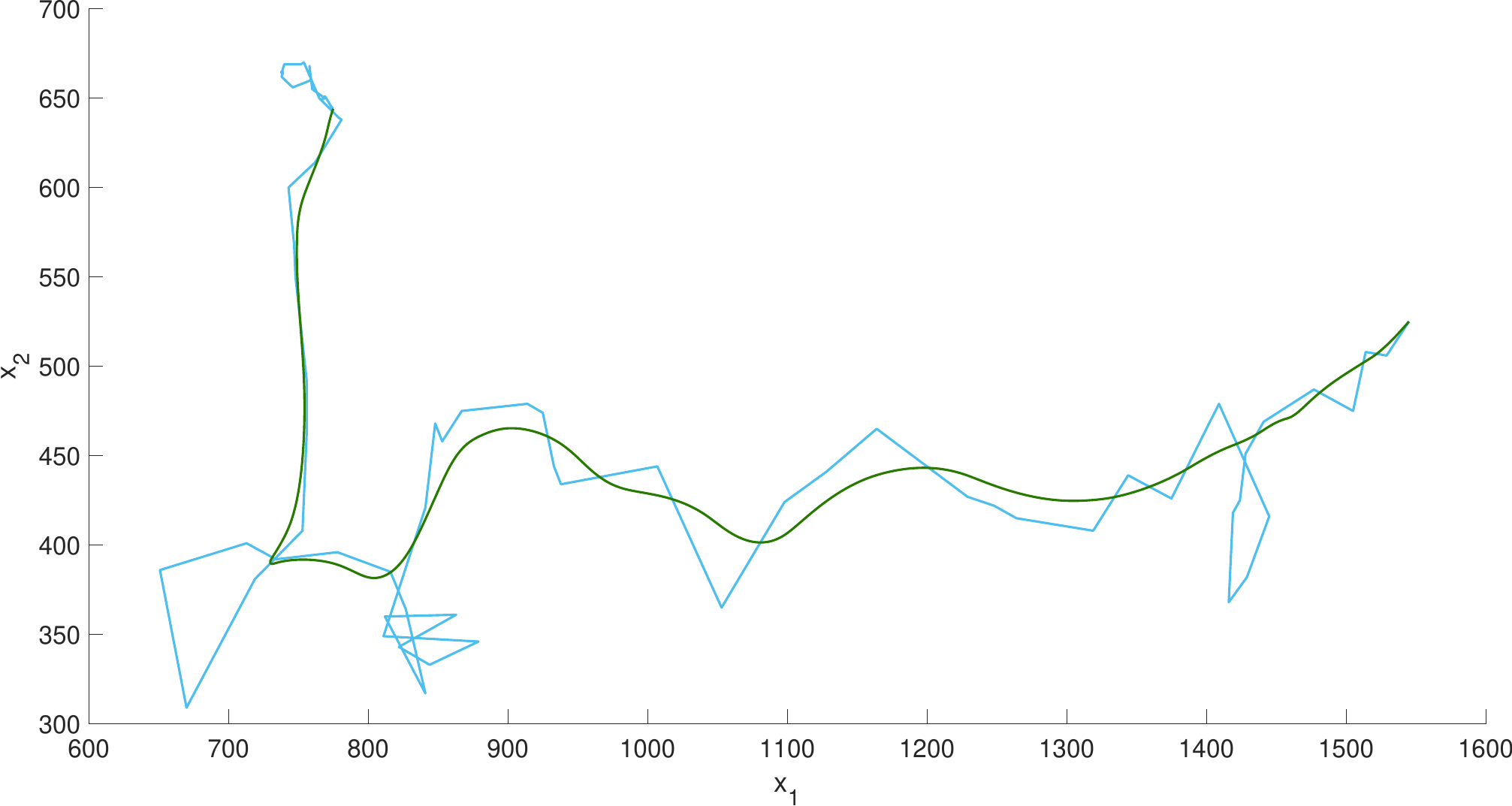} 
\includegraphics[width=10cm]{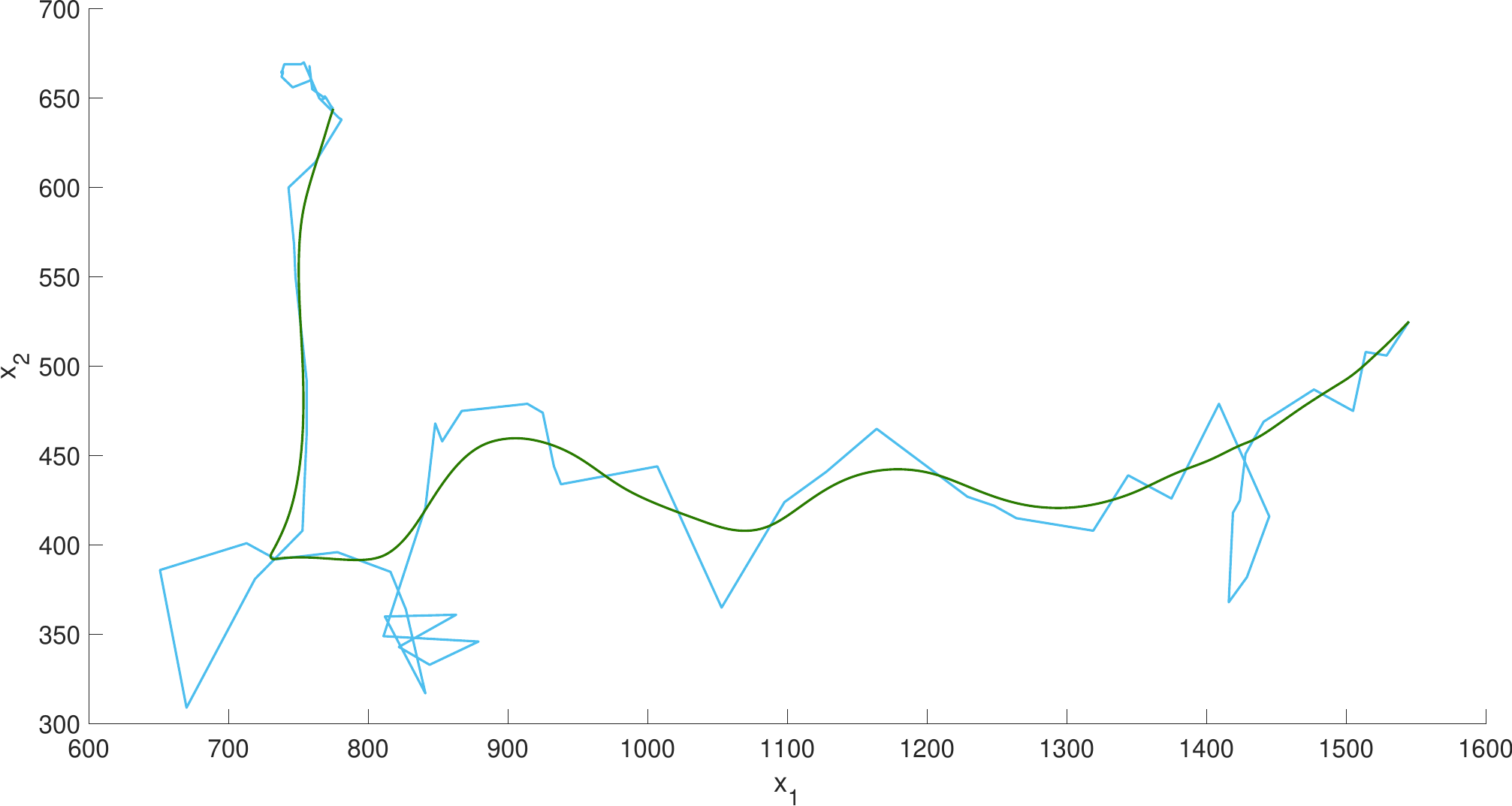}
 \caption{Original trajectory (light blue line) and results of smoothing (green lines) for $\delta=0.003$ (top) and $\delta=0.01$ (bottom). When no further self-intersections are detected, the algorithm performs $50$ additional curve evolution steps before stopping.}
 \label{F:T2}
\end{figure}

 \begin{figure}[htbp]
\centering
   \includegraphics[width=10cm]{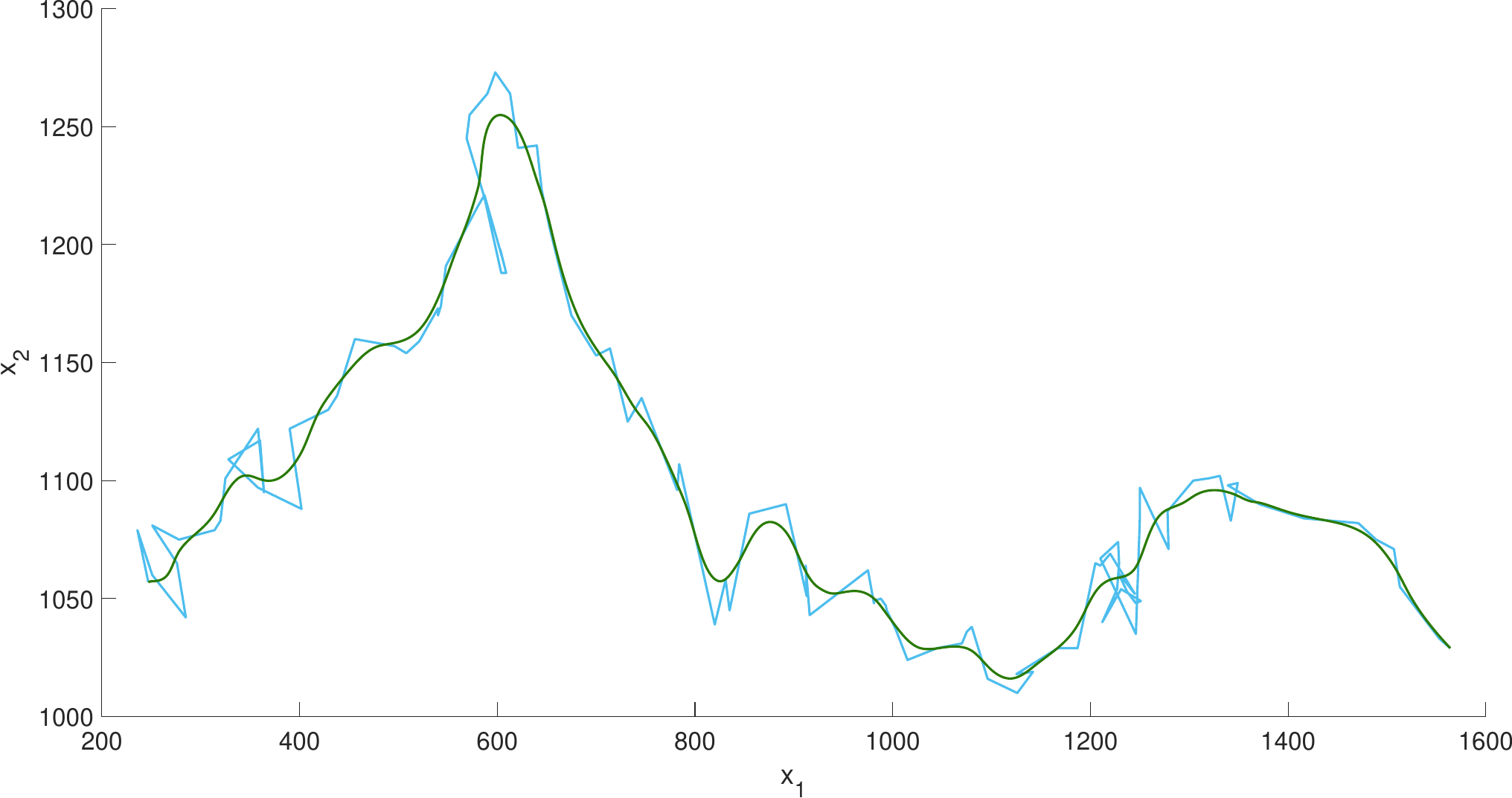} 
\includegraphics[width=10cm]{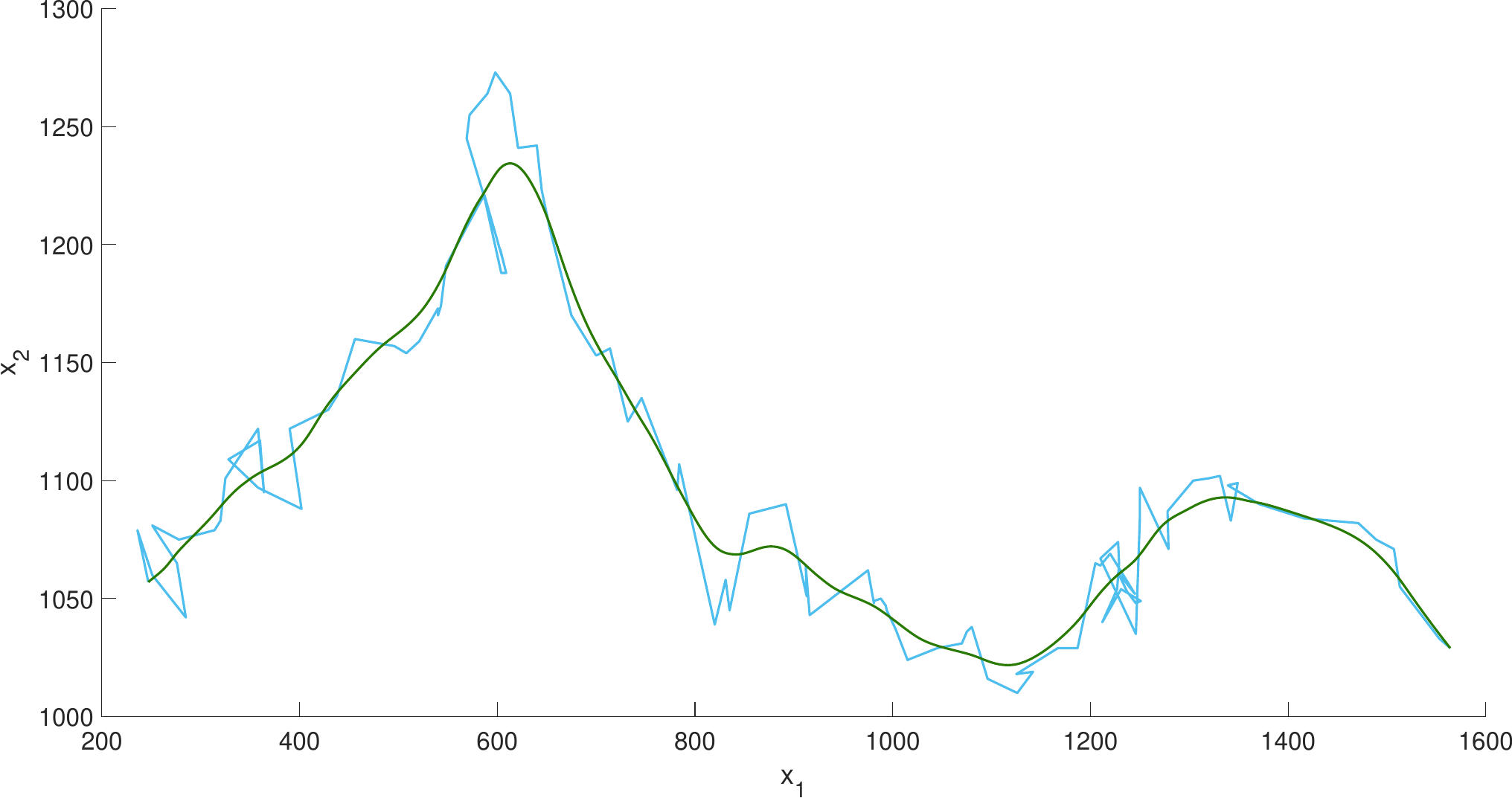}
 \caption{Original trajectory (light blue line) and results of smoothing (green lines) for $\delta=0.003$ (top) and $\delta=0.01$ (bottom). When no further self-intersections are detected, the algorithm performs $50$ additional curve evolution steps before stopping.}
 \label{F:T27}
\end{figure}

 \begin{figure}[htbp]
\centering
   \includegraphics[width=10cm]{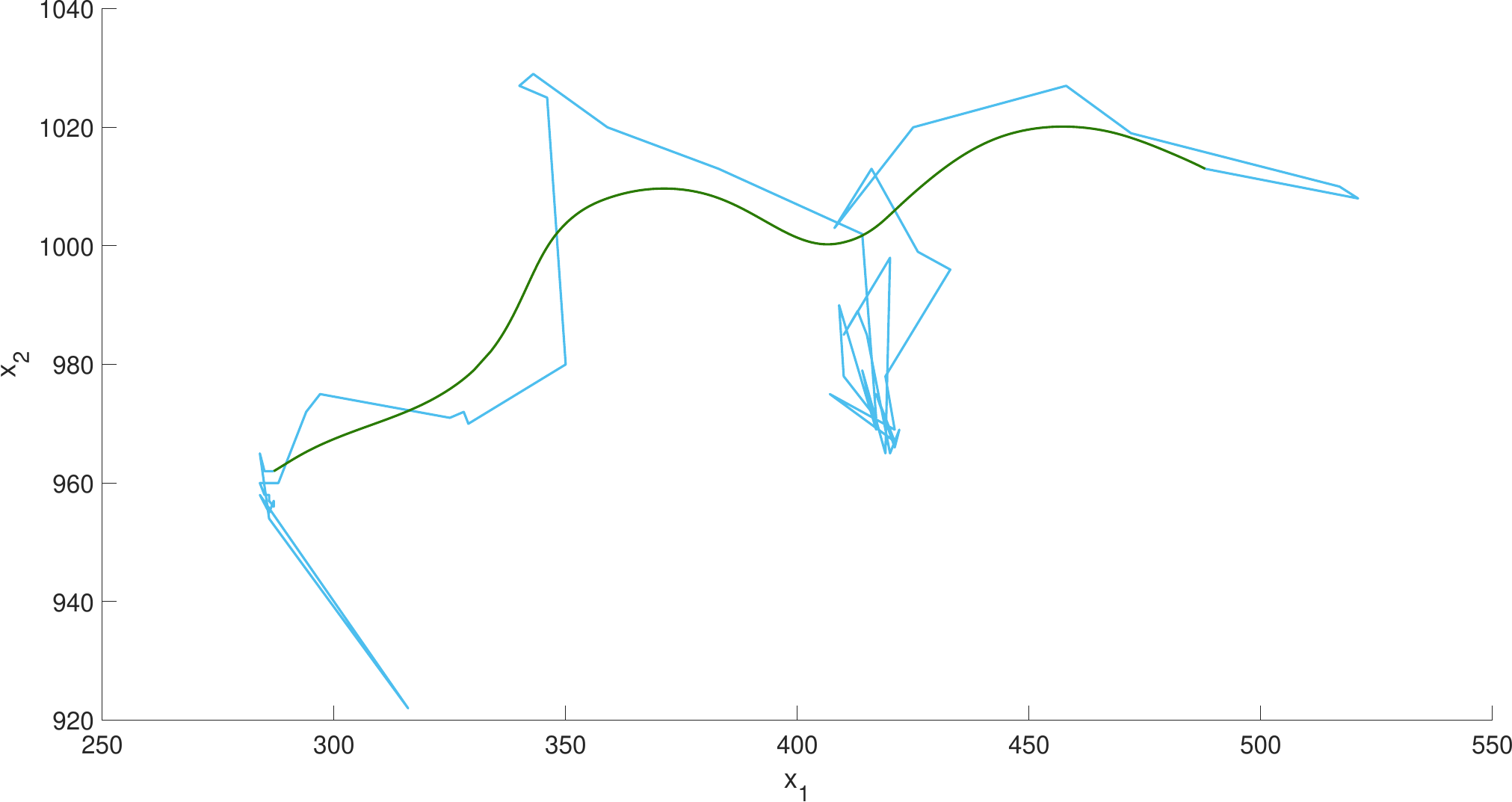} 
\includegraphics[width=10cm]{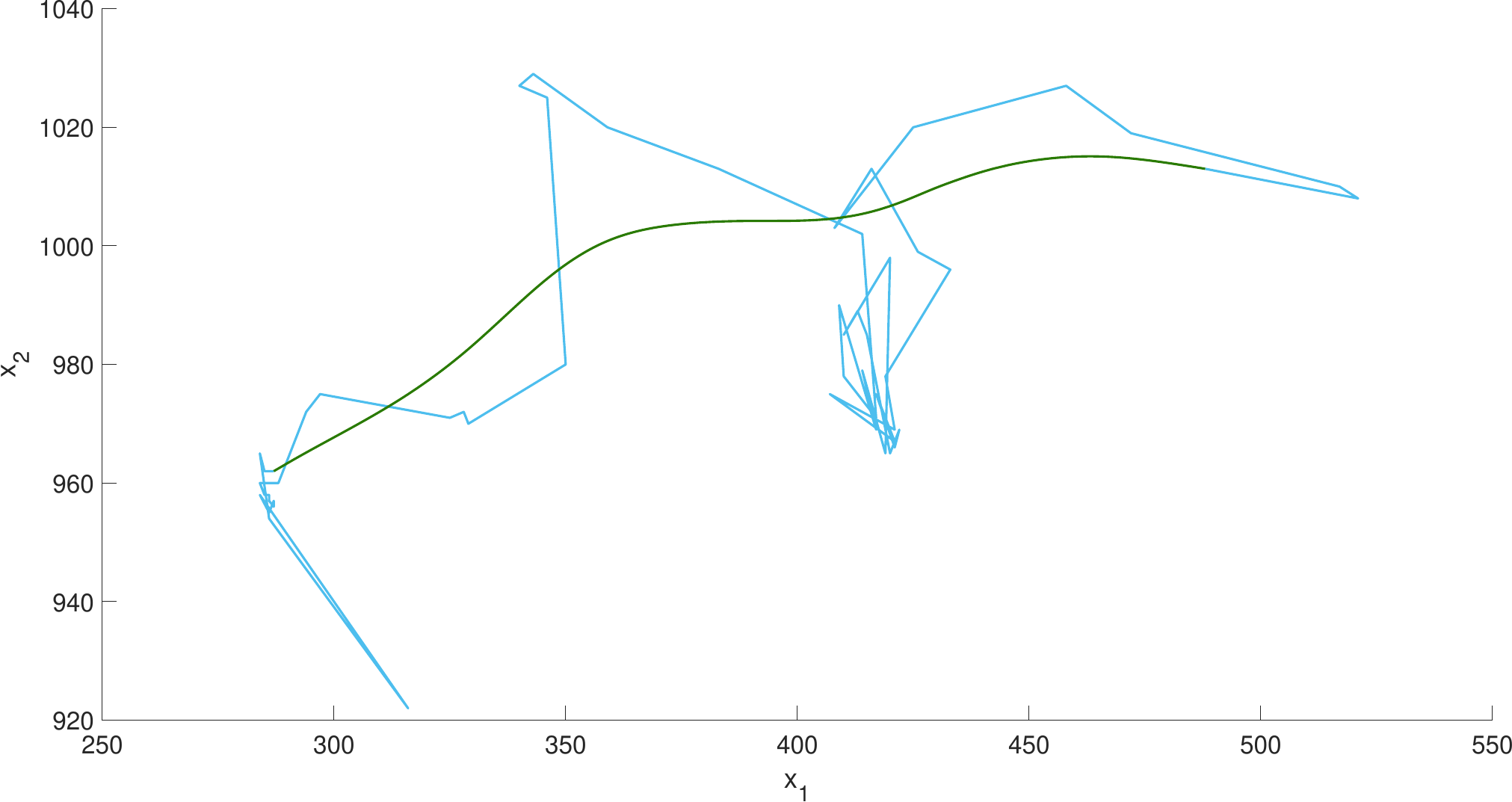}
 \caption{Original trajectory (light blue line) and results of smoothing (green lines) for $\delta=0.003$ (top) and $\delta=0.01$ (bottom). When no further self-intersections are detected, the algorithm performs $50$ additional curve evolution steps before stopping.}
 \label{F:T63}
\end{figure}
Figs. \ref{F:T2}, \ref{F:T27}, \ref{F:T63} show the result of the smoothing using the curve self-intersections as a stopping criterion. When no more self-intersections are detected, the algorithm runs for $50$ more curve evolution steps and stops. We added this step to avoid having sharp angles in the regions where the self-intersections are resolved. In Figs. \ref{F:T2}, \ref{F:T27}, \ref{F:T63} we have compared the results for $2$ different values of $\delta$, namely $\delta=0.003$ and $\delta=0.01$. The other parameters were chosen as follows: $\lambda=20$, $\tau=0.000001$, and $\omega=50$. The light blue lines are the original trajectories, while the green curves are the smoothed curves; the result of the curve evolution with $\delta=0.003$ is shown at the top, while that with $\delta=0.01$ is shown at the bottom. As one can see, a smaller value of $\delta$ keeps the trajectories closer to the original ones.
 \begin{figure}[htbp]
\centering
   \includegraphics[width=10cm]{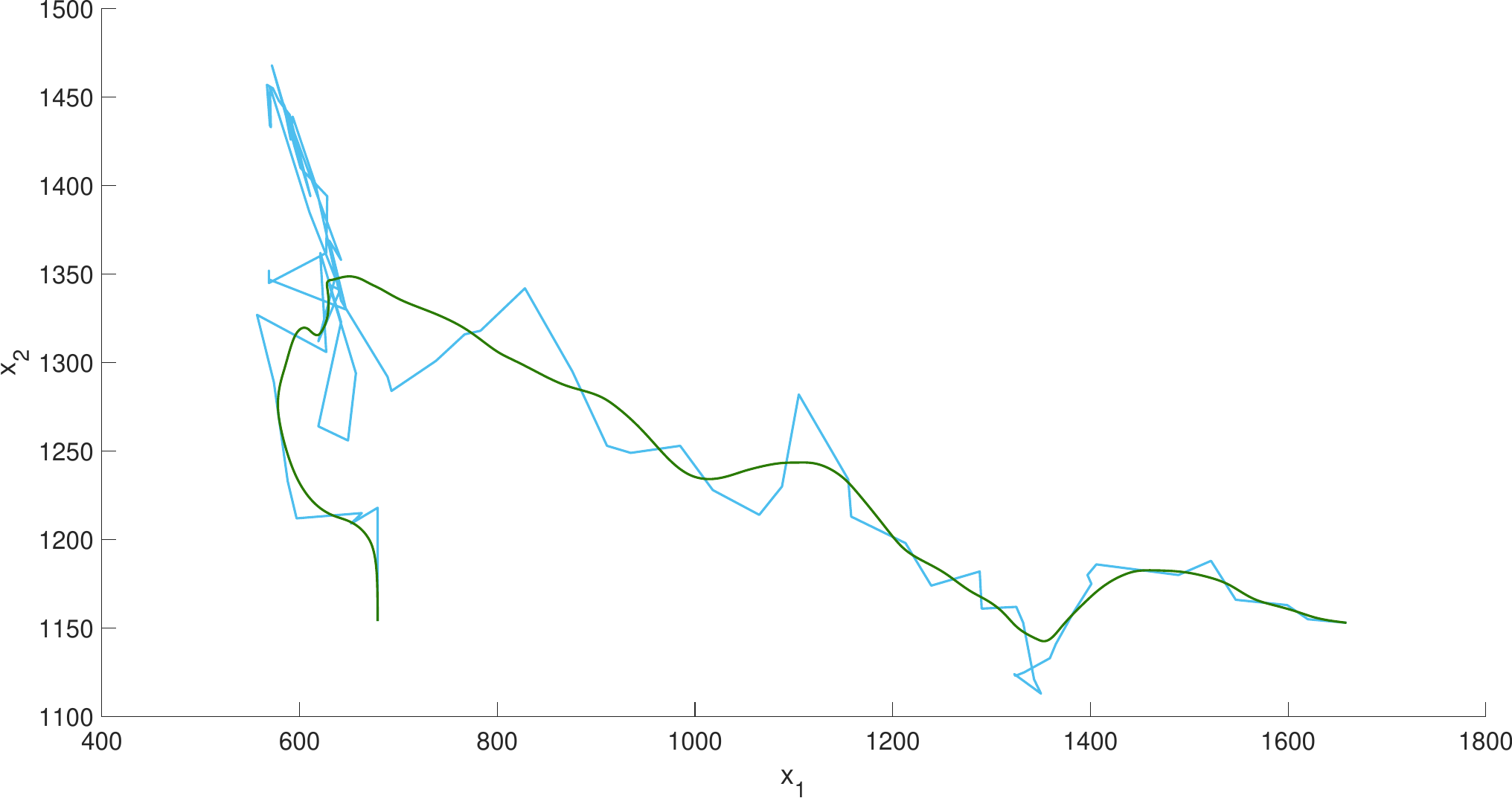} 
\includegraphics[width=10cm]{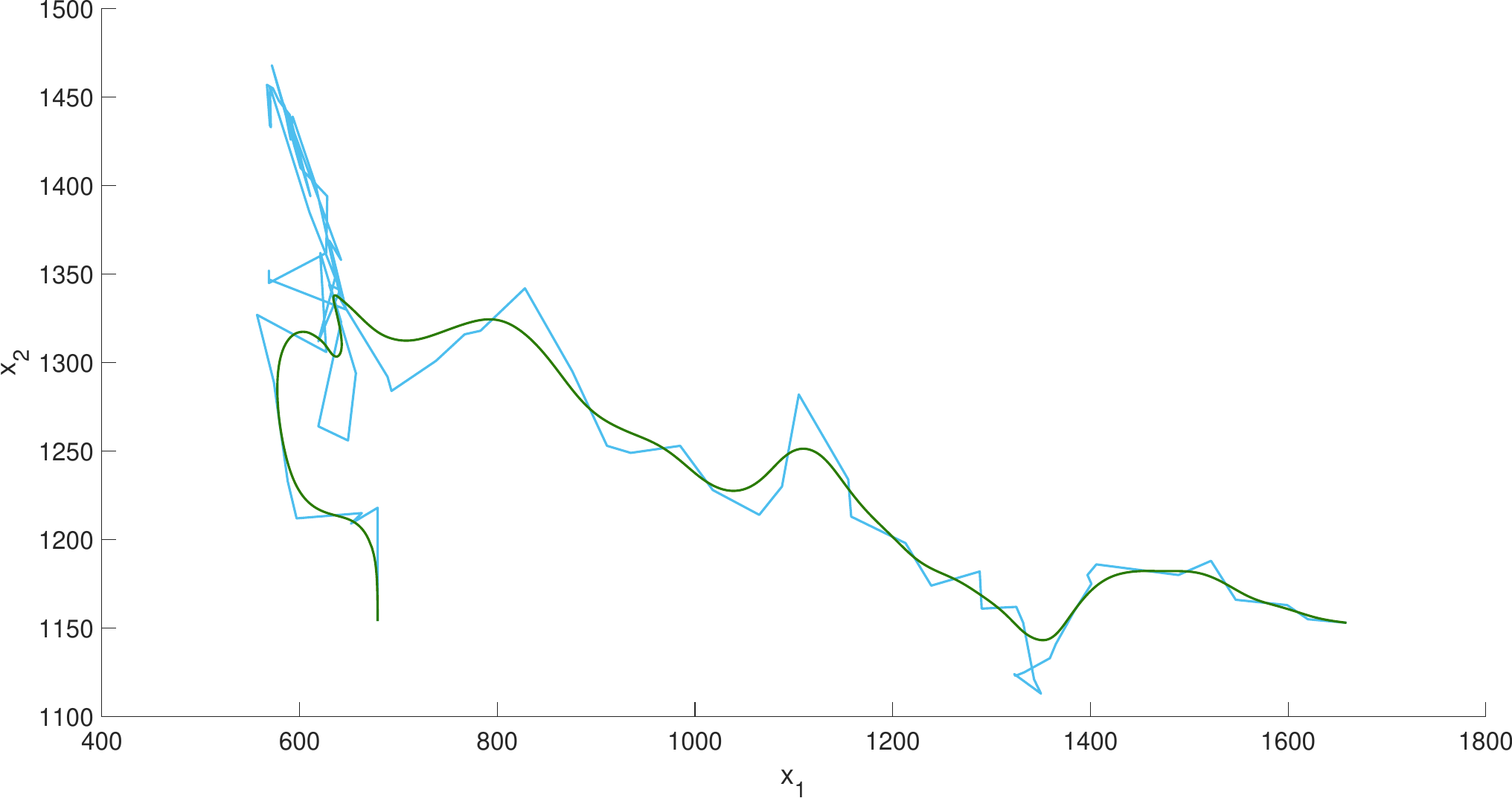}
 \caption{Original trajectory (light blue line) and results of smoothing (green lines) for $\delta=0.003$ (top) and $\delta=0.01$ (bottom). When no further self-intersections are detected, the algorithm performs $50$ additional curve evolution steps before stopping.}
 \label{F:T32}
\end{figure}
However, this is not always the case. In fact, the resulting evolved curve may be smoother when considering a smaller value of $\delta$ than when considering a larger value of $\delta$. This is the case, for example, for the trajectory shown in Fig. \ref{F:T32}. Also in this figure, the result of the curve evolution with $\delta=0.003$ is shown at the top, while the one with $\delta=0.01$ is shown at the bottom. As one can notice, this time the curve obtained with a smaller value of $\delta$ is smoother. This is due to the number of curve evolution steps needed to solve the self-intersections. For $\delta=0.003$ the algorithm ran $24638$ curve evolution steps, while for $\delta=0.01$ the number of curve evolution steps was $6728$. 
Since the algorithm takes longer to remove the self-intersections for small values of $\delta$, the resulting curve is too smooth and we lose details of the directional parts of the motion. For this reason, we considered an adaptive choice of parameters. The details of this choice are explained in the next section.

\subsubsection{Adaptive choice of parameters}
Once we have detected the self-intersections, we can use the results not only as a stopping criterion but also to adaptively select the smoothing parameters. In fact, our goal is to smooth the trajectories in the self-intersecting parts and have almost no smoothing for the directional parts. Therefore we choose $\delta$ and $\lambda$ in the following way: consider $0\leq\delta_{min}<\delta_{max}$, and $\lambda_{max}>0$. Let us indicate by $i_1$ and $i_2$ the indexes of a detected self-intersection, with $i_1<i_2$. Then
\begin{equation}
\left\{
\begin{array}{l}
\delta_i=\frac{(6-i_1+i)\delta_{max}}{6},~ \lambda_i=\frac{(6-i_1+i)\lambda_{max}}{6}~~~for~ i=i_1-5,...,i_1 \\
\delta_i=\delta_{max},~\lambda_i=\lambda_{max}~~~~~~~~~~~~~~~~~~~~~~for~i=i_1,...,i_2 \\
\delta_i=\frac{(6+i_2-i)\delta_{max}}{6},~ \lambda_i=\frac{(6+i_2-i)\lambda_{max}}{6}~~~for~i=i_2,...,i_2+5 \\
\delta_i=0,~ \lambda_i=0~~~~~~~~~~~~~~~~~~~~~~~~~~~~~~~~otherwise.
\end{array}
\right.
\end{equation}
Notice that, for the grid points belonging to directional parts of motion, we will have $\beta_i=0$ but $\alpha_i\neq 0$, so that the points will move in the tangential direction and redistribute along the curve. Once no more self-intersections are detected, the algorithm runs for $50$ more curve evolution steps with $\delta_i=\delta_{min}$ and $\lambda_i=\lambda_{max}$ for $i=1,...,n$. The reason for this is that we want to remove the sharp angles in the original trajectories and get a smooth curve. \\
For the smoothing of macrophage trajectories, the parameters were chosen as follows: $\delta_{min}=0.003$, $\delta_{max}=0.01$, $\lambda_{max}=20$, $\tau=0.000001$, and $\omega=50$. 
 \begin{figure}[htbp]
\centering
   \includegraphics[width=10cm]{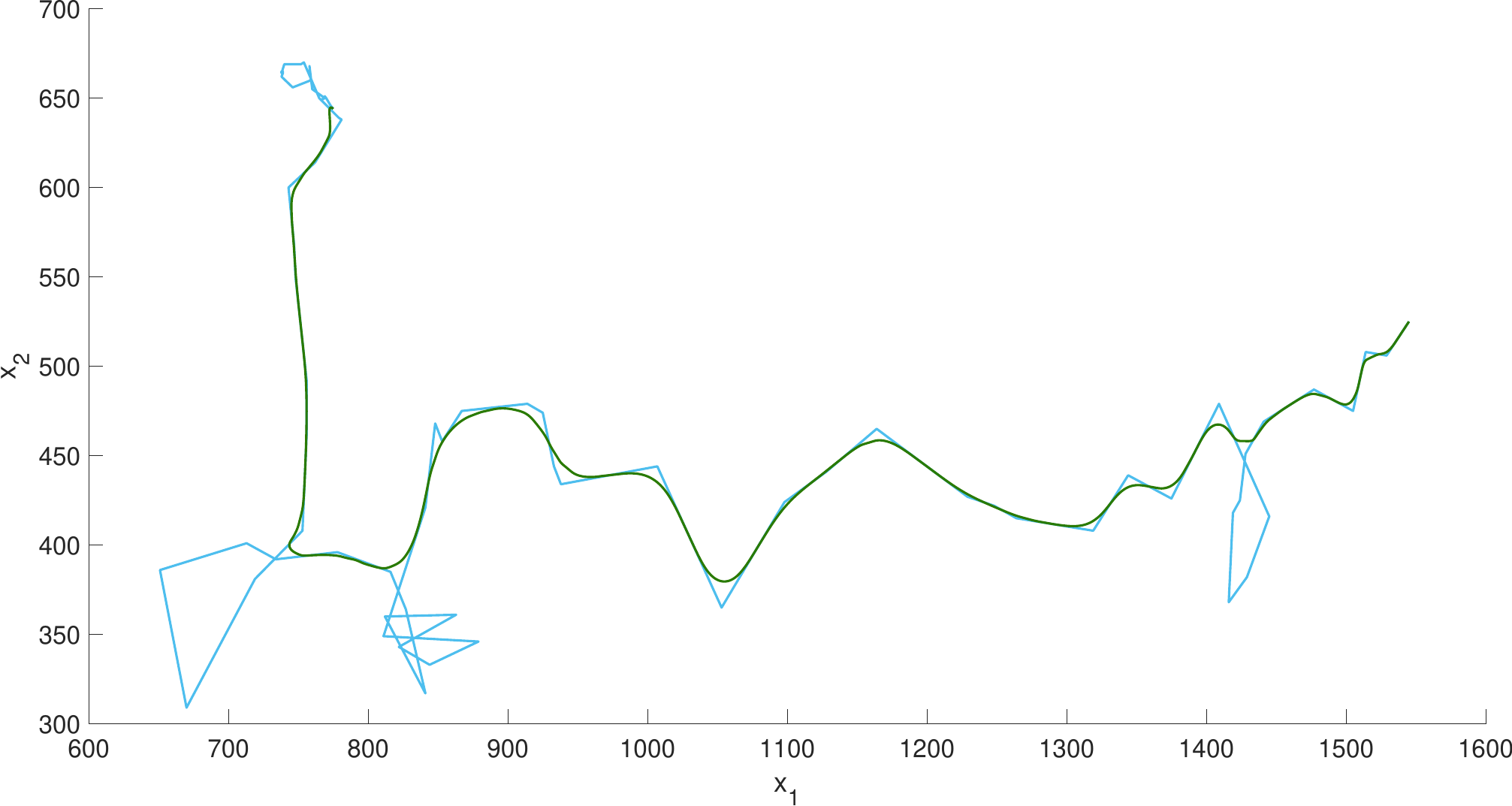} 
 \caption{Original trajectory (light blue line) and results of smoothing for adaptive choice of parameters (green line). When no further self-intersections are detected, the algorithm performs $50$ additional curve evolution steps before stopping.}
 \label{F:T2 Mix}
\end{figure}

 \begin{figure}[htbp]
\centering
   \includegraphics[width=10cm]{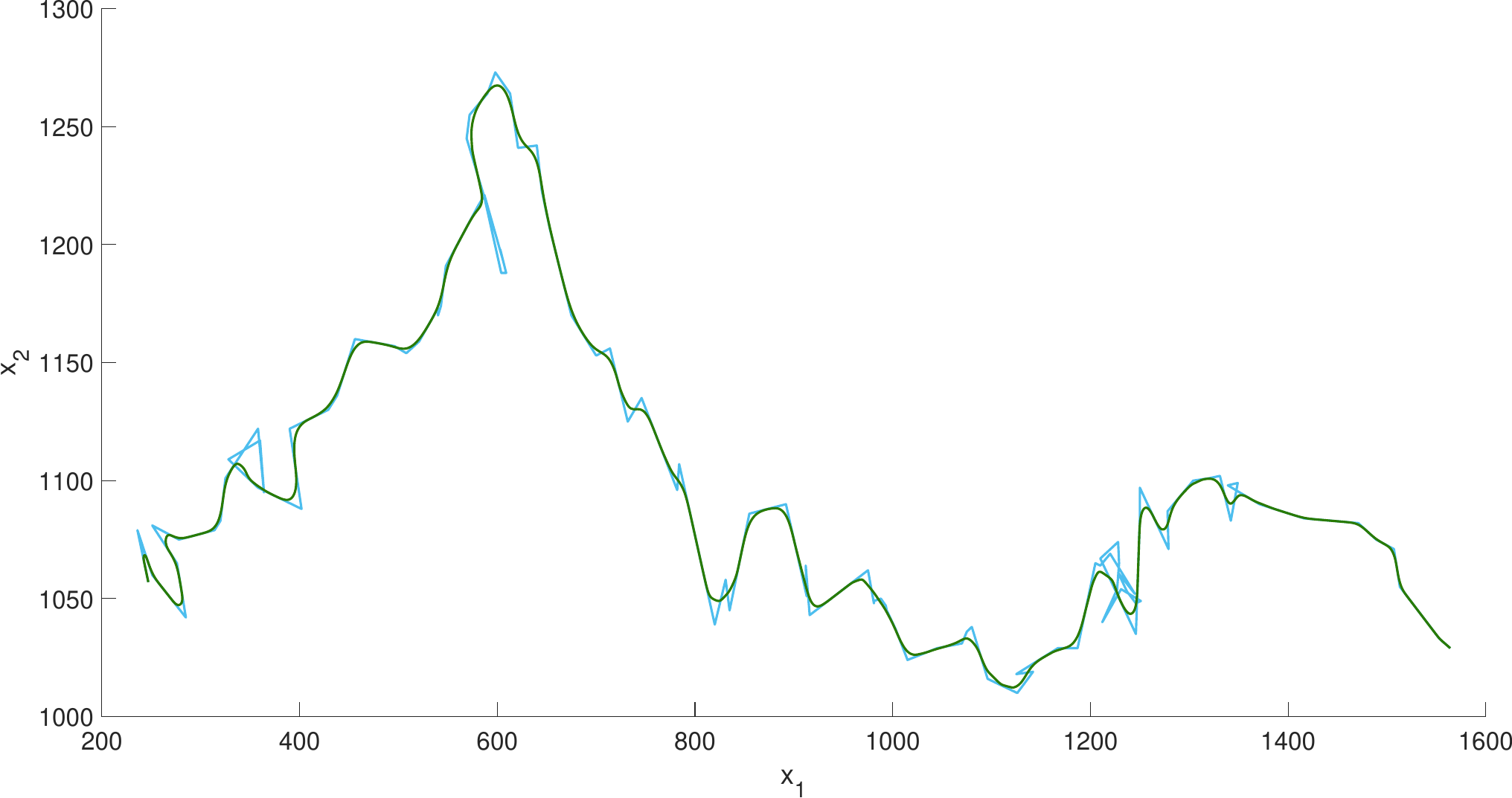} 
 \caption{Original trajectory (light blue line) and results of smoothing for adaptive choice of parameters (green line). When no further self-intersections are detected, the algorithm performs $50$ additional curve evolution steps before stopping.}
 \label{F:T27 Mix}
\end{figure}
 \begin{figure}[htbp]
\centering
   \includegraphics[width=10cm]{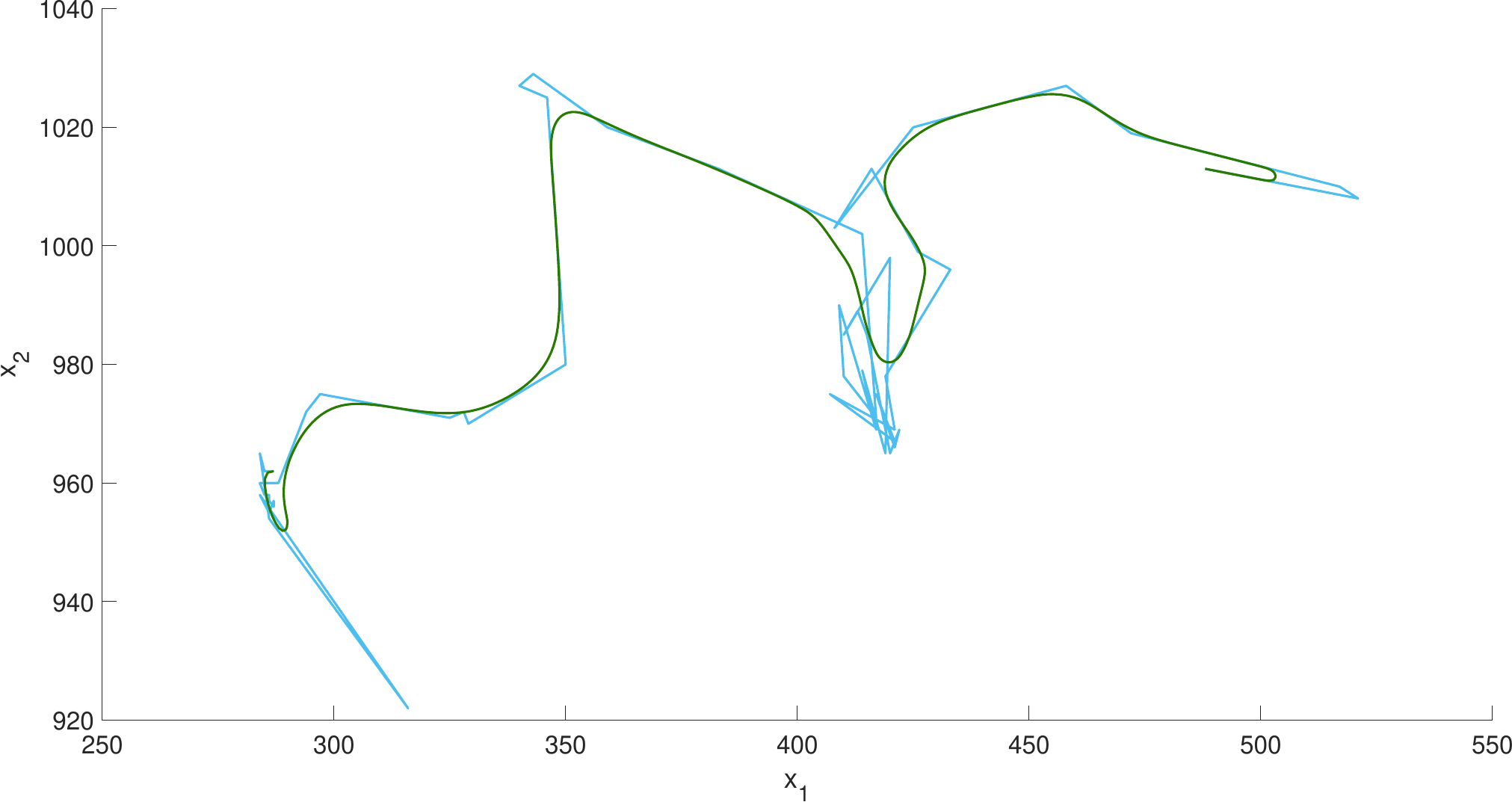} 
 \caption{Original trajectory (light blue line) and results of smoothing for adaptive choice of parameters (green line). When no further self-intersections are detected, the algorithm performs $50$ additional curve evolution steps before stopping.}
 \label{F:T63 Mix}
\end{figure}
 \begin{figure}[htbp]
\centering
   \includegraphics[width=10cm]{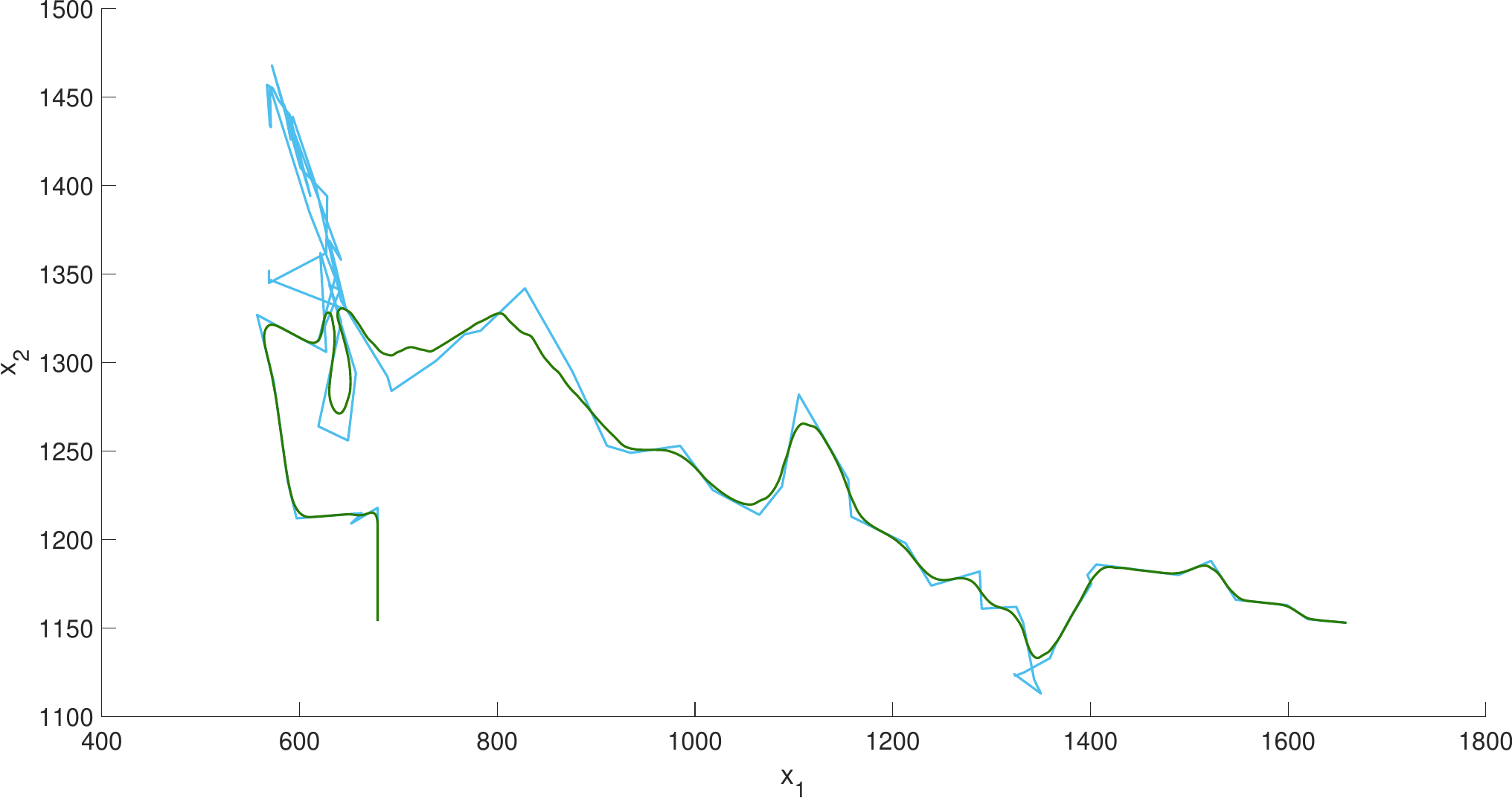} 
 \caption{Original trajectory (light blue line) and results of smoothing for adaptive choice of parameters (green line). When no further self-intersections are detected, the algorithm performs $50$ additional curve evolution steps before stopping.}
 \label{F:T32 Mix}
\end{figure}
Figs. \ref{F:T2 Mix}, \ref{F:T27 Mix}, \ref{F:T63 Mix}, \ref{F:T32 Mix} show the results of the trajectories considered in Figs. \ref{F:T2}, \ref{F:T27}, \ref{F:T63} and \ref{F:T32} respectively. The light blue curve is the original trajectory, while the green one is the smoothed curve. The adaptive choice of parameters allows us to keep the curve close to the original trajectory in the directional parts and to remove the self-intersections.

\section{Analysis of random parts}
\label{S: Analysis of random parts}
\begin{figure}[htbp]
    \centering
    \includegraphics[width=0.8\linewidth]{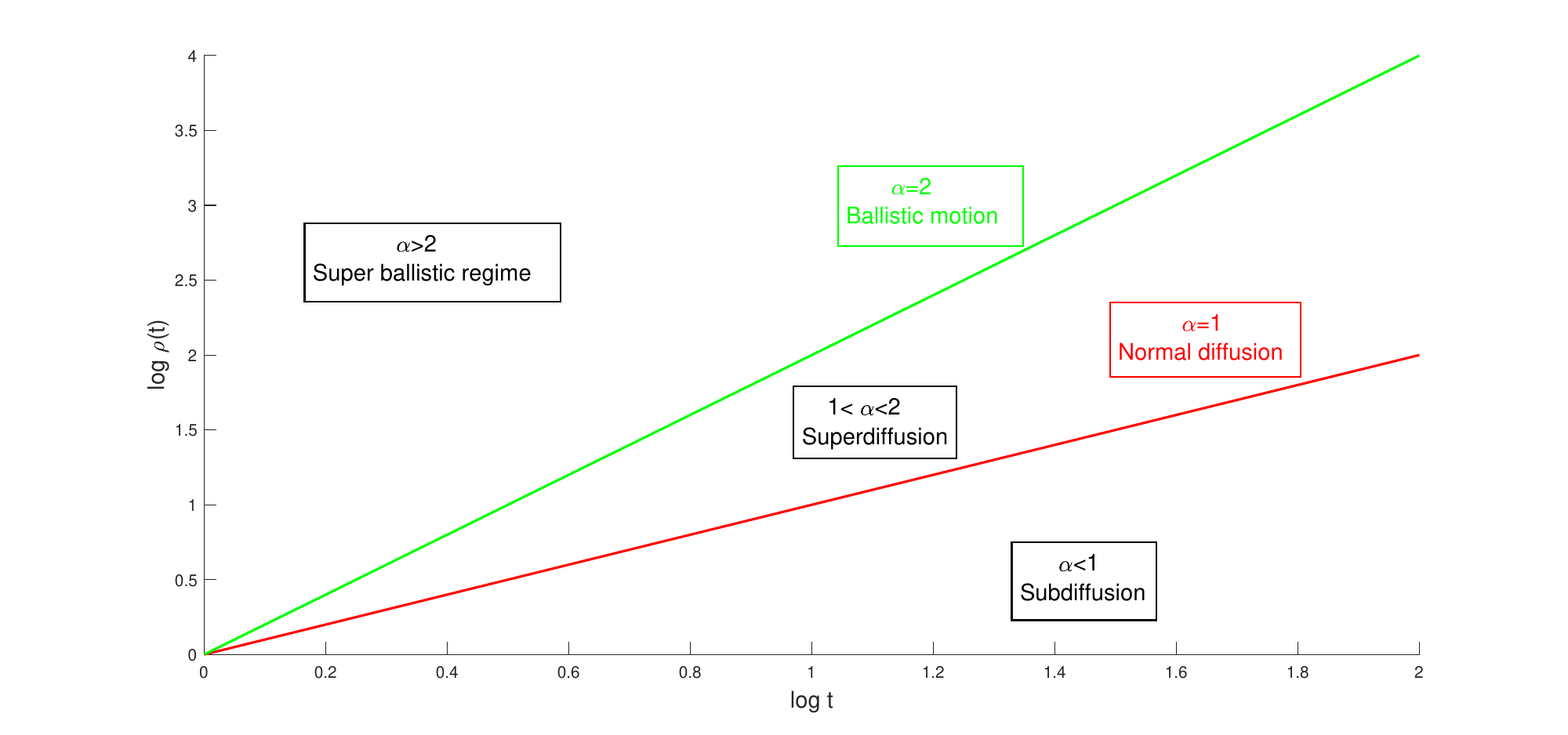}
    \caption{Plot of MSD versus time in log-log scale. Normal diffusion (red line) corresponds to slope $\alpha=1$ ($H=\frac{1}{2}$), and ballistic motion (green line) corresponds to slope $\alpha=2$ ($H=1$). $\alpha>2$ gives the superballistic regime, while $1<\alpha<2$ gives superdiffusion. Subdiffusion is obtained for $0<\alpha<1$.}
    \label{F: MSD}
\end{figure}

Once the trajectories have been smoothed, we want to analyze the random parts. A common tool to analyze random walks is the so-called mean squared displacement (MSD); indeed, it allows us to distinguish the type of diffusion the random walker shows. The motion can be either diffusive, subdiffusive, or superdiffusive. For example, Brownian motion is a type of diffusive process; examples of anomalous diffusion, superdiffusion, and subdiffusion are the Levy flights or the random walks in crowded environments, respectively \cite{viswanathan2011physics,sokolov2012models}.\\
The mean squared displacement is a measure of the deviation of a particle concerning a reference point. It is defined as
\begin{equation}
    \rho(t)  = \langle\lvert\textbf{x}(t)-\textbf{x}^0\rvert^2 \rangle,
\end{equation}
where $\textbf{x}^0$ is the reference point that can be either the initial point of the trajectory or just a point of the trajectory as it will be discussed in the next section. \\
Consider first normal diffusion; the probability density function of the positions of particles is Gaussian and the variance grows linearly with time. It can be proved \cite{sokolov2012models} that the MSD grows as
\begin{equation}
\label{EQ: MSD BM}
    \rho(t) \sim t,
\end{equation}
in particular
\begin{equation}
   \rho(t)=2dDt,
\end{equation}
where $d$ is the dimension of the space, $D$ is the diffusion coefficient and $t$ is time. \\
To characterize anomalous diffusion one usually defines the \textit{Hurst exponent} $H$
\begin{equation}
\label{EQ: Hurst exponent}
    \rho(t) \sim t^{2H}
\end{equation}
to see how the MSD grows with time. Then
\begin{itemize}
    \item For $0<H<\frac{1}{2}$ the random motion is subdiffusive.
    \item For $H=\frac{1}{2}$ one obtains normal diffusion.
    \item For $H>\frac{1}{2}$ the random motion is superdiffusive.
    \item For $H=1$ one obtains ballistic motion.
    \item For $H>1$ the random motion is superballistic.
\end{itemize}
Therefore, to distinguish the type of diffusion followed by the random walker one has to find $H$. Considering the logarithm
\begin{equation}
    \log \rho(t) \sim 2H \log t,
\end{equation}
and the dependence on time of the MSD in a log-log scale, results in $\alpha=2H$ being the slope of the line obtained in the double-log scaling. Fig. \ref{F: MSD} shows the plot of MSD versus time in the log-log scale.
In the next section we will give two different discrete definitions of $\rho(t)$ and discuss their differences.\\

\subsection{Mean squared displacement}
\label{SubS: msd}
The first way to define the MSD is the so-called \textit{ensemble-averaged} MSD (EAMSD). Consider a set of trajectories. Every trajectory has an initial point: for the $i$-th trajectory it will be indicated as $\textbf{x}_{i}^{0}$ and will be considered the reference point for that trajectory. Then, the EAMSD at time $t$ is defined as
\begin{equation}
    \rho(t)=\frac{1}{N}\sum_{i=1}^{N}\lVert \textbf{x}_i(t)-\textbf{x}_{i}^{0} \rVert^2,
\end{equation}
where $N$ is the total number of trajectories \cite{regner2013anomalous}. Notice that, in our case, trajectories do not have the same length therefore for different times $t$ we get different values of $N$. We consider then the same formula but generalized as
\begin{equation}
    \rho(t)=\frac{1}{N_t}\sum_{i=1}^{N_t}\lVert \textbf{x}_i(t)-\textbf{x}_{i}^{0} \rVert^2,
    \label{E: EAMSD}
\end{equation}
where $N_t$ is the number of trajectories at time $t$ and $\textbf{x}_{i}^{0}$ is the initial point of the trajectory.\\
Getting statistically significant results would require many trajectories and long observation times: this is not usually the case in studying biological systems. For this reason, a different definition for the MSD is introduced, the so-called \textit{time-averaged} MSD (TAMSD) \cite{kubala2021diffusion, regner2013anomalous}. In this case, instead of time, we consider \textit{time lags}. Consider the $i$-th trajectory with intial time $t_{i}^{0}$ and final time $T_i$. Consider the data acquisition time interval $\Delta T$: for the $i$-th trajectory we have the cell position at times $(n\Delta T)+t_{i}^{0}$ for $n=0,1,2.., \frac{T_i-t_{i}^{0}}{\Delta T}$; here $n$ is the time lag.\\
Let us denote by $K_i=\frac{T_i-t_{i}^{0}}{\Delta T}$ the final time lag for the $i-$th trajectory. The TAMSD is defined as
\begin{equation}
    \bar{\rho}_i(n\Delta T)=\frac{1}{K_i-n+1}\sum_{j=0}^{K_i-n} \lVert \textbf{x}_i((j+n)\Delta T)-\textbf{x}_{i}(j\Delta T) \rVert^2.
\end{equation}
The intuitive idea is to consider many small trajectories inside one single trajectory and make an average of the displacement. In this way, we obtain more values for the MSD from a single trajectory. Notice that with increasing time lag $n$ the number of sub-trajectories inside the original one decreases and this results in fewer data points and more variance in the value of the TAMSD. Therefore, to get meaningful results, the TAMSD for the $i$-th trajectory should always be calculated for time lags $n$ in the following interval \cite{ruthardt2011single}
\begin{equation}
    n\leq\frac{T_i-t_{i}^{0}+1}{4}.
\end{equation}
One can then consider the \textit{ensemble average of the time-averaged} MSD (EATAMSD) \cite{regner2013anomalous}. As in Eq. \ref{E: EAMSD}, we consider the generalized formula
\begin{equation}
\label{E: MSD time all}
    \bar{\rho}(n\Delta T)=\frac{1}{N_n}\sum_{i=1}^{N_n}\bar{\rho}_i(n\Delta T),
\end{equation}
where $N_n$ is the number of trajectories for which $\bar{\rho}_i(n\Delta T)$ exists.
We consider a similar restriction as the one for $n$ for $N_n$. A too-small value of $N_n$ would lead to high variance in the value of $\bar{\rho}(n\Delta T)$, therefore we consider only
\begin{equation}
N_n\geq\frac{N}{4},
\end{equation}
where $N$ is the total number of trajectories.\\
Notice that, for $\bar{\rho}(n\Delta T)$ it holds similar asymptotic behaviors as in (\ref{EQ: MSD BM}) for $\rho(t)$ and the classification of anomalous diffusion with the Hurst exponent as in (\ref{EQ: Hurst exponent}) with the only difference that for $\bar{\rho}(n\Delta T)$, we consider the time lag $n$ and not the time $t$ \cite{regner2013anomalous}.\\

We developed two different methods to detect the random parts, also referred to as random sub-trajectories throughout the text. One considers the evolution of the lengths of the original segments while the other considers the self-intersections in the trajectory. 

\subsection{Random parts extracted by evolving curves}
\label{SS: Random parts extracted by evolving curves}

In this approach, if a segment disappears during smoothing, we consider it to belong to a random part of the motion. Indeed, random segments are usually shorter and in areas of high curvature; this causes their lengths to vanish as they evolve.
\begin{figure}[H]
    \centering
    \includegraphics[width=8cm]{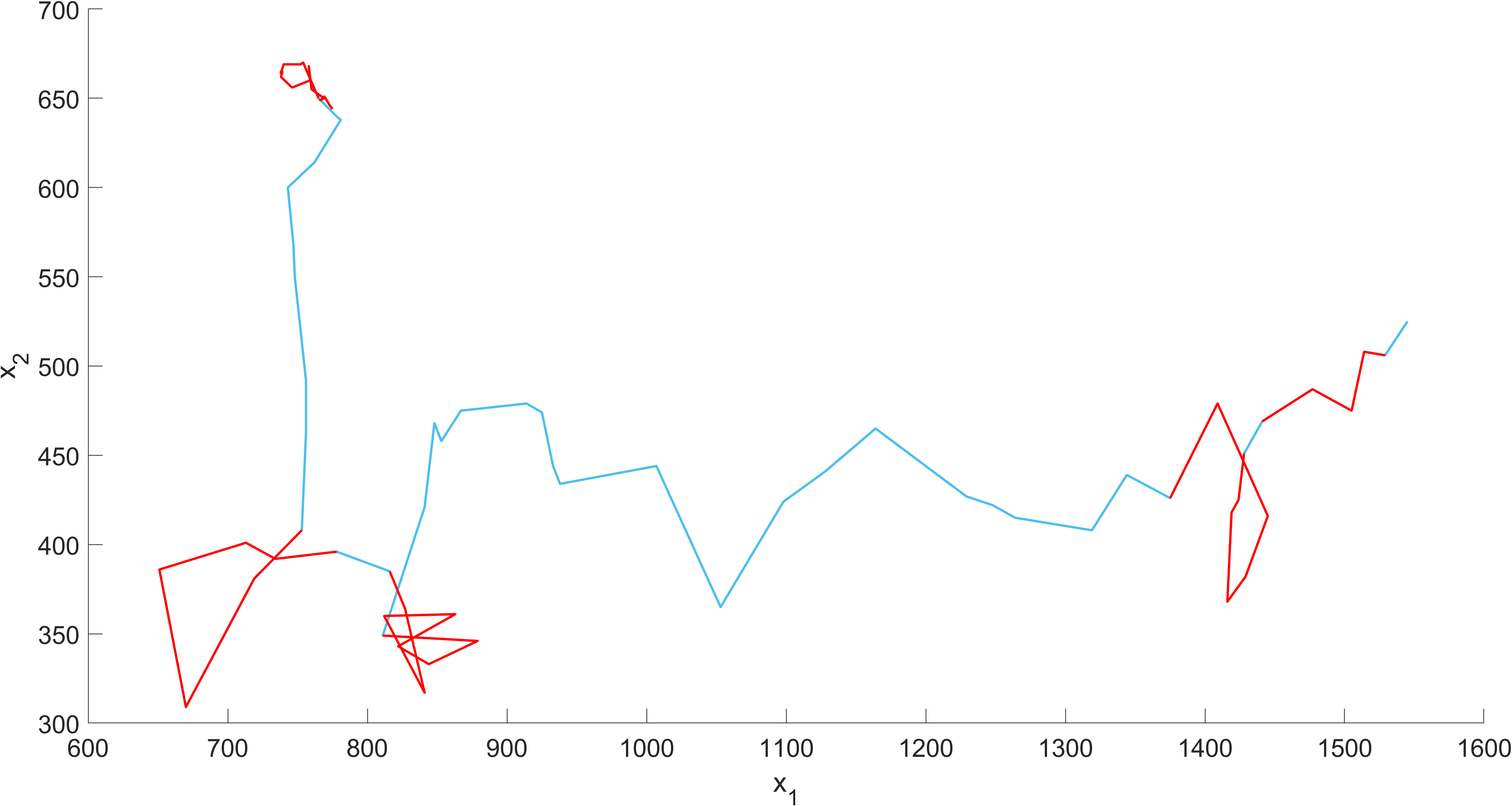}
    \includegraphics[width=8cm]{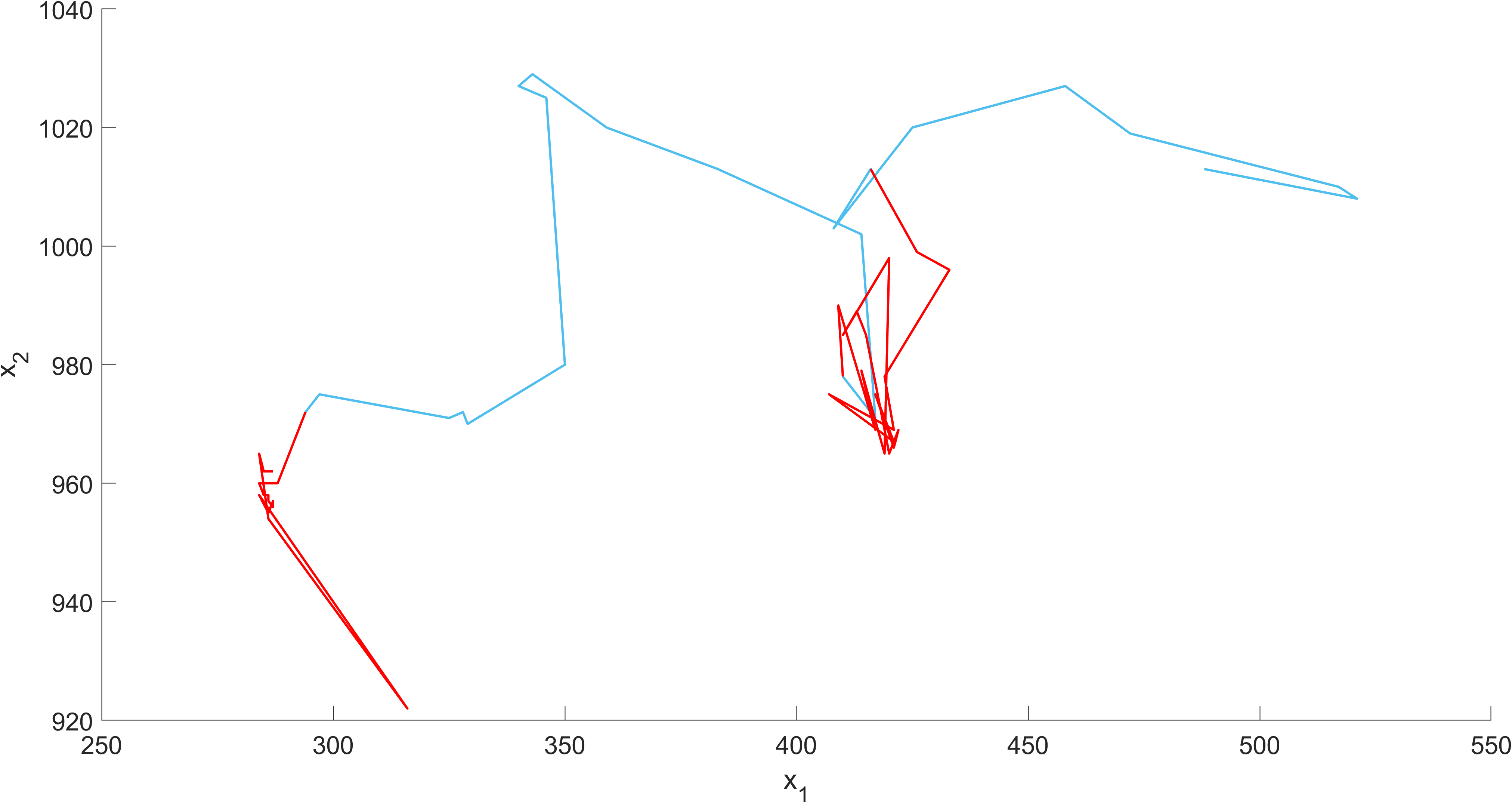}
    \caption{Extraction of random parts by evolving curves: original trajectory (light blue line), and extracted random parts (red lines) by evolving curves.}
    \label{F: RandomParts evolving curve}
\end{figure}
A random sub-trajectory is obtained when more consecutive segments disappear. To capture the random motion for all the points in the disappeared segments, we included also the point immediately preceding the start and the point immediately following the end of the random sub-trajectory. Since our goal is to analyze the diffusive behavior and do a statistical analysis of the extracted parts, we need a minimum number of points for each random sub-trajectory; therefore, we chose to consider only trajectories with at least $5$ points. \\
We considered $3$ different datasets as described in Section \ref{SS: Data acquisition}, extracted the random parts, and analyzed them separately. For all the datasets we chose the same parameters for the trajectory smoothing as follows: $\delta=0.003$, $\lambda=20$, $\omega=50$, and $\tau=0.000001$. The parameters $\delta$ and $\lambda$ were constant and, when no more self-intersections were detected, the algorithm performed $50$ more iterations and stopped. Fig. \ref{F: RandomParts evolving curve} shows the result of the random parts detection for the trajectories shown in Figs. \ref{F:T2}, \ref{F:T63}: the light blue lines are the original trajectories, and the red parts are the extracted random parts. \\
\begin{figure}[htbp]
    \centering
    \subfigure[]{\includegraphics[width=6.5cm]{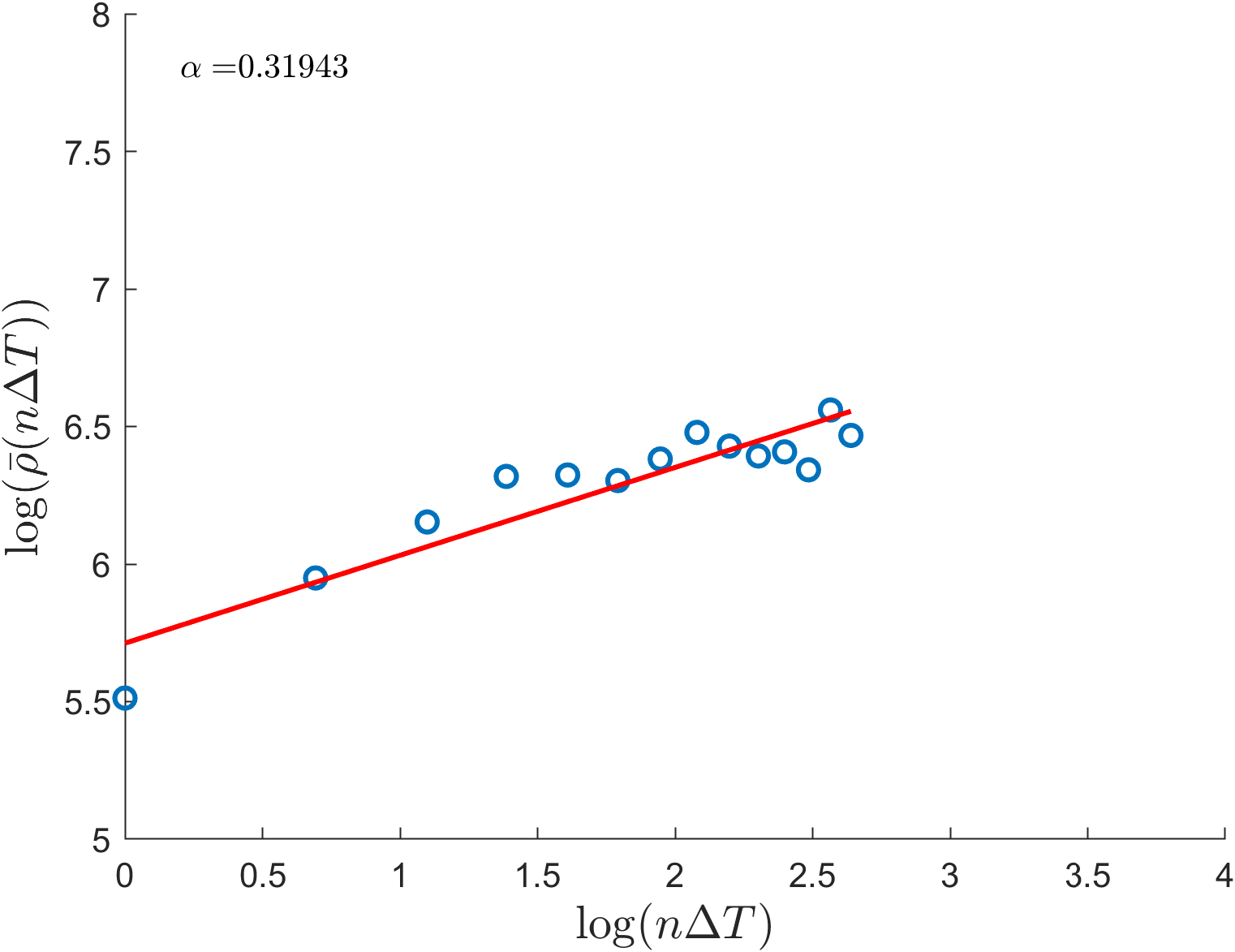} \label{F: MSD EC Resul3}}
    \subfigure[]{\includegraphics[width=6.5cm]{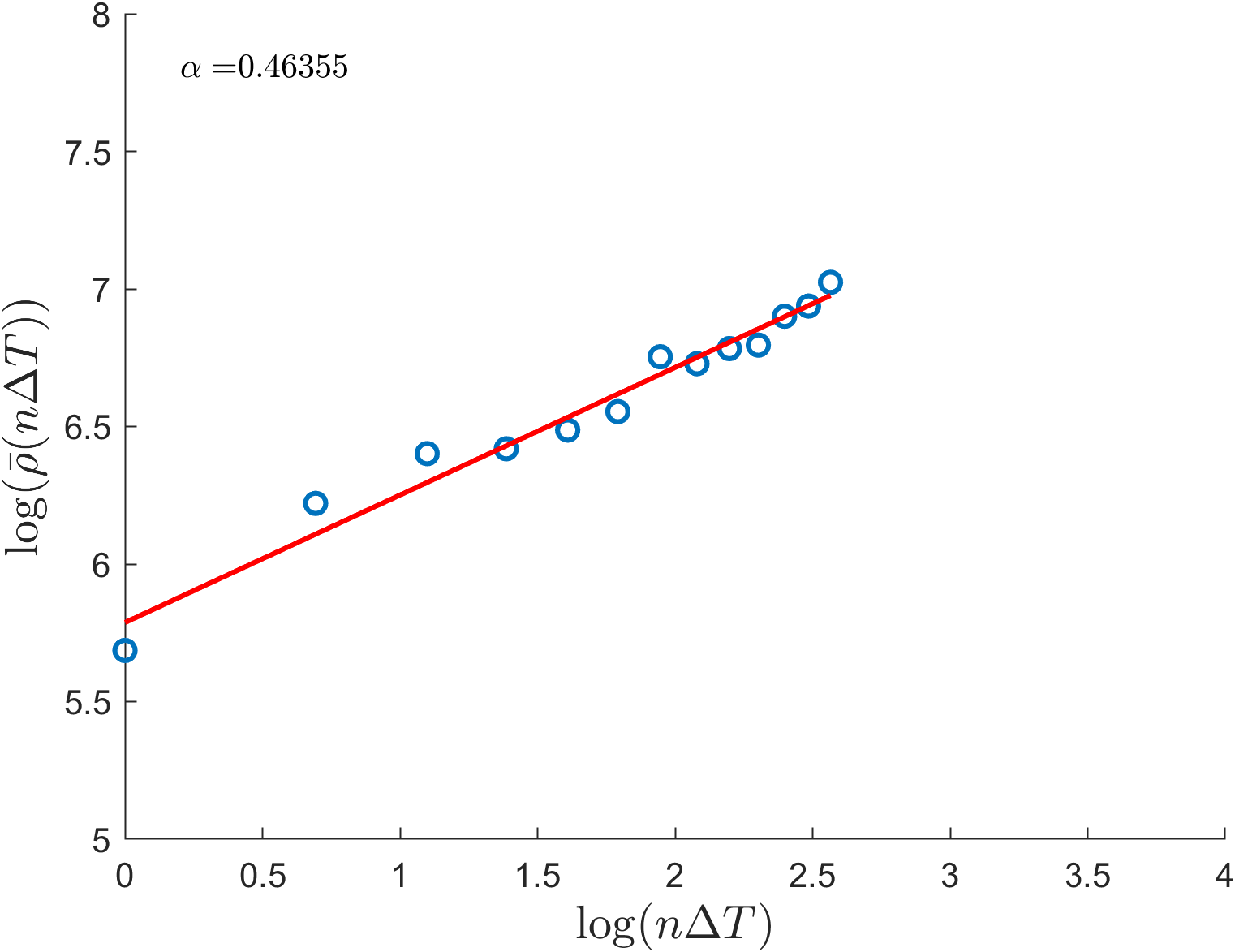}\label{F: MSD EC Resul4}}
    \subfigure[]{\includegraphics[width=6.5cm]{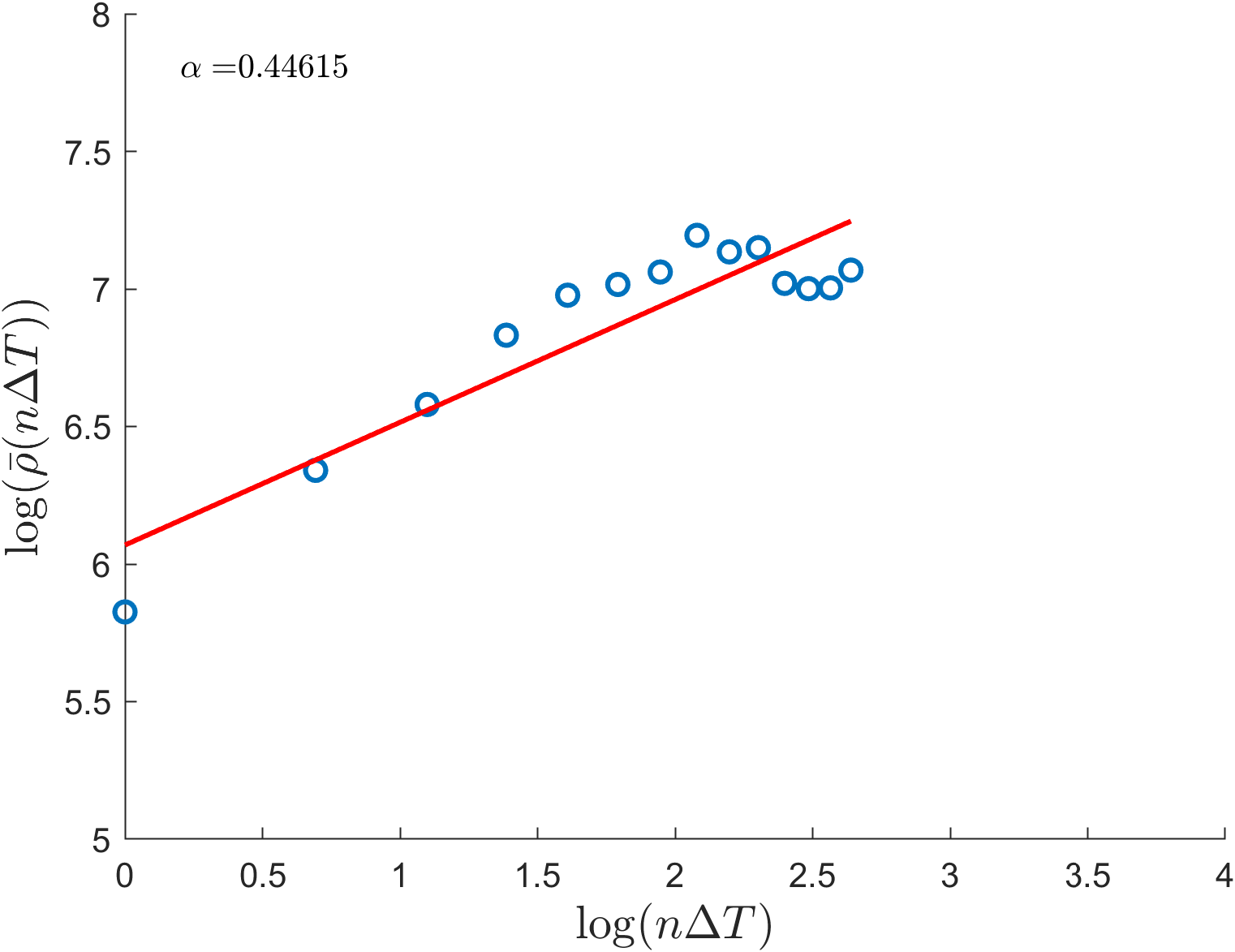}\label{F: MSD EC Resul5}}
    \caption{log-log plot of $\bar{\rho}(n\Delta T)$ versus $n\Delta T$. Results from the $3$ datasets for random parts extracted by evolving curves. Blue dots represent the values of EATAMSD, the red line is the result of the linear regression; $\alpha$ is the slope of the red line.}
    \label{F: MSD for all Evolving Curves}
\end{figure}
 We then applied the formula (\ref{E: MSD time all}) to the $3$ different datasets, and calculated the values of $\alpha$. We considered the log-log dependence on time of the EATAMSD and computed linear regression using the \textit{polyfit} function in MATLAB. Fig. \ref{F: MSD for all Evolving Curves} shows the results for the different datasets: the blue dots represent the values of the EATAMSD, while the red line is the result of the linear regression. The extracted random parts were $181$ for the dataset in Fig. \ref{F: MSD EC Resul3}, $105$ for the one in Fig. \ref{F: MSD EC Resul4}, and $77$ for the one in Fig. \ref{F: MSD EC Resul5}.
 As one can see, the EATAMSD shows subdiffusive motion for the $3$ datasets analyzed. 
 
\subsection{Random parts extracted from curve self-intersections}
The second approach to finding the random parts considers the curve self-intersections. Indeed, the random segments are in areas of high curvature where the cell does a back-and-forth motion causing the random segments to intersect. Thanks to the algorithm described in Section \ref{SS: Detect self intersections}, we can detect the two grid points where the intersection starts and ends. We define as random segments the segments to which the two points belong and all the segments in between. If the grid point where the start of the self-intersection is detected belongs to the $2$ segment, i.e. is the endpoint of those segments, we consider as initial segment the one in the direction of parametrization. Similarly, if the grid point where the end of the self-intersection is detected belongs to the $2$ segment, we consider as the last segment of the random sub-trajectory the segment in the direction opposite to parametrization.
 We then obtain random sub-trajectories; as in Section \ref{SS: Random parts extracted by evolving curves}, we consider only trajectories with at least $5$ points.
\begin{figure}[htbp]
    \centering
    \includegraphics[width=8cm]{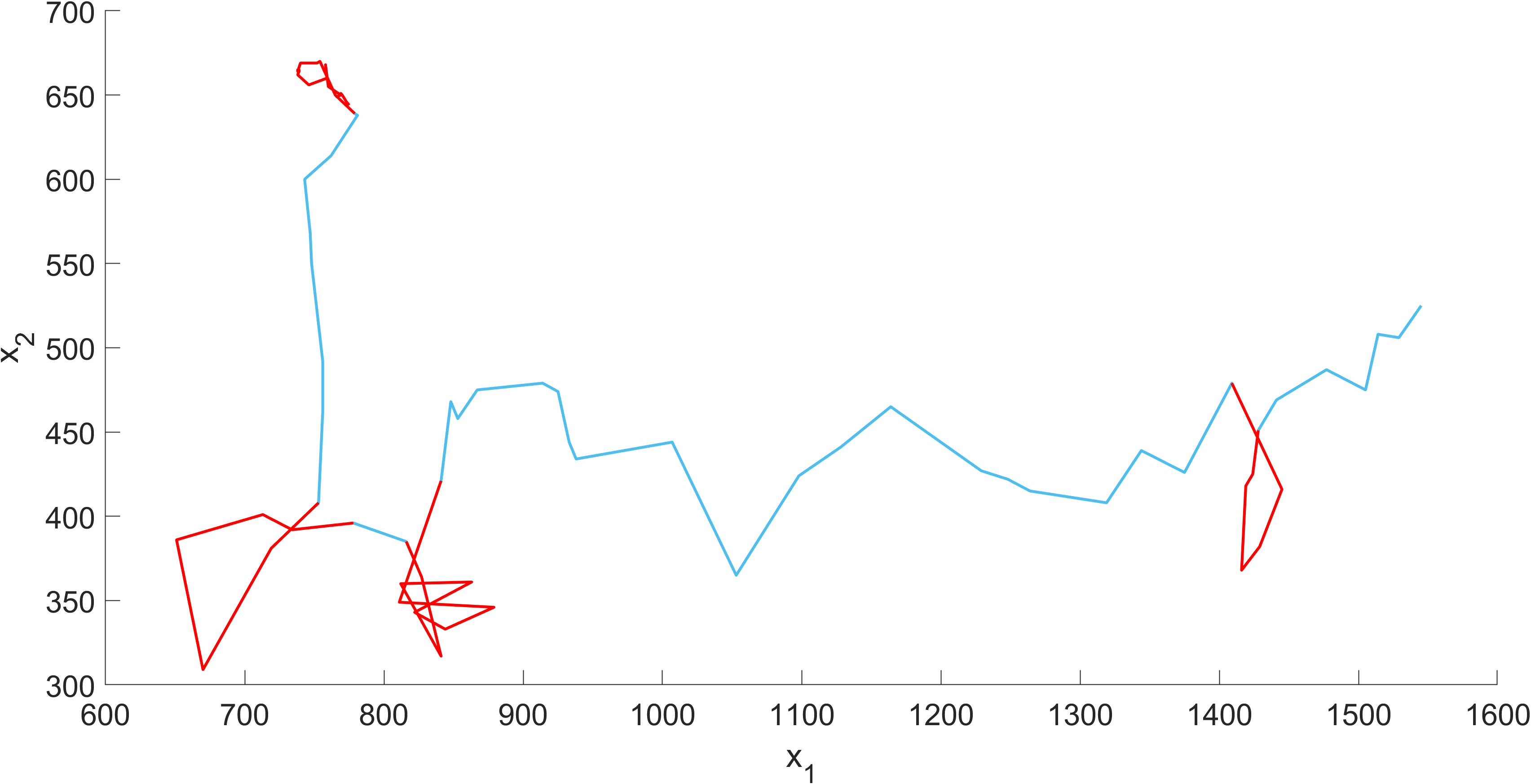}
    \includegraphics[width=8cm]{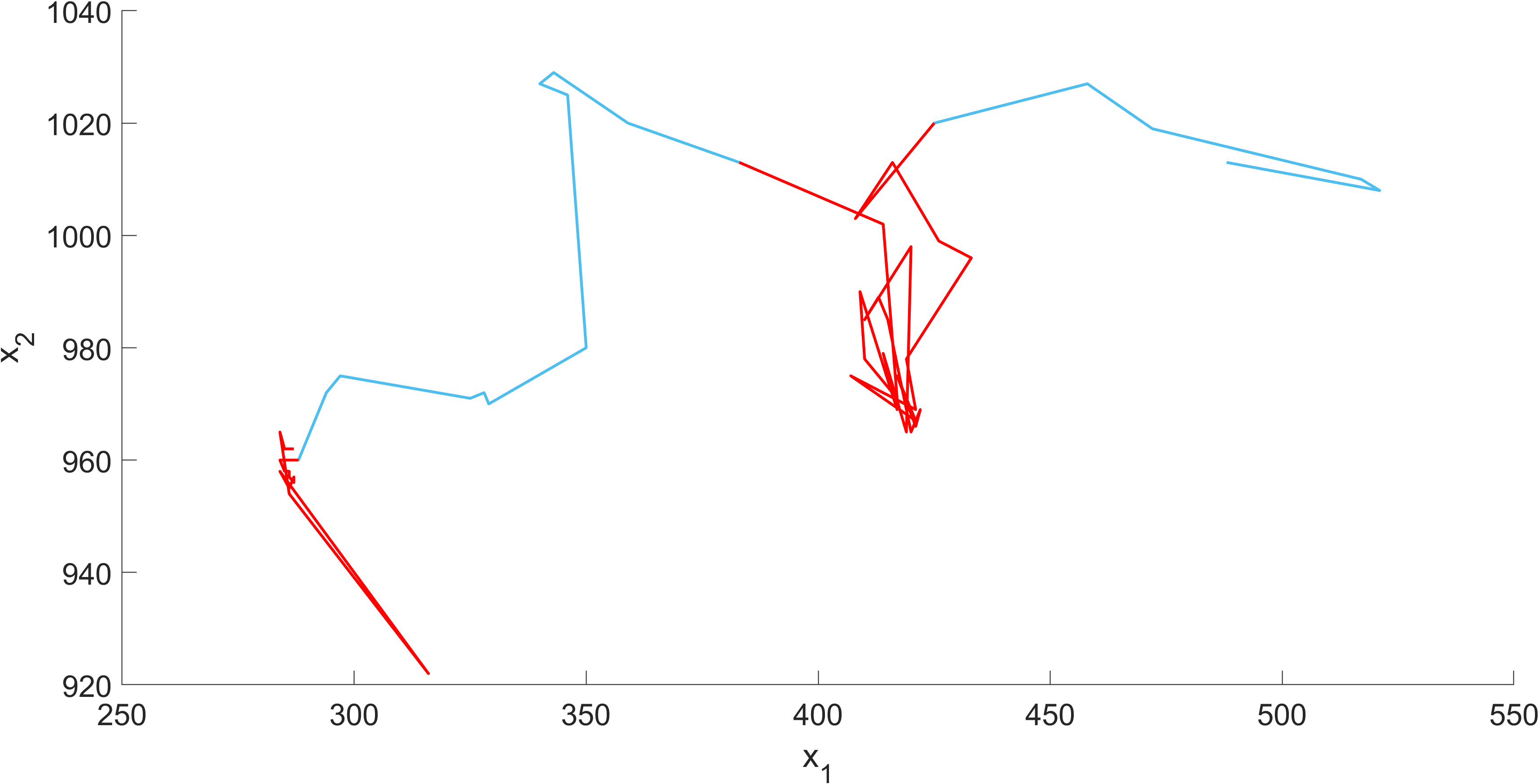}
    \caption{Extraction of random parts by evolving curves: original trajectory (light blue line), and extracted random parts (red lines) from curve self-intersections.}
    \label{F: RandomParts self intersections}
\end{figure}
The detection of the random parts, in this case, can be done directly from the original curves and does not depend on the smoothing parameters. Fig. \ref{F: RandomParts self intersections} shows the result for the trajectories shown in Figs. \ref{F:T2 Mix}, \ref{F:T63 Mix}: the light blue line is the original trajectory, and the red parts are the extracted random parts. \\ 
\begin{figure}[htbp]
    \centering
    \subfigure[]{\includegraphics[width=6.2cm]{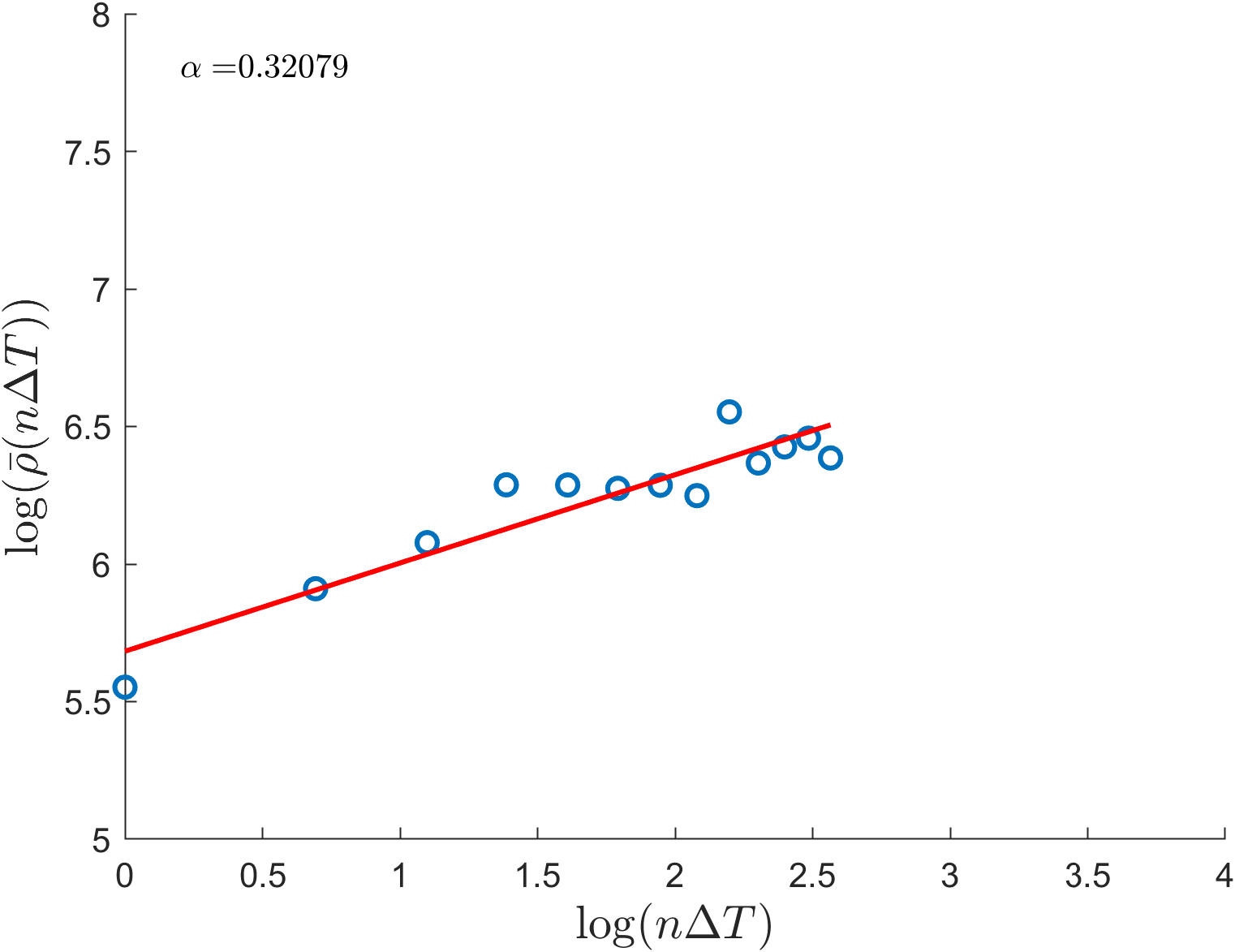}\label{F: MSD INT Resul3}}
    \subfigure[]{\includegraphics[width=6.2cm]{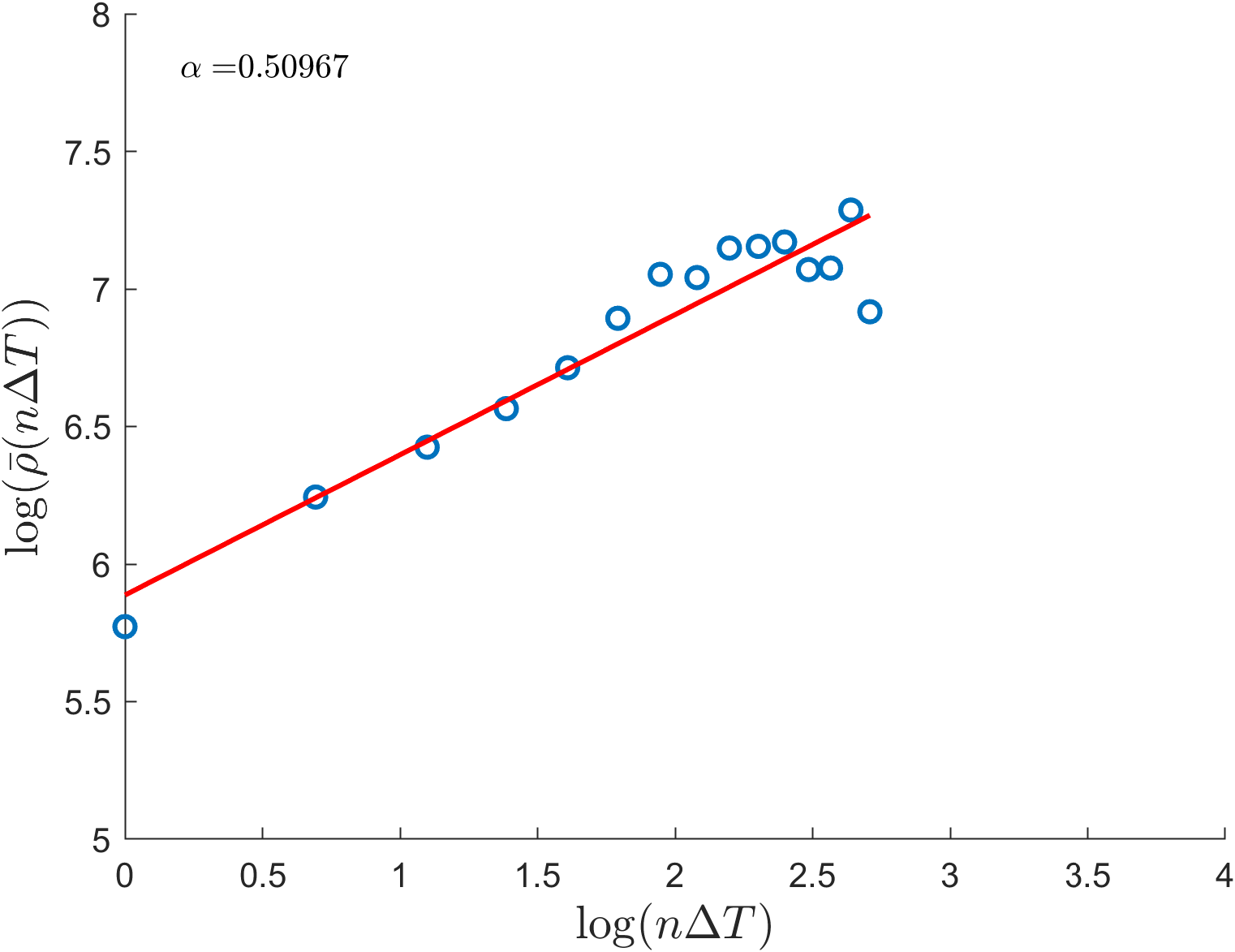}\label{F: MSD INT Resul4}}
    \subfigure[]{\includegraphics[width=6.2cm]{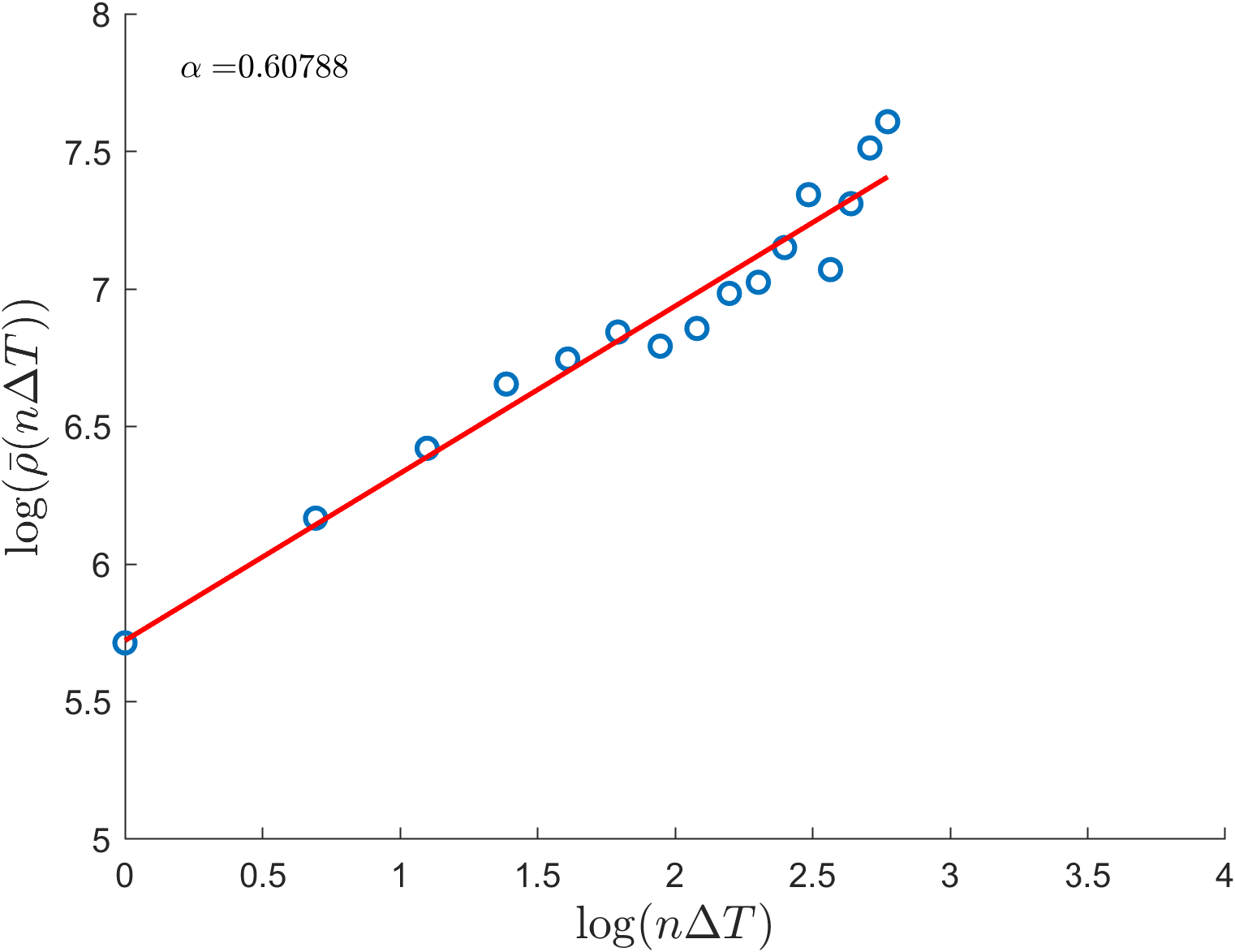}\label{F: MSD INT Resul5}}
    \caption{log-log plot of $\bar{\rho}(n\Delta T)$ versus $n\Delta T$. Results from the $3$ datasets for random parts extracted from curve self-intersections. Blue dots represent the values of EATAMSD, the red line is the result of the linear regression; $\alpha$ is the slope of the red line.}
    \label{F: MSD for all Intersections}
\end{figure}
We then applied the formula (\ref{E: MSD time all}) to the $3$ different datasets, and calculated the values of $\alpha$ as in Section \ref{SS: Random parts extracted by evolving curves}. Fig. \ref{F: MSD for all Intersections} shows the results for the different datasets: the blue dots represent the values of the EATAMSD, while the red line is the result of the linear regression, and $\alpha$ is the slope of the line. In this case, the extracted random parts were $205$ for the dataset in Fig. \ref{F: MSD INT Resul3}, $107$ for the one in Fig. \ref{F: MSD INT Resul4}, and $83$ for the one in Fig. \ref{F: MSD INT Resul5}. As one can see, also with this approach the EATAMSD shows subdiffusive motion for the $3$ datasets analyzed. 
\begin{table}
\centering
\begin{tabular}{ | m{4cm} | m{3cm}| m{3cm} | } 
  \hline
  Dataset& $\alpha$ for disappearing segments approach & $\alpha$ for \newline 
  self-intersections \newline approach \\ 
  \hline
  Dataset 220224$\_$05hpa$\_$6hpa & 0.31943 & 0.32079 \\ 
  \hline
  Dataset 221006$\_$05hpa$\_$6hpa & 0.46355 & 0.50967 \\ 
  \hline
  Dataset 230225$\_$1hpa$\_$6hpa & 0.44615 & 0.60788 \\ 
  \hline
\end{tabular}
 \vspace{0.4cm}
\caption{Slope $\alpha$ of the line in the log-log scale for $3$ different datasets. The first column shows the results for the disappearing segments approach, while the second column shows the results for the self-intersection approach.}
\label{T: Table alpha}
\end{table}
Finally, Table \ref{T: Table alpha} illustrates the results of the slopes $\alpha$ for the 3 datasets and the 2 different approaches. As one can notice, the values for the different approaches are similar, with the maximum difference in the slope values between the two approaches being $0.16173$ for the last dataset.

\section{Vector field reconstruction}
\label{S: Vector Field}
The mathematical model we propose to reconstruct the $2$D vector field is based on interpolating/extrapolating the given information about the vector field. In particular, we use the smoothed velocities on the smoothed trajectories as sparse samples to reconstruct the wound attractant field. Initially, we consider the image domain of the 2D microscopy video of macrophages moving toward the wound. Then, we define the model domain by cutting off from the image domain some pixels, i.e. squares, where the vector field is $0$ respectively where the sparse samples are located.

\begin{figure}[htbp]
    \centering
    \subfigure[]{\includegraphics[width=0.45\textwidth]{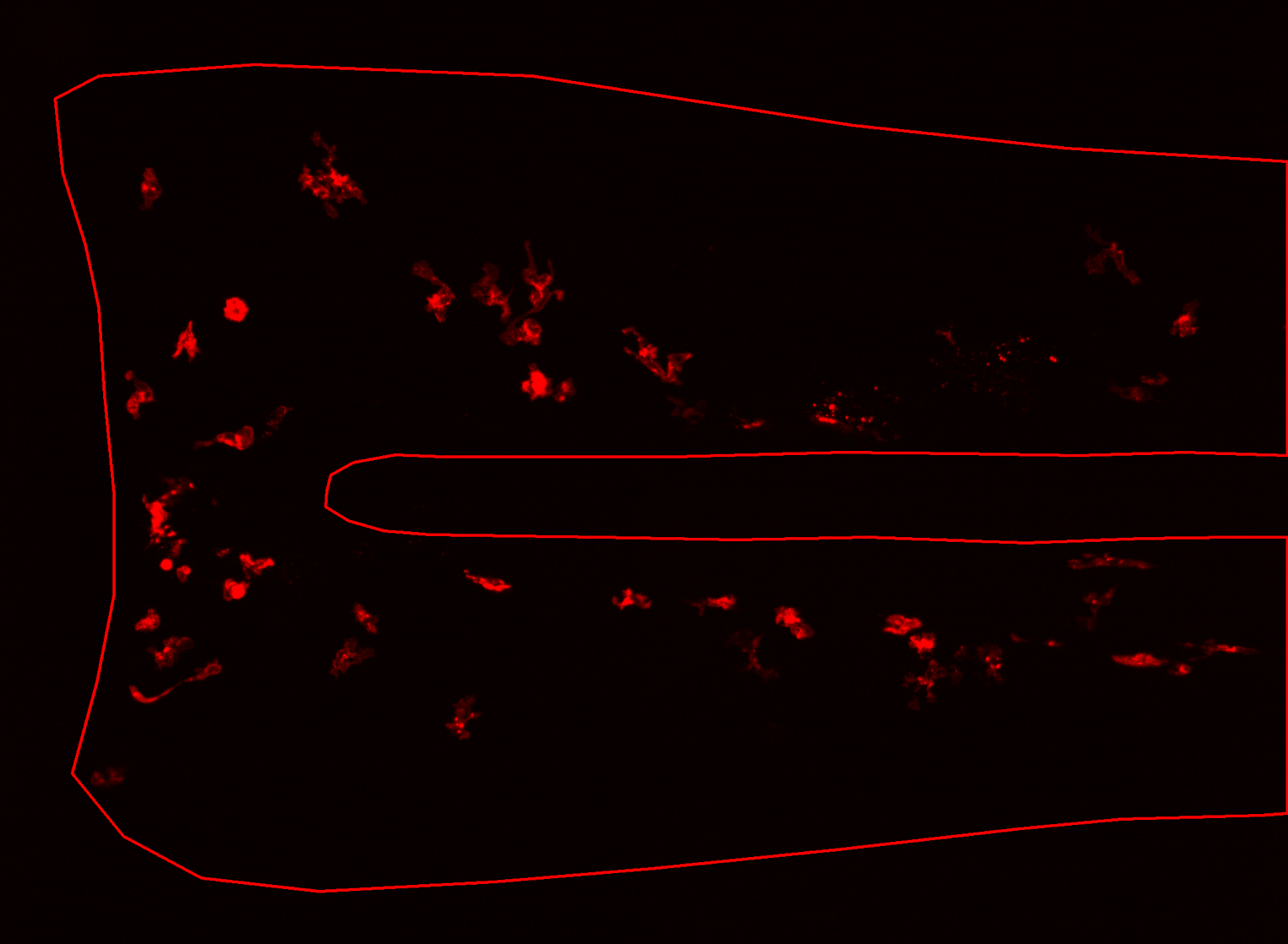} \label{F: Resul3A_THEORY} }
    \subfigure[]{\includegraphics[width=0.45\textwidth]{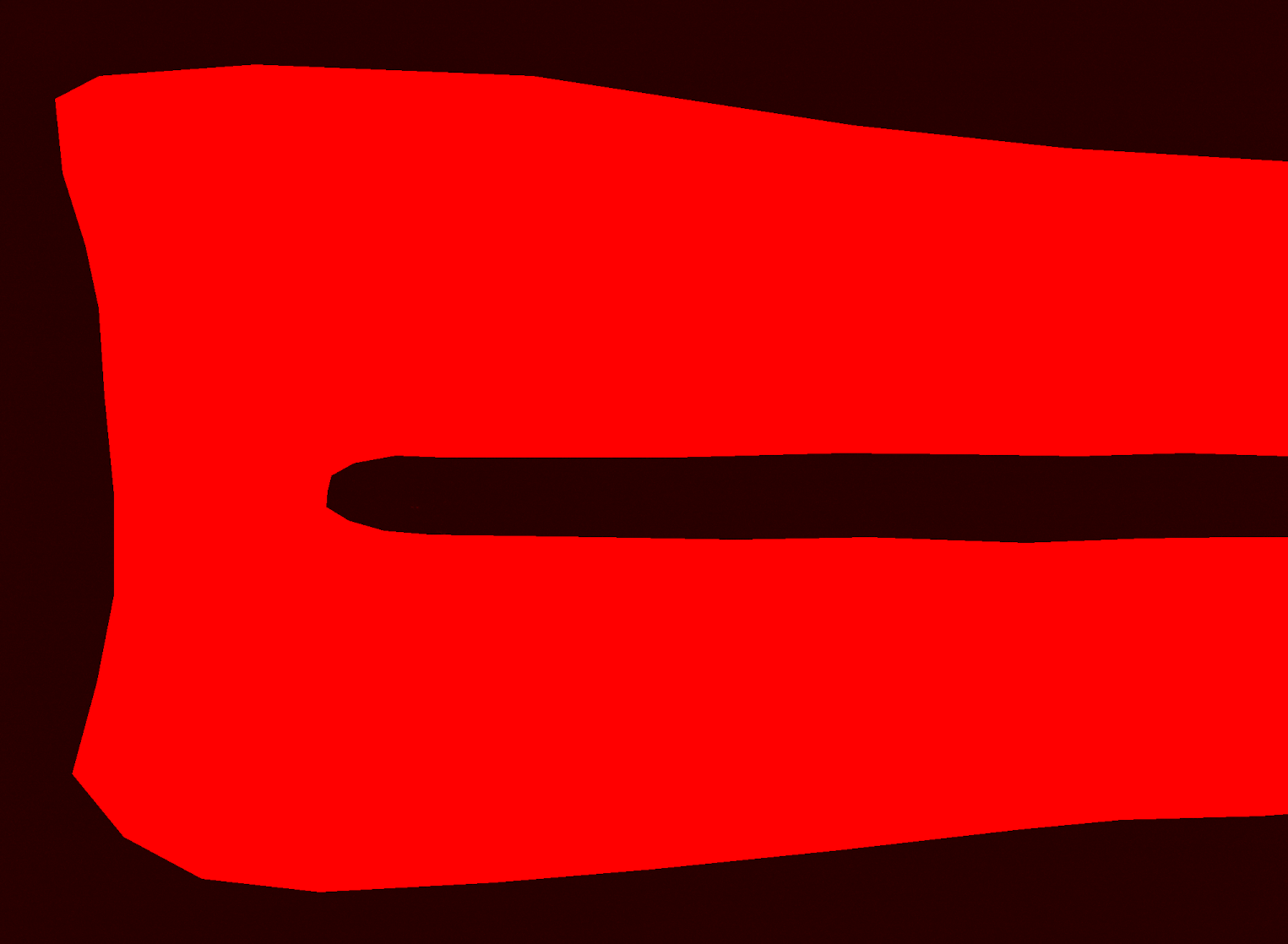}\label{F: Resul3B_THEORY}}
    \caption{\subref{F: Resul3A_THEORY} Black area indicates the image domain. The red curve indicates the boundary of $\tilde{\Omega}_{cut}$. \subref{F: Resul3B_THEORY} $\tilde{\Omega}_{cut}$, region of the image domain where the wound attractant field is reconstructed.}
    \label{F: Resul3_theory}
\end{figure}

\begin{figure}[htbp]
    \centering
    \subfigure[]{\includegraphics[width=0.45\textwidth]{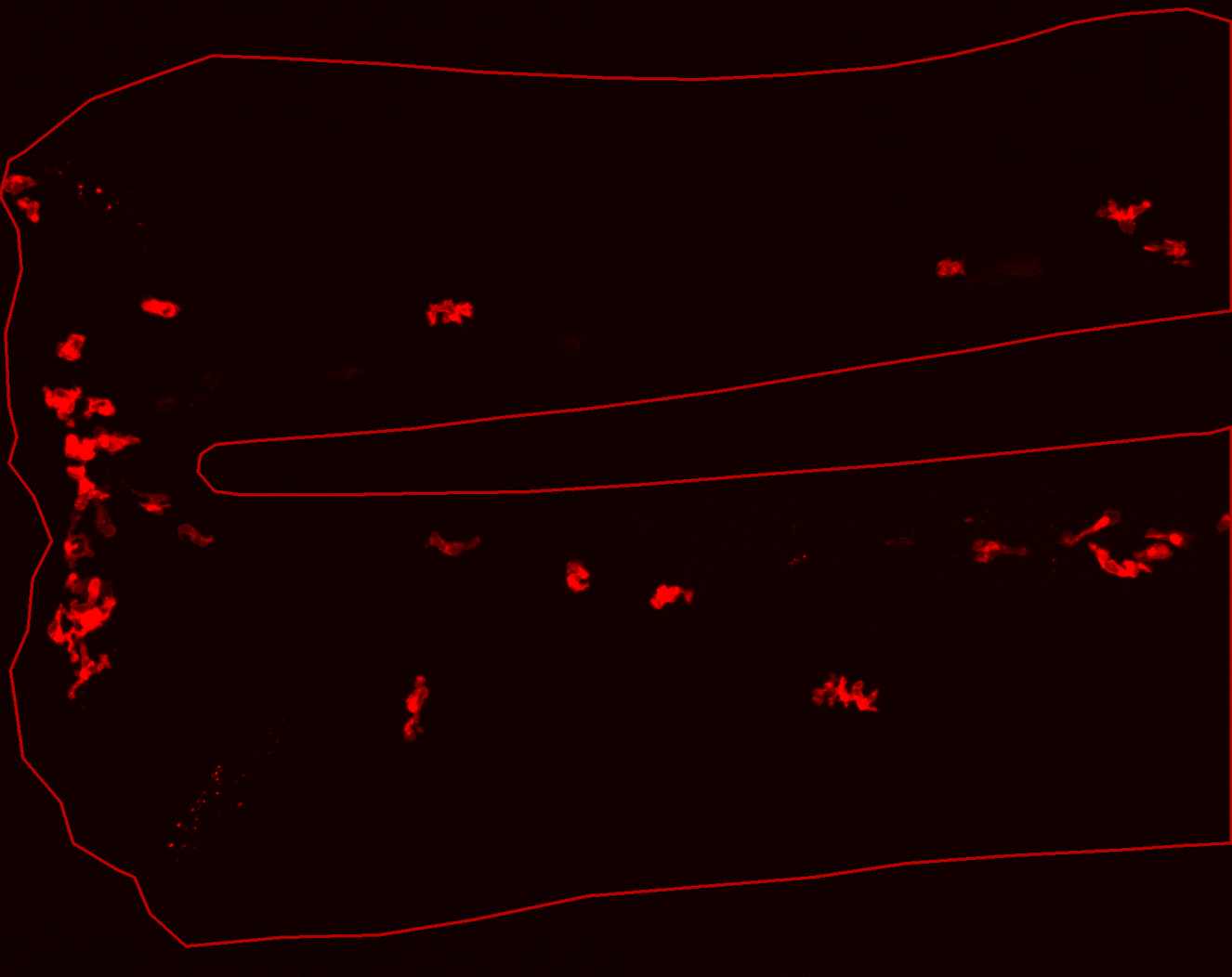} \label{F: Resul4A_THEORY} }
    \subfigure[]{\includegraphics[width=0.45\textwidth]{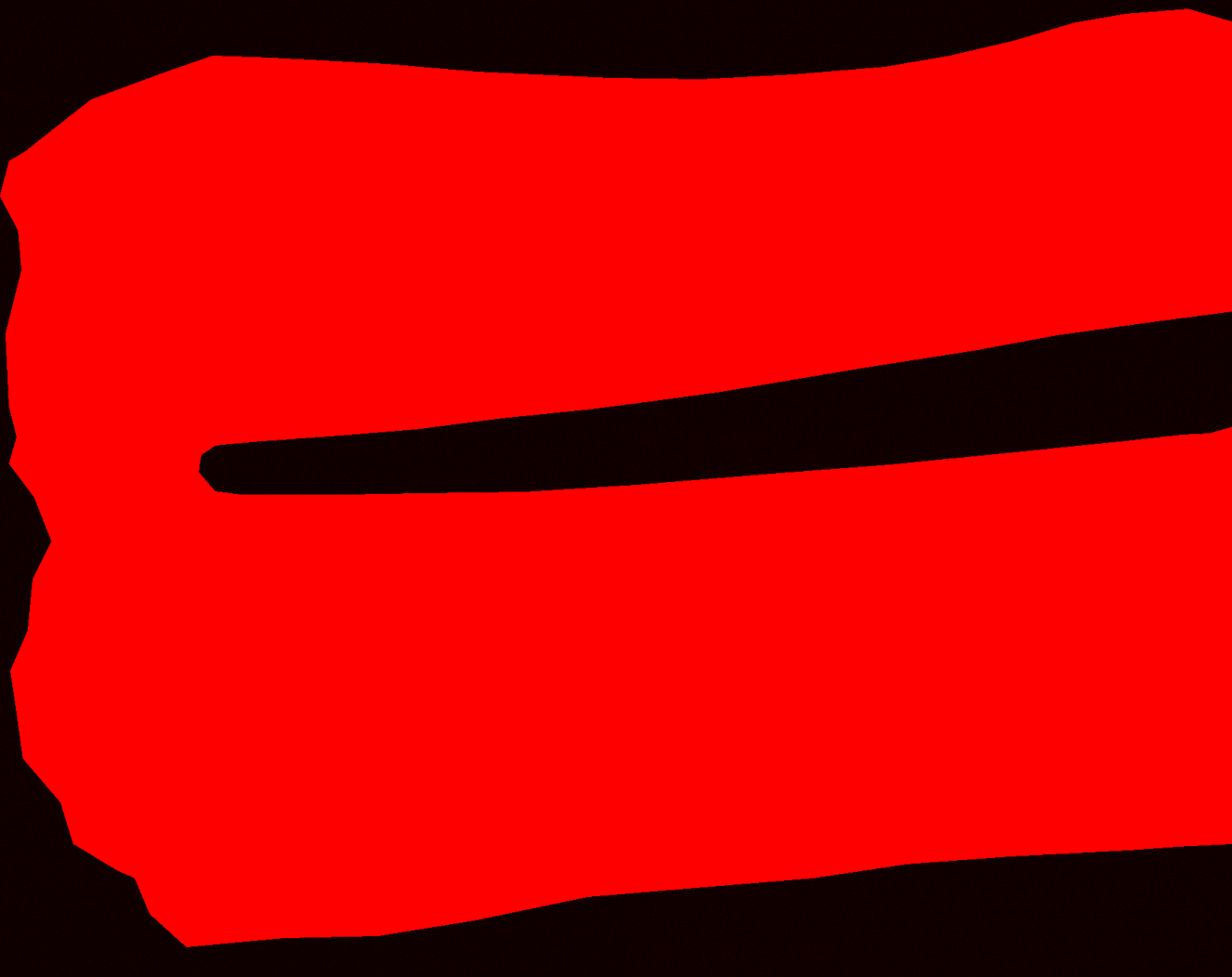}\label{F: Resul4B_THEORY}}
    \caption{\subref{F: Resul4A_THEORY} Black area indicates the image domain. The red curve indicates the boundary of $\tilde{\Omega}_{cut}$. \subref{F: Resul4B_THEORY} $\tilde{\Omega}_{cut}$, region of the image domain where the wound attractant field is reconstructed.}
    \label{F: Resul4_theory}
\end{figure}

\begin{figure}[htbp]
    \centering
    \subfigure[]{\includegraphics[width=0.45\textwidth]{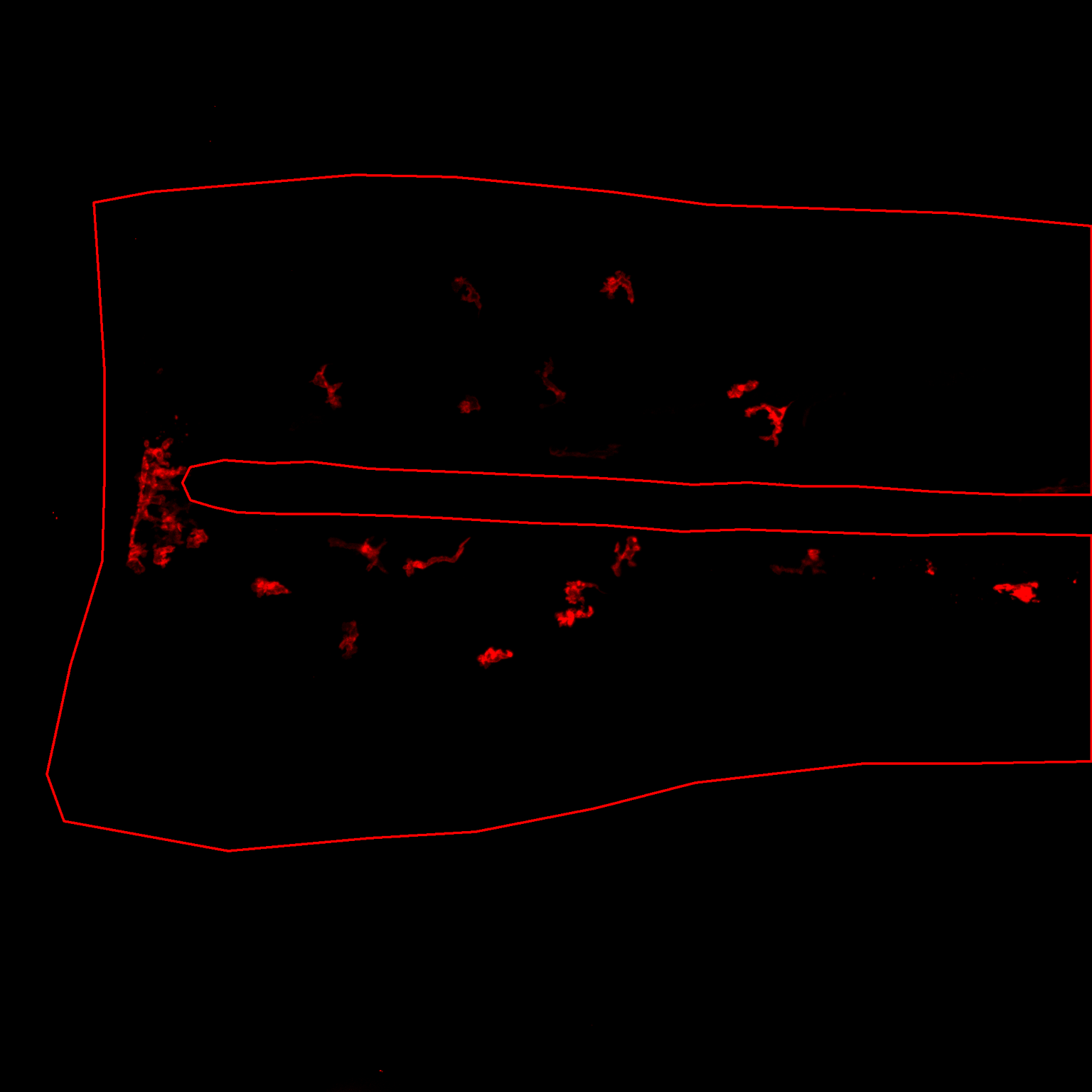} \label{F: Resul5A_THEORY} }
    \subfigure[]{\includegraphics[width=0.45\textwidth]{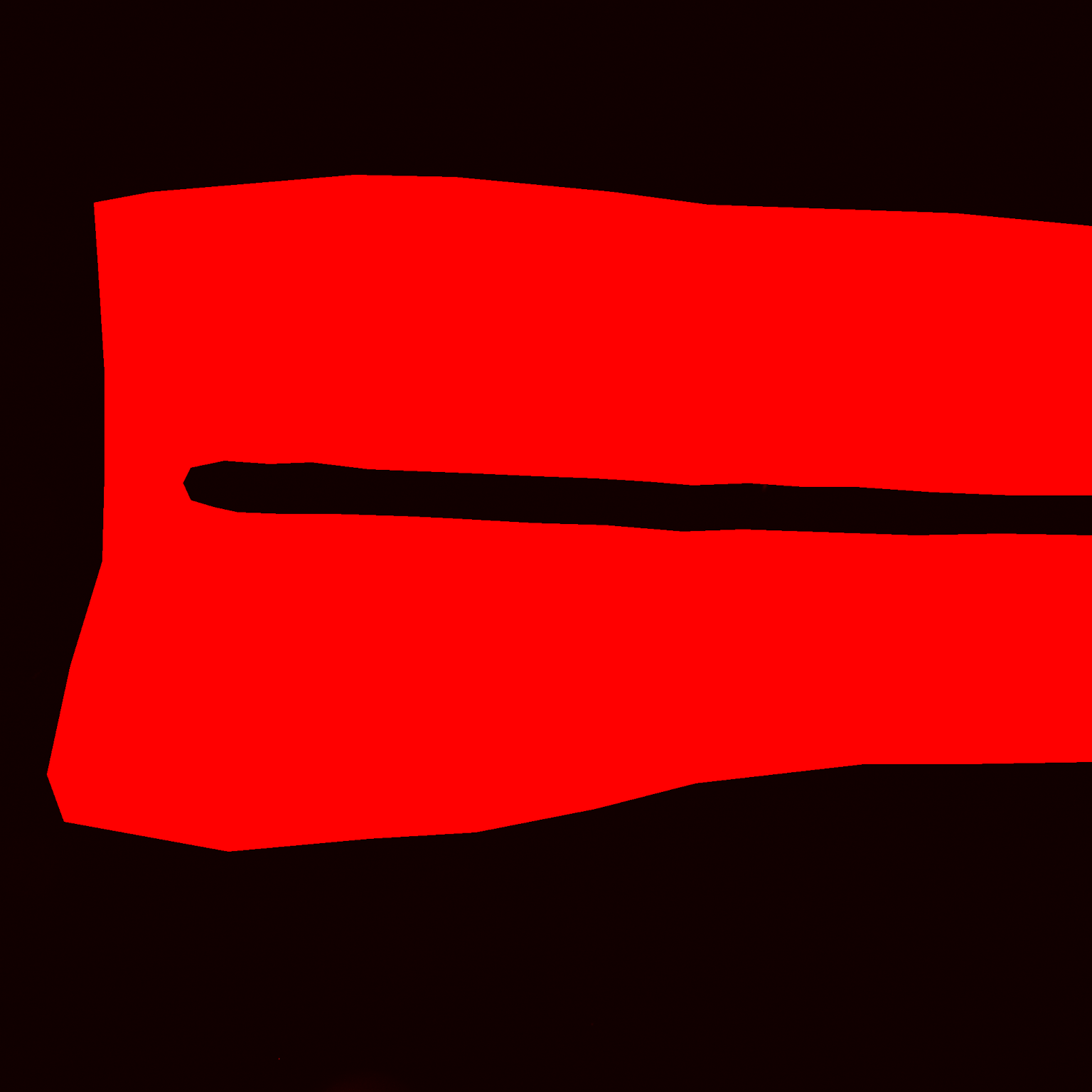}\label{F: Resul5B_THEORY}}
    \caption{\subref{F: Resul5A_THEORY} Black area indicates the image domain. The red curve indicates the boundary of $\tilde{\Omega}_{cut}$. \subref{F: Resul5B_THEORY} $\tilde{\Omega}_{cut}$, region of the image domain where the wound attractant field is reconstructed.}
    \label{F: Resul5_theory}
\end{figure}

Consider $\tilde{\Omega}$ the auxiliary domain, i.e. the image domain. In our specific application, we know that there is a part of the domain where the vector field is $0$ because either is outside of the fish body or because the macrophages can not enter that part of the fish body. 
Therefore, we first consider the cut domain $\tilde{\Omega}_{cut}$ as indicated in Figs. \ref{F: Resul3B_THEORY}, \ref{F: Resul4B_THEORY}, \ref{F: Resul5B_THEORY}. In the figures, the black area represents the image domain $\tilde{\Omega}$ while the red area represents $\tilde{\Omega}_{cut}$.
Let us define
 \begin{equation}
    \Gamma^2=\partial \tilde{\Omega}_{cut}
 \end{equation} 
the domain boundary, cf. Figs. \ref{F: Resul3A_THEORY}, \ref{F: Resul4A_THEORY}, \ref{F: Resul5A_THEORY}. 
By the construction, $\partial \tilde{\Omega}_{cut}$ is Lipschitz.  
\begin{figure}[htbp]
    \centering
    \includegraphics[width= 0.6\linewidth]{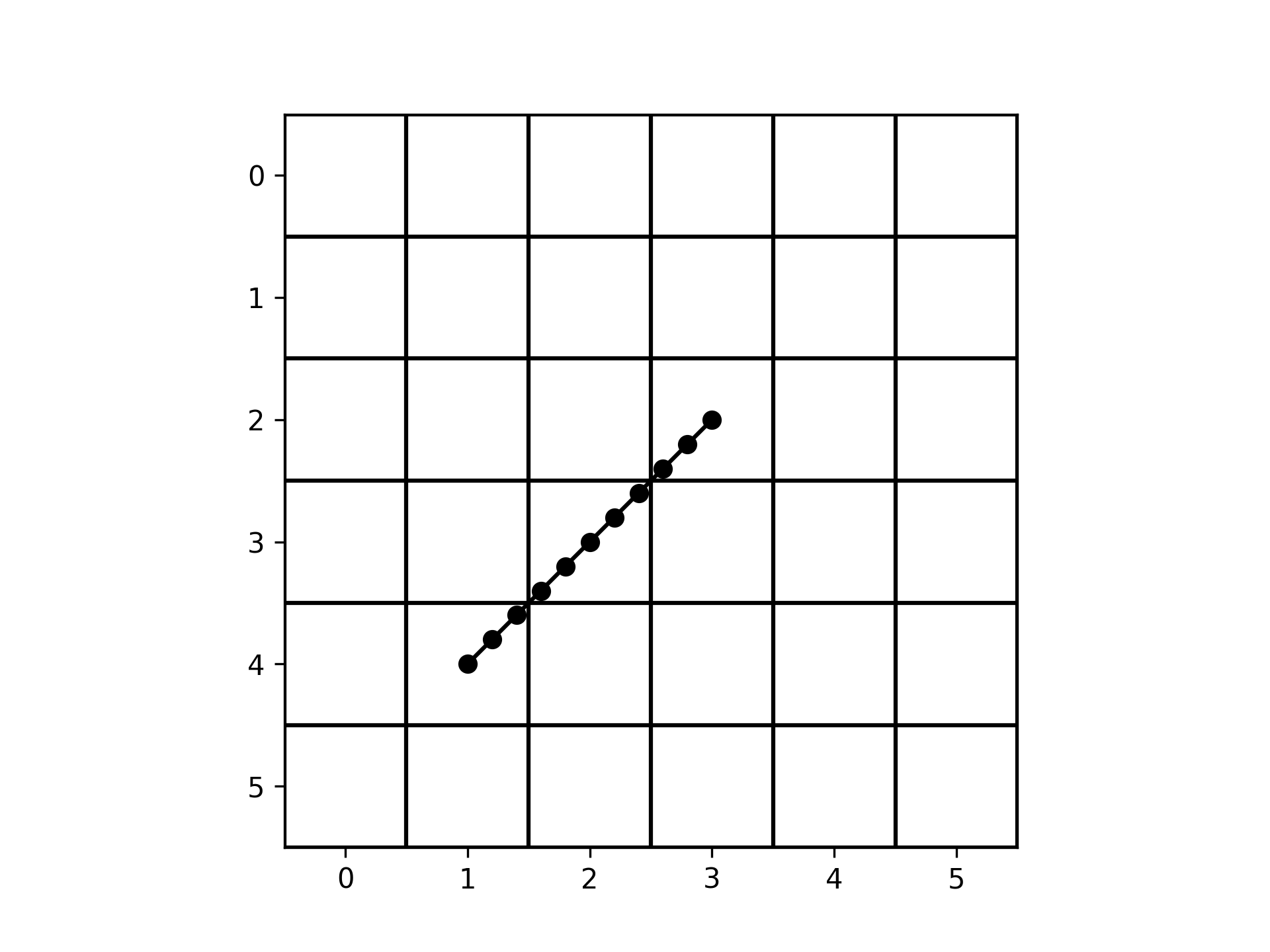} 
    \caption{Example of trajectory (black line) in a pixel grid. The black dots represent the trajectory's grid points.}
    \label{F: VF1}
\end{figure}

\begin{figure}[htbp]
    \centering
    \includegraphics[width= 0.6\linewidth]{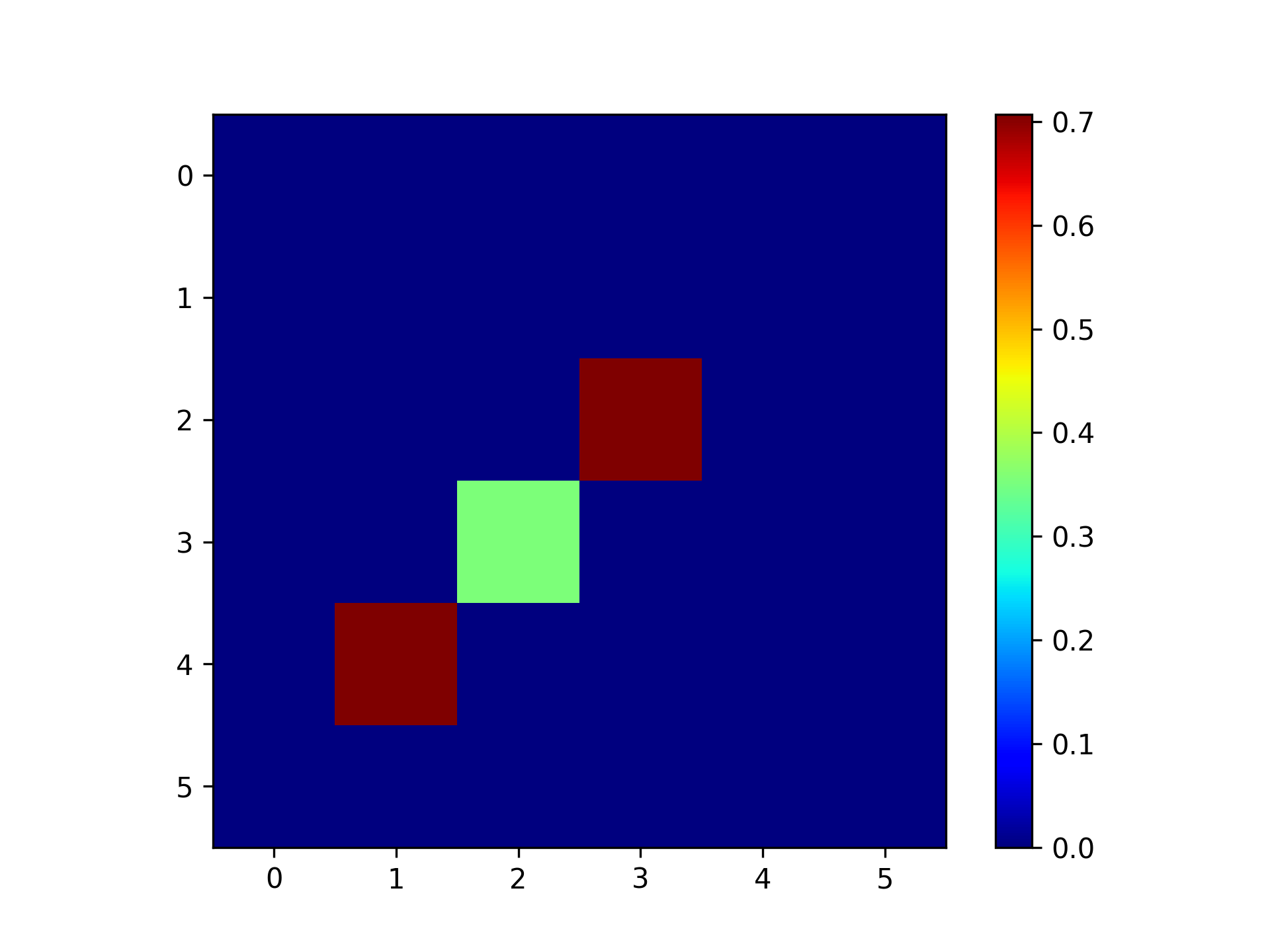} 
    \caption{Visualization of the values of the squares in Fig. \ref{F: VF1}: the value is the mean of the values of the grid points belonging to that square. The other squares have value $0$ for visualization.}
    \label{F: VF2}
\end{figure}

\begin{figure}[htbp]
    \centering
    \includegraphics[width= 0.6\linewidth]{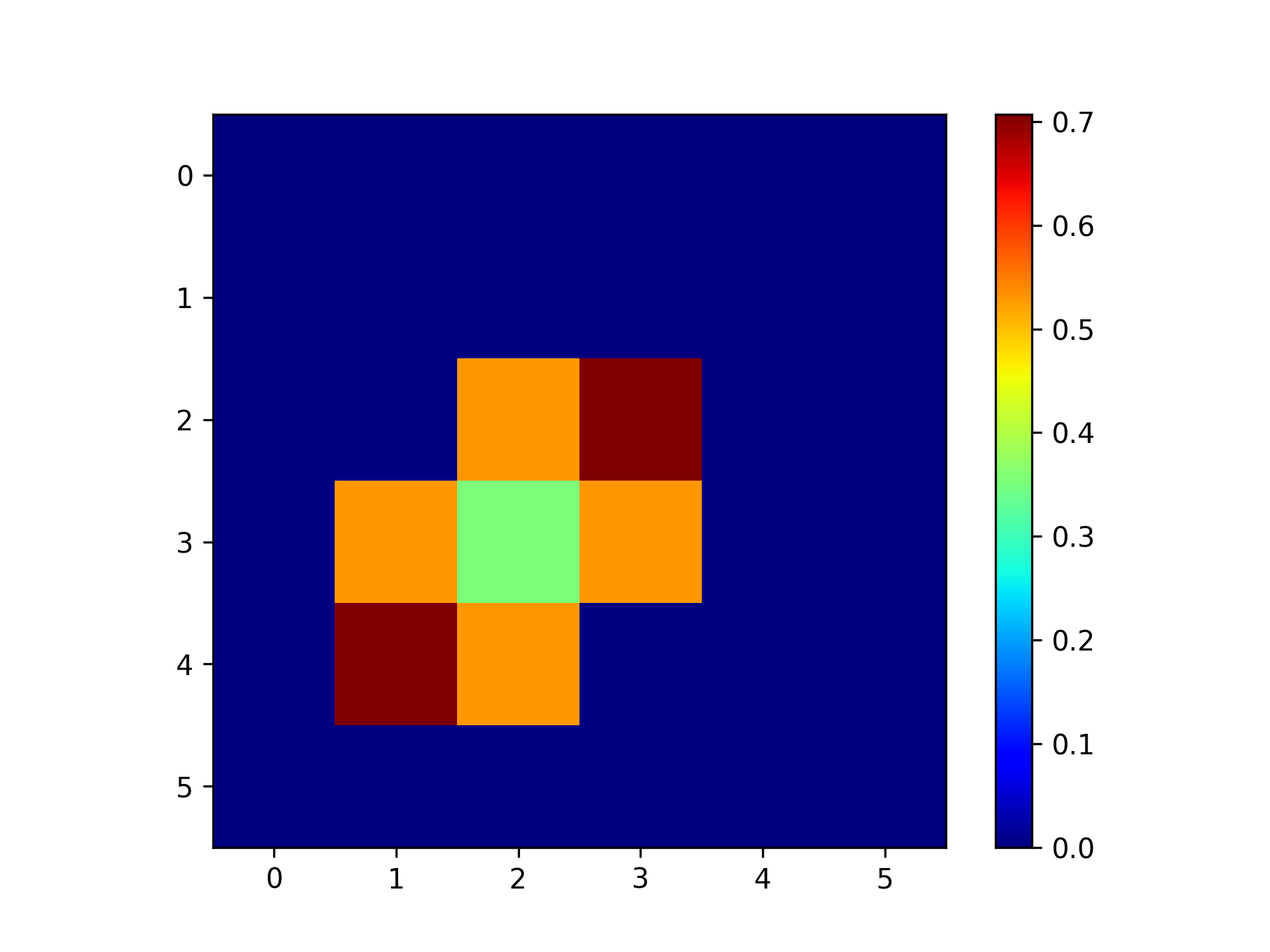} 
    \caption{Lipschitz domain associated to the domain in Fig. \ref{F: VF2}. The squares that have in common the vertex are considered and their value is set as the mean of the other two squares.}
    \label{F: VF3}
\end{figure}

\begin{figure}[htbp]
    \centering
    \includegraphics[width= 0.6\linewidth]{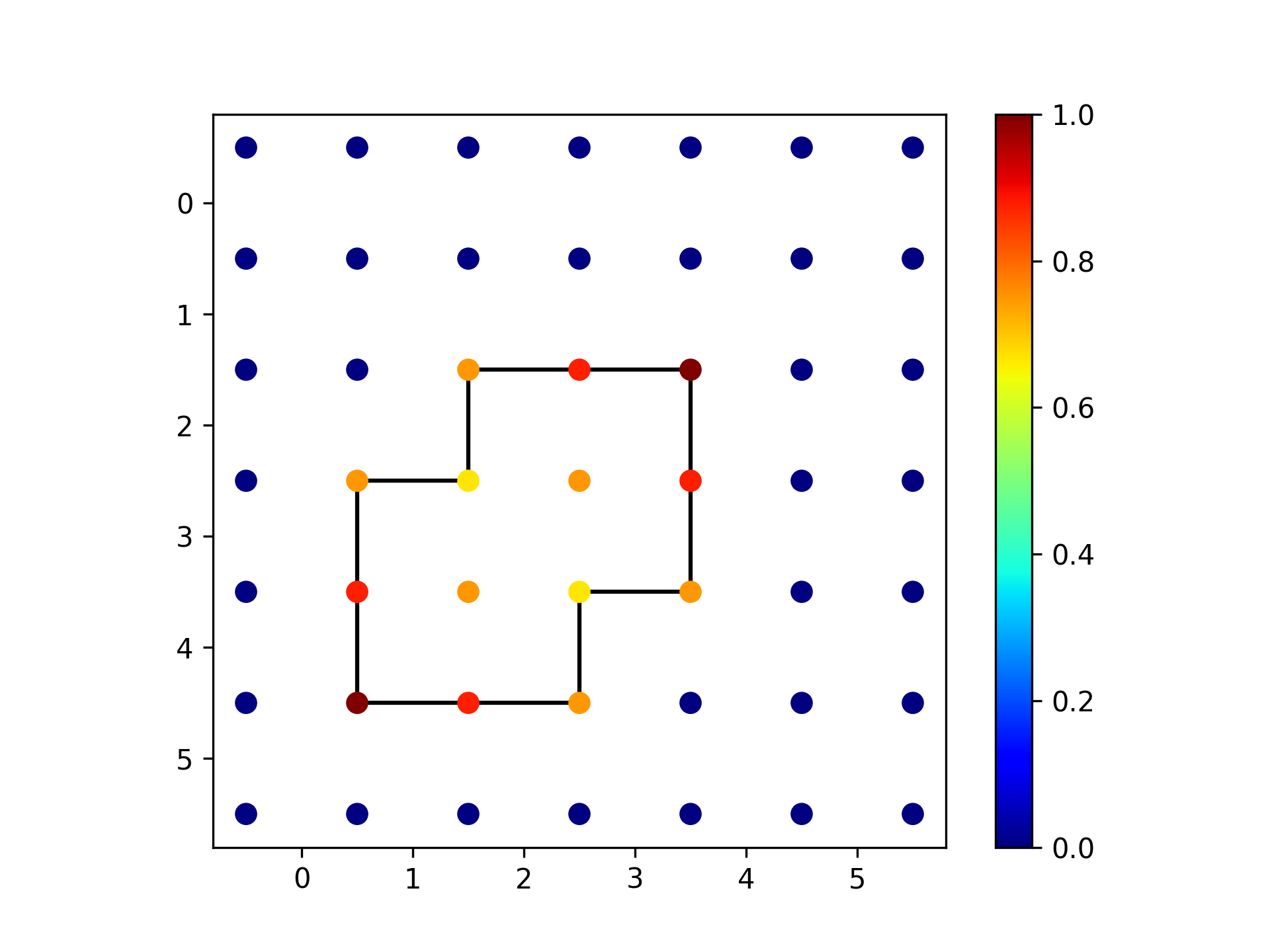} 
    \caption{Visualization of the values of the vertices of the squares in Fig. \ref{F: VF3}. The values of the vertices $u_{i,j}$ are set as the mean of the values of $v_{i,j}$.}
    \label{F: VF4}
\end{figure}

Let us now consider the smoothed trajectories: each grid point on the curve is associated with a smoothed velocity vector. 
We consider a pixel grid of dimension $n_1\times n_2$ with grid size $h$ and check the pixels to which the grid points of the trajectories belong as shown in Fig. \ref{F: VF1}. We then use those squares in the definition of Dirichlet conditions inside $\tilde{\Omega}_{cut}$. We will consider the 2 components of the vectors and the vector lengths separately; throughout the text, the term sparse sample will indicate the position and value of one of those. Let us denote by $S_{p}$, $p=1,...,s$ the squares of size $h$ where the sparse samples are located. We also assume that none of $S_{p}$ touches the boundary of $\tilde{\Omega}_{cut}$.
Let us define
\begin{equation}
\tau_p=\{v_{i,j}: v_{i,j} \text{ is value in center of } S_p \text{ for some } p\in\{1,...s\} \},
\end{equation}
and 
\begin{equation}
\tau_h=\{u_{i,j}: u_{i,j} \text{ is value in vertex of } S_p \text{ for some } p\in\{1,...s\} \}.   
\end{equation}
The values of the Dirichlet conditions will be constructed by using the values of the sparse samples; if more grid points of the curve(s) belong to the same pixel, we consider as $v_{i,j}$ the mean value of the sparse samples belonging to that pixel. The visualization of the described process for the domain in Fig. \ref{F: VF1} is shown in Fig. \ref{F: VF2}; the values of the other pixels are set to 0 for better visualization.
To prove the existence and uniqueness of a weak solution for the minimization problem we will consider, we need to have a Lipschtiz domain; this is not the case for the domain shown in Fig. \ref{F: VF2} if we would cut away those 3 squares. Therefore, we add another step to the construction of the model domain. If the squares that we would cut away would have in common only one vertex, we add another two squares and set their values as the mean value of the two squares having in common that vertex. In Fig. \ref{F: VF3} is visualized the result for the domain shown in Fig. \ref{F: VF2}.  
 Let us define
\begin{equation}
    S=\bigcup_{p=1}^{s} \bar{S}_p
\end{equation}
and
\begin{equation}
\Gamma^1=\partial S,
\end{equation}
Indicate by
\begin{equation}
    \Omega=\tilde{\Omega}_{cut}/S,
\end{equation}
so that
\begin{equation}
    \partial \Omega=\Gamma^{2}\cup \Gamma^1.
\end{equation}
The last step leading to the definition of Dirichlet boundary conditions is the definition of the values $u_{i,j}\in\tau_h$. To do that, we consider the squares $S_p$ of which $u_{i,j}$ is the value of one of the vertices; to each such $S_p$ a value $v_{i,j}\in\tau_p$ is associated. Therefore, we set $u_{i,j}$ as the mean of those values $v_{i,j}$. In Fig. \ref{F: VF4} is shown the result for the example in Fig. \ref{F: VF3}.\\
We finally define the continuous function
\begin{equation}
    g:\Gamma_1\rightarrow\mathbb{R},
    \label{E: Def g}
\end{equation}
on the boundary $\Gamma_1$ in the following way. If two vertices $u_{i,j}$ of an edge have the same value, then the value of $g$ is set as constant on the edge. If the values of the vertices $u_{i,j}$ differ, then we define $g$ on the edge as linearly changing from one value to the other.

Once we have defined the model domain, we consider the following minimization problem for the two vector components and the vector lengths.

\begin{equation}
\label{E: minimization problem}
    \min_{u} \frac{1}{2}\int_{\Omega} \lVert \nabla u \rVert^2,
\end{equation}
where $u: \Omega\rightarrow \mathbb{R}$ is a function defined on $\Omega \subseteq \mathbb{R}^2$ and $E[u]=\frac{1}{2}\int_{\Omega} \lVert \nabla u \rVert^2$ is the Dirichlet energy. Solving (\ref{E: minimization problem}) leads to
\begin{equation}
-\Delta u^i=0 ~~~ i=1,2,3
\end{equation}
with appropriate boundary conditions. Here, $u^1=v_x$ and $u^2=v_y$ are the two vector components, while $u^3=L$ is the vector length. In particular, we set
\begin{itemize}
    \item Zero Neumann boundary conditions on $\Gamma^{2}$
    \begin{equation}
        \frac{\partial u}{\partial \nu}(x)=0 ~~~~ x\in \Gamma^{2},
    \end{equation}
    where $\nu$ is the unit outer normal vector to the boundary $\Gamma^{2}$.
    \item Dirichlet boundary conditions on $\Gamma^1$
    \begin{equation}
    \begin{split}
        u(x)&=g(x) ~~~~ x\in \Gamma^1,  
    \end{split}
    \end{equation}
    where $g$ is defined in (\ref{E: Def g}).
\end{itemize}
 The proof of the existence and uniqueness of a weak solution can be done in the same way as in \cite{LupiVF}. Indeed, we constructed the Lipschtiz domain $\Omega$, therefore we can apply the theory from \cite{Necas2012DirectMI, rektorys2012variational}.
When the double-Laplacian operator for the two vector components is used, the direction of the vectors changes smoothly from one to another, depending on the Dirichlet conditions. However, it is important to note that the Laplace equation averages the values of neighboring points. As previously discussed in \cite{DesignTangentVF,LupiVF}, while the double-Laplacian effectively maintains the desired directionality of the vectors, it tends to shorten their lengths. Since we want to find the vector field from the provided data, it is logical to consider reconstructing the vector lengths based on the given Dirichlet conditions. 
Consequently, we propose to apply the same minimization problem (\ref{E: minimization problem}) for determining the vector lengths as well.

The discretization of the boundary value problem
 \begin{equation}
 \left \{
\begin{array}{l}
-\Delta u(x)=0, ~~~x\in \Omega,\\
\frac{\partial u}{\partial \nu}(x)=0, ~~~~~~~ x\in \Gamma^{2},\\
u(x)=g(x), ~~~~ x\in \Gamma^1,\\
\end{array}  
\right.
\end{equation}
is based on \cite{article}. In particular, for the vertices $u_{i,j}$ which are not vertices of any $S_p$, $p\in\{1,...,s\}$ the standard discretization of the Laplace equation is used. On the other hand, in vertices of $\Gamma^1$ we prescribe Dirichlet values, and on the boundary $\Gamma^{2}$ we consider the discretization of zero Neumann boundary condition similar to \cite{article}. 
%Note that, in our case, the pixels corresponding to $\Gamma^2$ are those intersected by the boundary $\Gamma^2$, rather than the pixels on the image boundary. Consequently, we adapt the discretization by applying the zero Neumann boundary conditions specifically to these pixels. 
In Section \ref{SS: Velocity} we recall how to find the smoothed velocities, cf.  \cite{Lupi2022MacrophagesTS}, and results for the reconstruction of the wound attractant field are shown in Section \ref{SS: VF Results}.

 \subsection{Velocity vectors on smooth trajectories}
\label{SS: Velocity}
To get the new velocity estimation consider $(\Delta_t)_j$ the time between the two endpoints $\textbf{x}(u_{j-1})$ and $\textbf{x}(u_j)$ of the segment $j$. In the beginning, every segment has the same value of $(\Delta_t)_j$; later, if a segment does not disappear, we consider the same $(\Delta_t)_j$ as in the beginning. On the other hand, if a segment  $j$ disappears, we add $\frac{(\Delta_t)_j}{2}$ to the first previous not disappeared segment $j_1<j$,
\[
 (\Delta_t)_{j_1}=(\Delta_t)_{j_1}+\frac{(\Delta_t)_j}{2},
\]
 and $\frac{(\Delta_t)_j}{2}$ to the first following not disappeared segment $j_2>j$
 \[
 (\Delta_t)_{j_2}=(\Delta_t)_{j_2}+\frac{(\Delta_t)_j}{2}.
\]
To find the new velocity estimation for the grid points on the evolved curve, we apply the following algorithm. Set $j=1$, $i=1$ and repeat until $j=M-1$:
\begin{itemize}
 \item If $r_{j}^{m+1}> 0$ 
 \begin{equation}
  \lvert \textbf{v} \rvert=\frac{L^{d}_{j}}{(\Delta_t)_j},
 \end{equation}
 where $L^{d}_{j}$ was defined in (\ref{E: L^d_j}).\\
 \item for $i=\mathcal{I}(u_{j-1}),...,\mathcal{I}(u_j)$
\begin{equation}
 \textbf{v}_{i}^{m+1}=\lvert \textbf{v} \rvert\frac{(\textbf{x}_{i}^{m+1}-\textbf{x}_{i-1}^{m+1})}{\lvert \textbf{x}_{i}^{m+1}-\textbf{x}_{i-1}^{m+1}\rvert}.
\end{equation}
\end{itemize}
If the segment does not disappear, we consider the grid points inside that segment to have constant velocity equal to the new length of the segment divided by the time interval $(\Delta_t)_j$. 
On the other hand, if the $j$-th segment disappears, the time $(\Delta_t)_j$ is distributed between the previous and following segments that have not disappeared.\\    

\subsection{Results}
\label{SS: VF Results}
For all the considered datasets we smoothed the trajectories with the model described in Section \ref{S: Mathematical model}. Then, we computed the smoothed velocities with the algorithm described in Section \ref{SS: Velocity} and used them as sparse samples to reconstruct the wound attractant field driving macrophages to the site of the wound. 
\begin{figure}[H]
    \centering
    \includegraphics[width=0.7\linewidth]{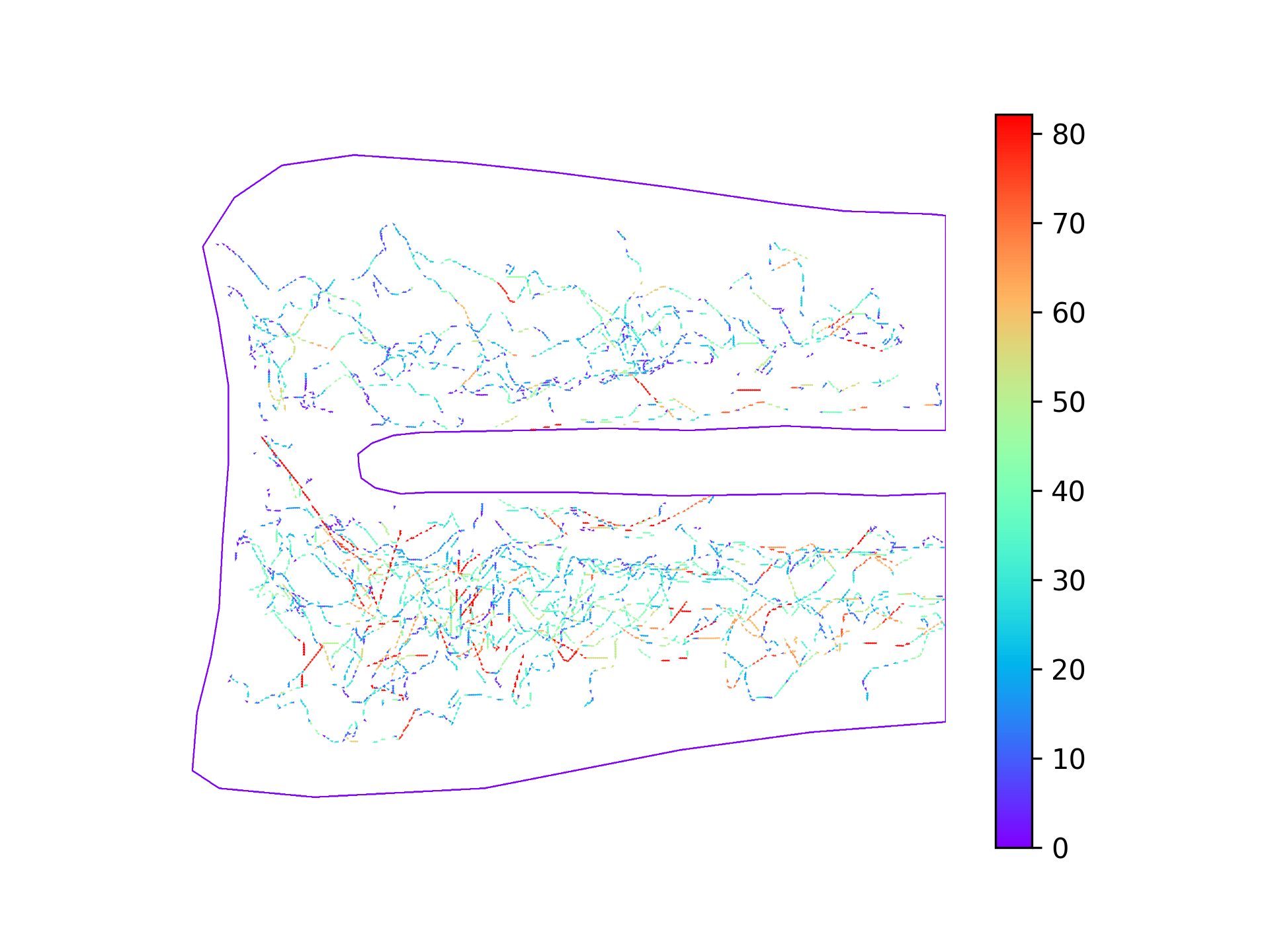}
    \includegraphics[width=0.7\linewidth]{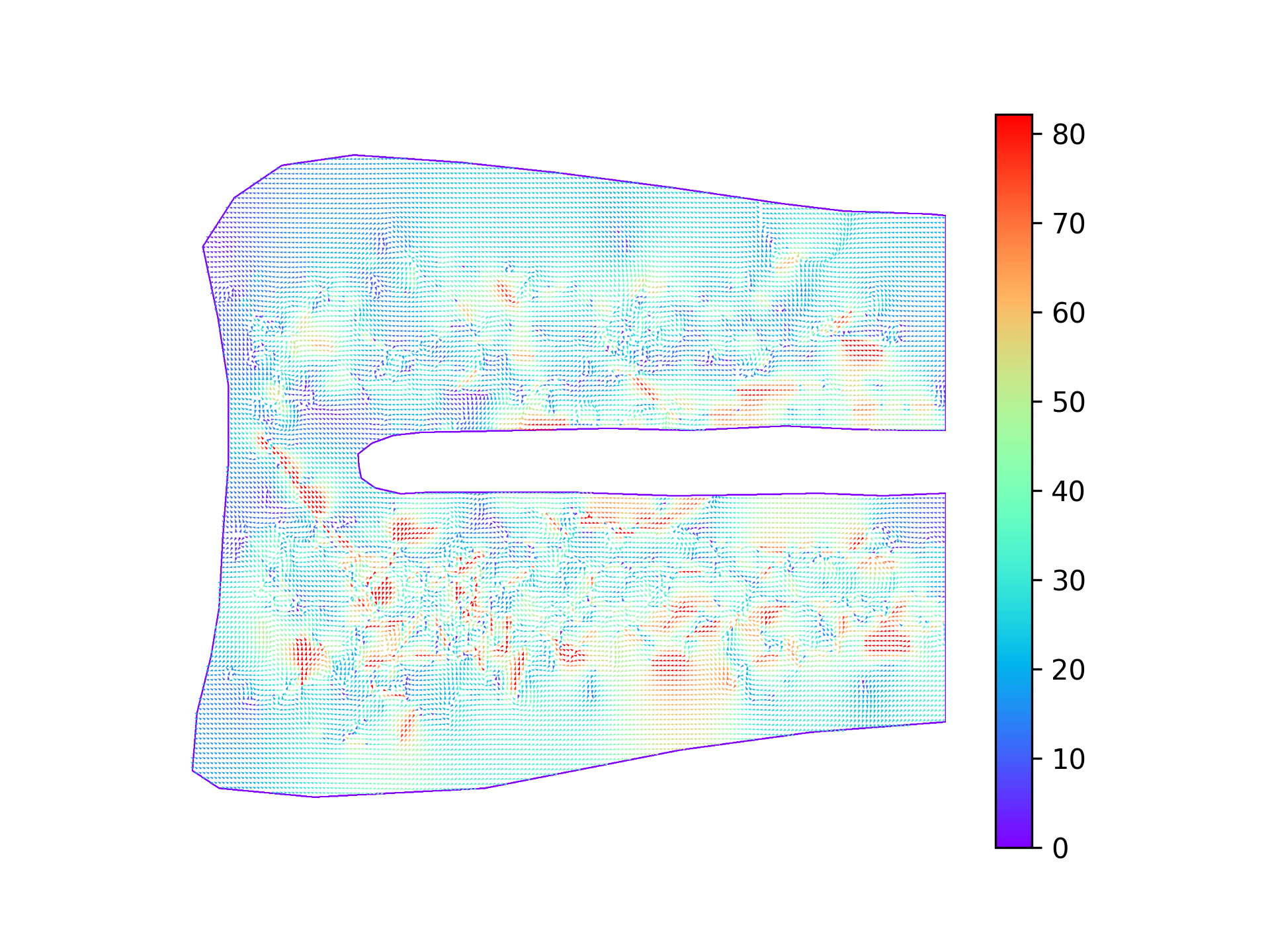} 
    \caption{Results for the dataset shown in Fig. \ref{F: Resul3_theory}. Top: Prescribed velocities on sparse samples obtained from smoothed macrophage trajectories. Bottom: Reconstructed vector field by triple-Laplacian.}
    \label{F: Resul3}
\end{figure}
\begin{figure}[H]
    \centering
   \includegraphics[width=0.7\linewidth]{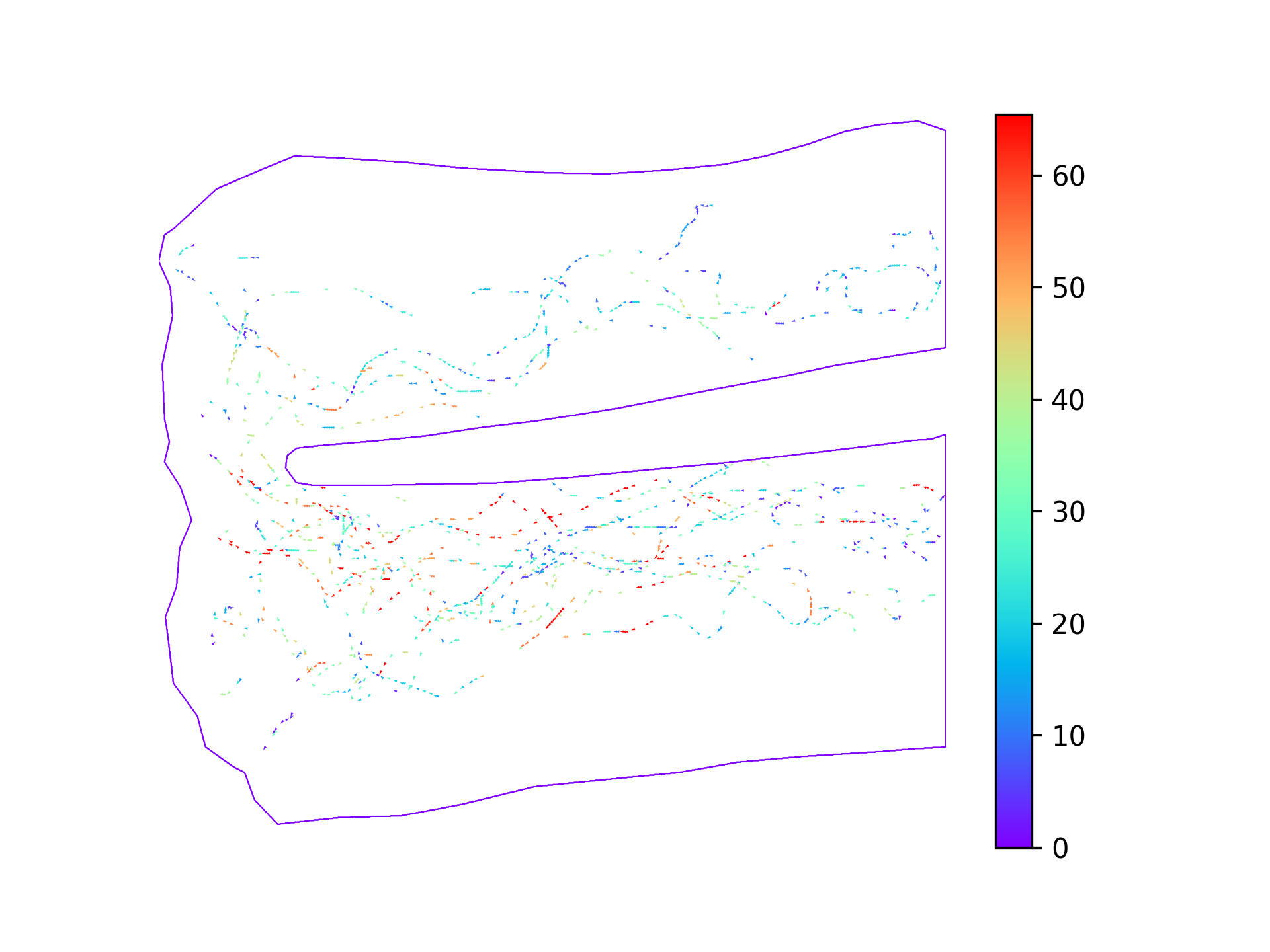}
    \includegraphics[width=0.7\linewidth]{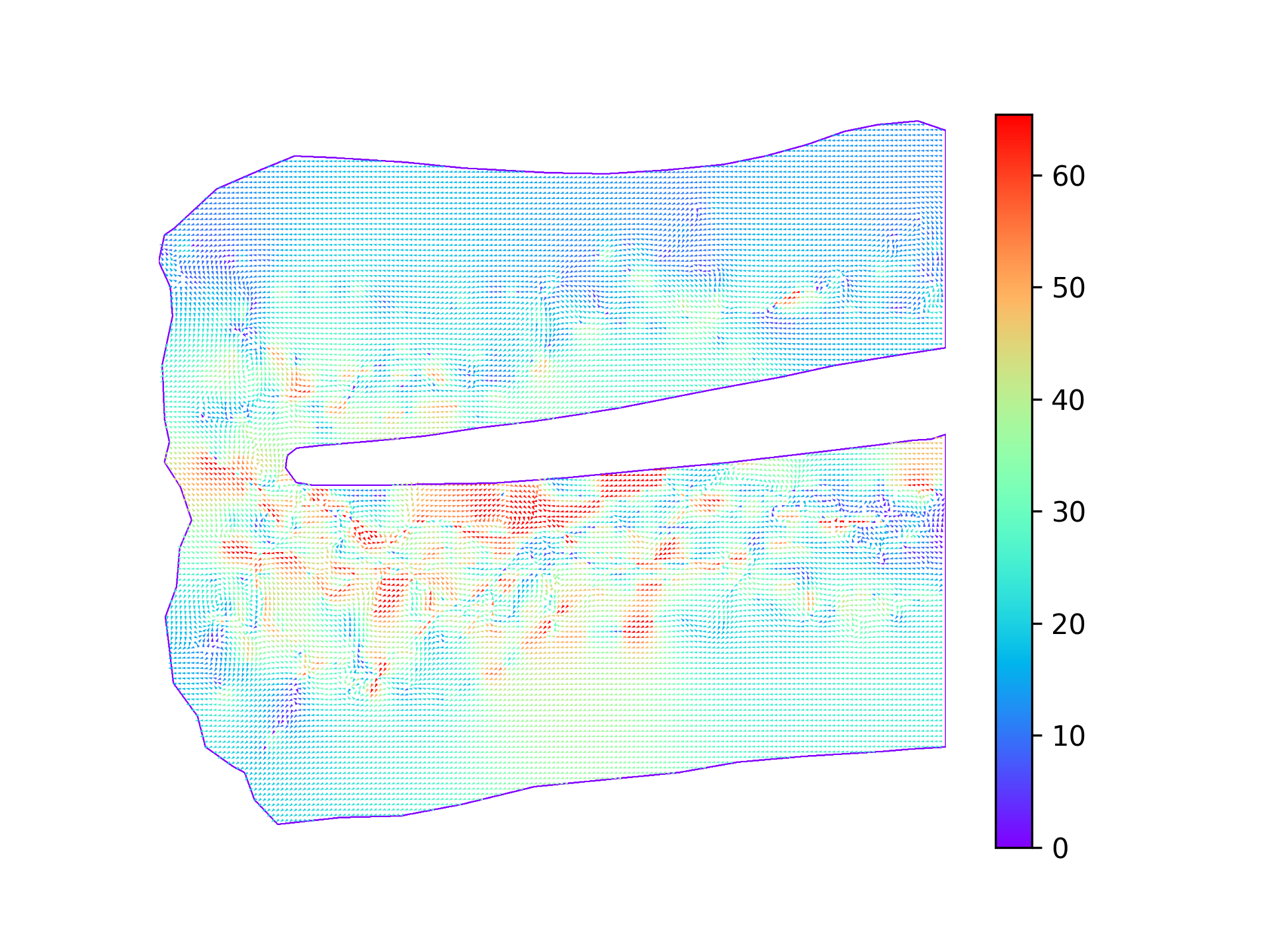}
    \caption{Results for the dataset shown in Fig. \ref{F: Resul4_theory}. Top: Prescribed velocities on sparse samples obtained from smoothed macrophage trajectories. Bottom: Reconstructed vector field by triple-Laplacian}
    \label{F: Resul4}
\end{figure}
\begin{figure}[H]
    \centering
    \includegraphics[width=0.8\linewidth]{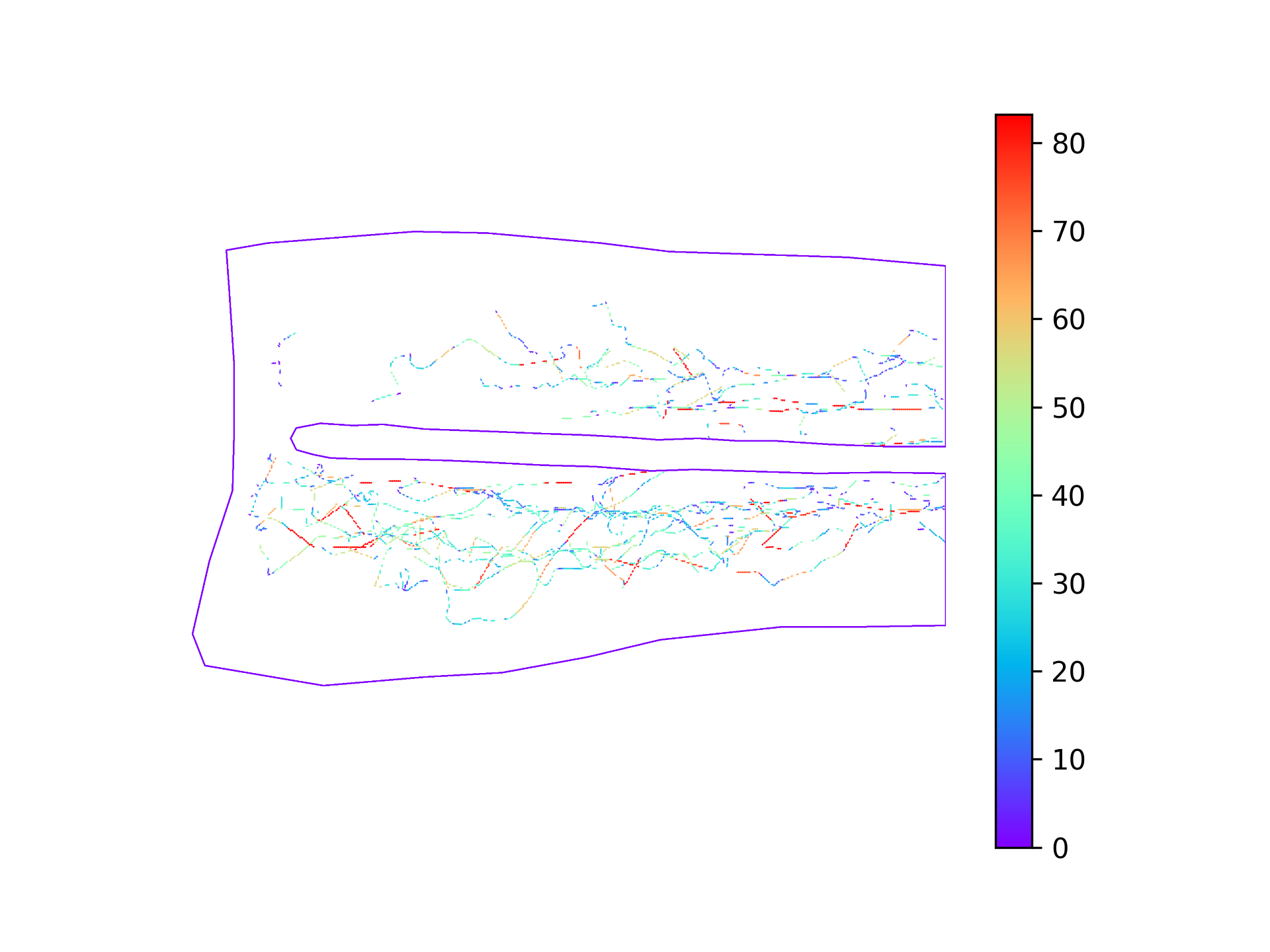}
    \includegraphics[width=0.8\linewidth]{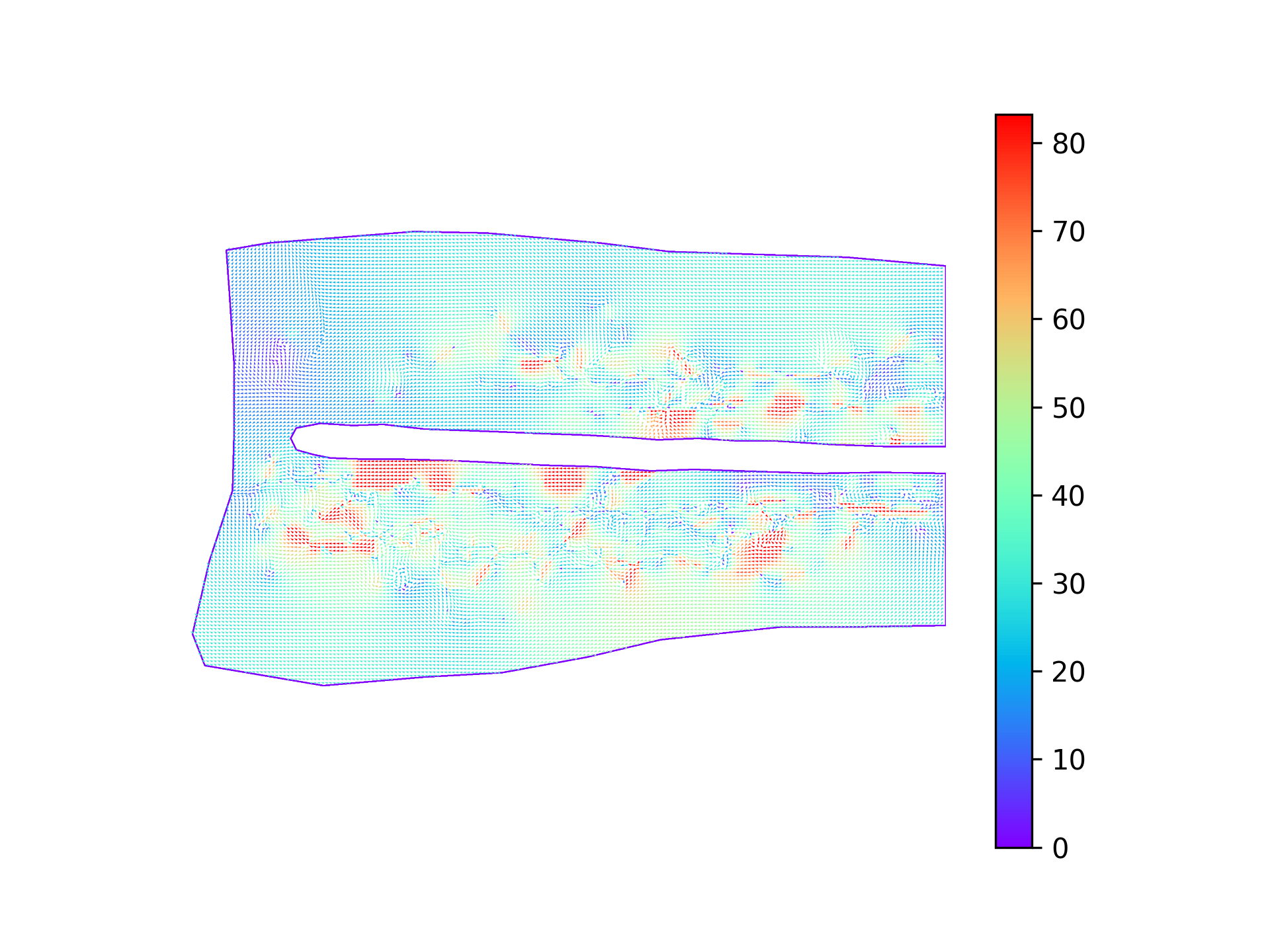}
    \caption{Results for the dataset shown in Fig. \ref{F: Resul5_theory}. Top: Prescribed velocities on sparse samples obtained from smoothed macrophage trajectories. Bottom: Reconstructed vector field by triple-Laplacian}
    \label{F: Resul5}
\end{figure}
Figs. \ref{F: Resul3}, \ref{F: Resul4}, \ref{F: Resul5} show the results for the $3$ different datasets whose model domain is visualized in Figs. \ref{F: Resul3_theory}, \ref{F: Resul4_theory}, \ref{F: Resul5_theory} respectively. For all of them, the wound is on the left side. On the top, are visualized the sparse samples, and on the bottom, one can see the reconstructed vector field; the color represents the length of the vectors.
\begin{figure}[H]
    \centering
    \includegraphics[width=0.7\linewidth]{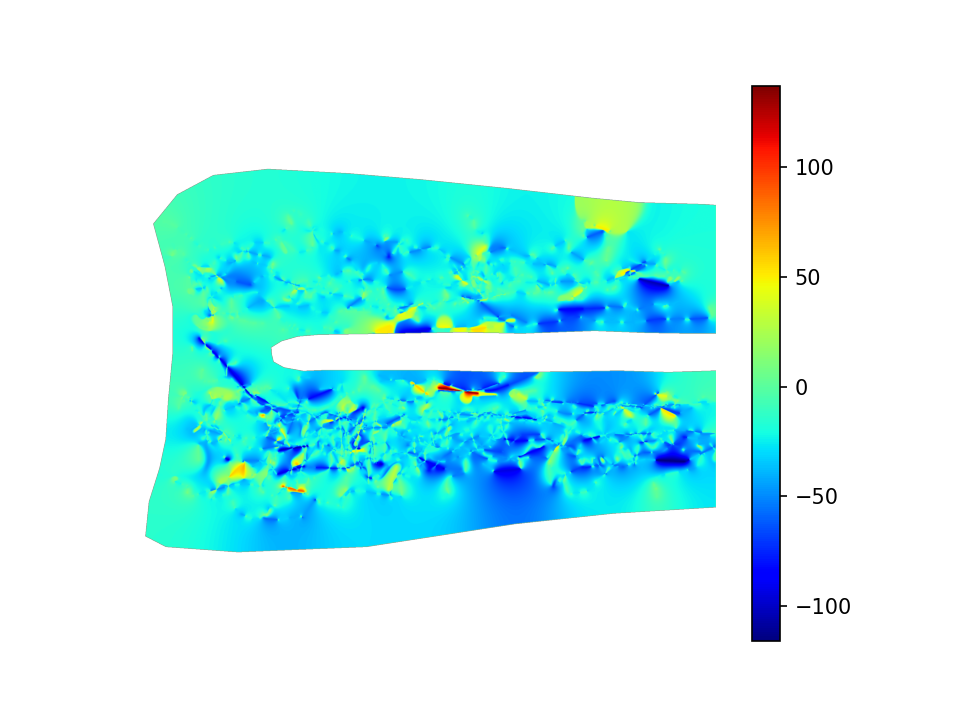}
    \includegraphics[width=0.7\linewidth]{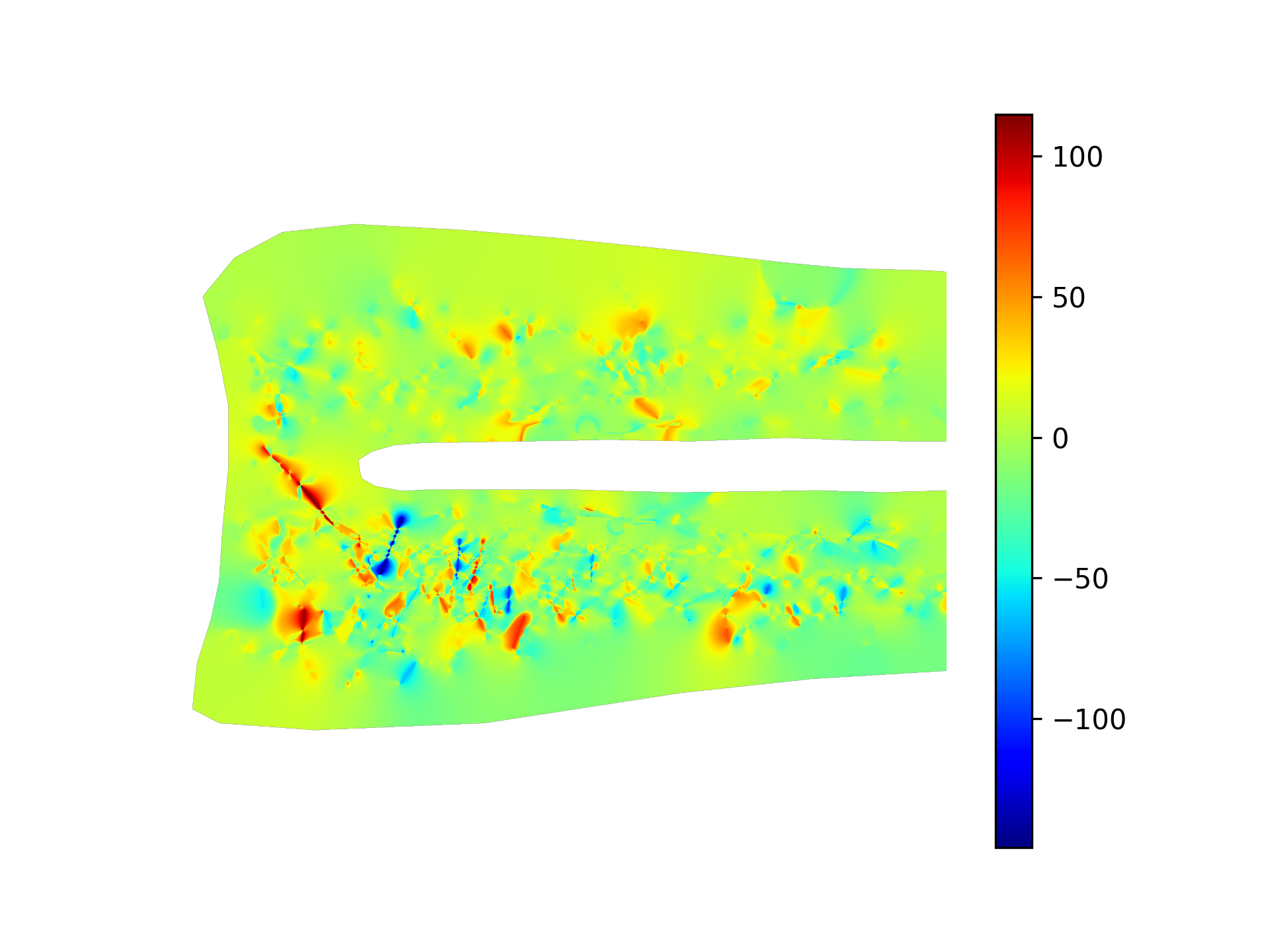}
    \caption{Projection in $x$ and $y$ direction of the reconstructed vector field shown in Fig. \ref{F: Resul3}. Top: $x$ direction. Bottom: $y$ direction.}
    \label{F: Resul3_XY}
\end{figure}
\begin{figure}[H]
    \centering
    \includegraphics[width=0.7\linewidth]{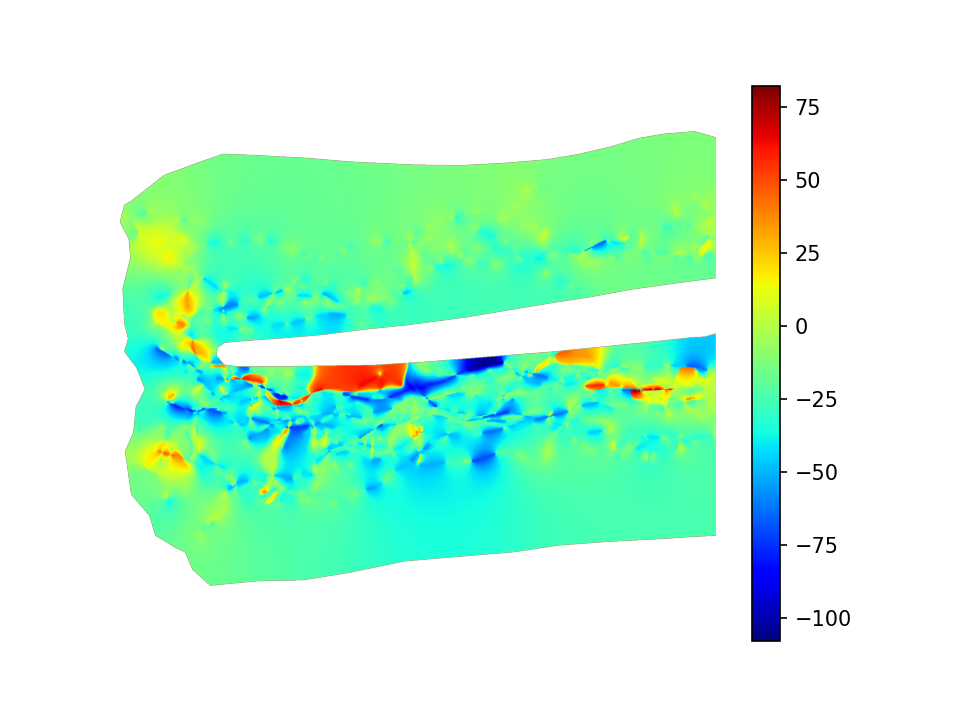}  \includegraphics[width=0.7\linewidth]{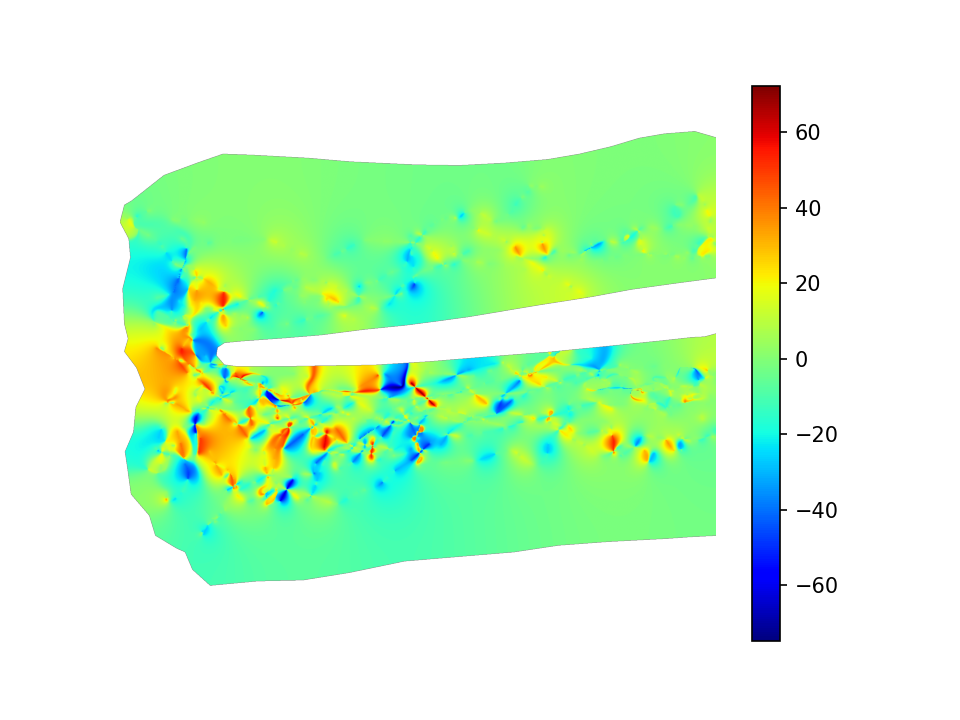}
    \caption{Projection in $x$ and $y$ direction of the reconstructed vector field shown in Fig. \ref{F: Resul4}. Top: $x$ direction. Bottom: $y$ direction.}
    \label{F: Resul4_XY}
\end{figure}
\begin{figure}[H]
    \centering
    \includegraphics[width=0.7\linewidth]{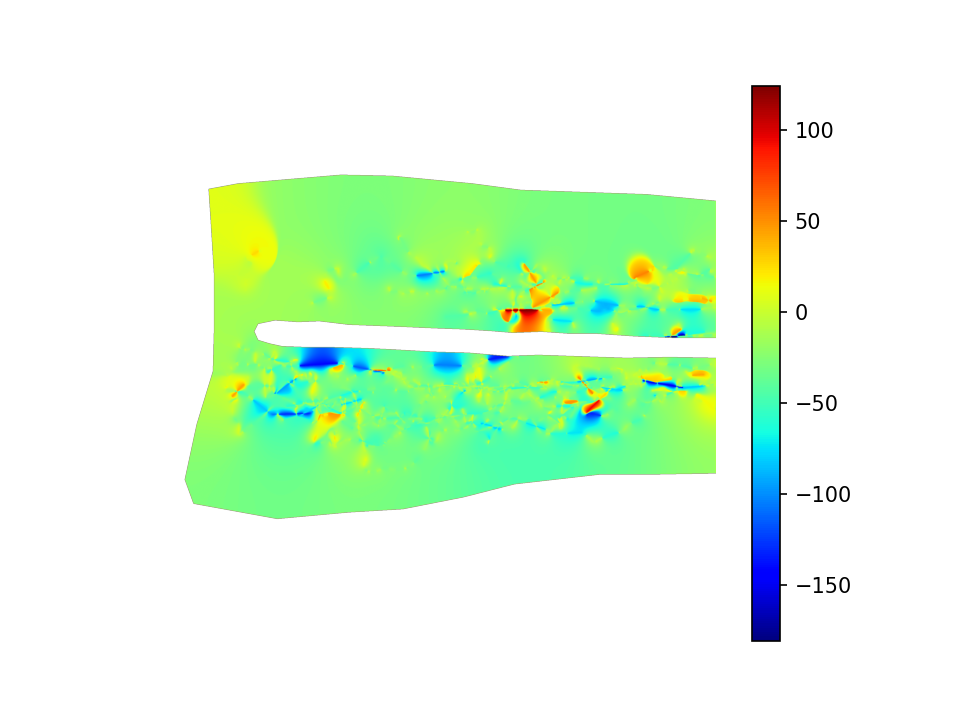}
    \includegraphics[width=0.7\linewidth]{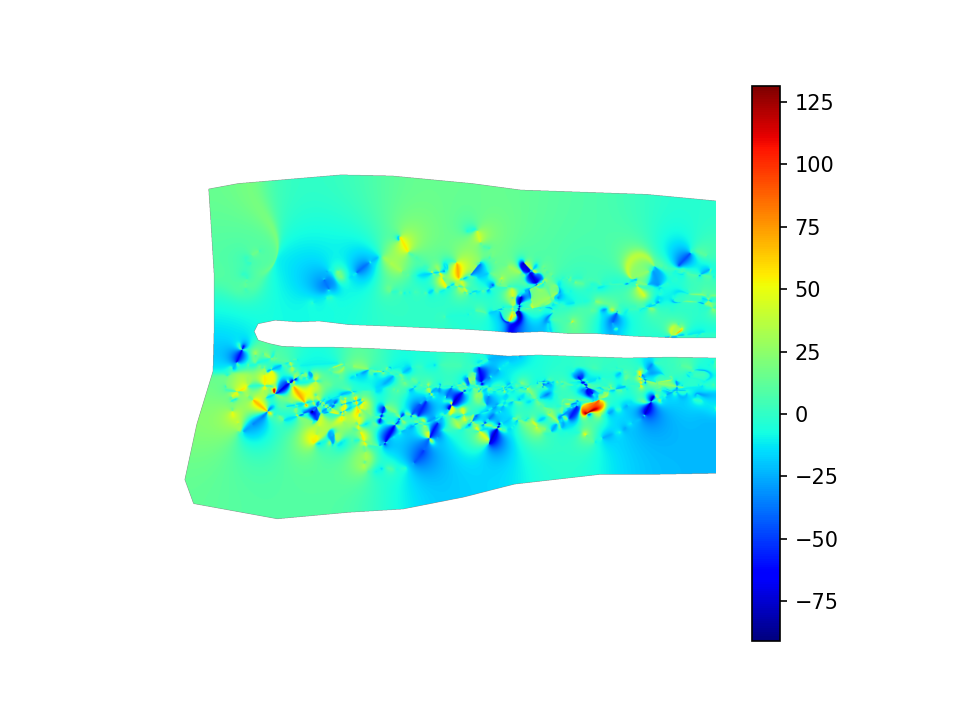}
    \caption{Projection in $x$ and $y$ direction of the reconstructed vector field shown in Fig. \ref{F: Resul5}. Top: $x$ direction. Bottom: $y$ direction.}
    \label{F: Resul5_XY}
\end{figure}
Figs. \ref{F: Resul3_XY}, \ref{F: Resul4_XY}, \ref{F: Resul5_XY} show the projection of vectors in Figs. \ref{F: Resul3}, \ref{F: Resul4}, \ref{F: Resul5} in $x$ (horizontal) and $y$ (vertical) direction. The images on top show the result for $x$ direction. The blue color represents the regions of higher directionality; indeed, since the wound is on the left side, the more negative the value the faster the cell is moving toward the wound. The images on the bottom show the results for the $y$ direction. The blue color represents fast motion from top to bottom. As one can notice, the fast motion toward the wound ($x$ direction) is in similar parts of the fish body for the $3$ datasets, on the bottom and close to the central part of the fish body where the macrophages can not enter. On the other side, the motion in the $y$ direction seems to be random. 
\section{Conclusions}
In this paper, we proposed a new workflow to analyze macrophage trajectories and reconstruct the velocity vector field driving macrophages toward the wound during wound healing. First, we smooth the trajectories; the smoothing algorithm allows us to separate the random parts from the directional parts of the motion. We considered two ways of defining the random parts. The first considers the evolution of the lengths of the original segments and defines a segment as random if it disappears during the smoothing process. This usually happens in regions of high curvature where the cell is moving back and forth without a clear directionality. The second approach to defining the random parts takes into account the self-intersections of the trajectory. We detect the points where the curve intersects itself and define as random the segments that belong to a self-intersecting part of the trajectory. Once we had the random sub-trajectories, we examined them using the mean squared displacement (MSD). This is a common tool used to characterize the type of diffusion motion followed by a random walker. The analysis showed that the random motion is subdiffusive; subdiffusion usually occurs in crowded environments where the random walk is not free to move in any direction but is constrained by the environment. Finally, we used the velocities found on the smoothed trajectories as sparse samples to reconstruct the wound attractant field that drives the macrophages during wound healing. The model we proposed minimizes the Dirichlet energy for the two vector components and the vector length with appropriate boundary conditions. In solving the minimization problem, we ended up solving the Laplace equation with Dirichlet conditions on the sparse samples and zero Neumann boundary conditions on the domain boundary. The result of the reconstruction is an interpolation/extrapolation process in which the information is extrapolated from the Dirichlet conditions where the vectors do not influence each other and interpolated in the regions where the vectors do influence each other. We developed such a model to extract the information about the velocity vector field only from the Dirichlet conditions, without the need to make any assumptions about the shape or properties of the wound attractant field.

Modeling a phenomenon such as motion during wound healing \textit{in vivo} is challenging due to the large number of possible parameters that can influence the motion of the cells. 
Whilst the motion is influenced by the wound attractant field, other causes may be responsible for the shape of the trajectories of these cells. Firstly, the motion of cells can vary significantly depending on the environment in which they are moving. Indeed, in \cite{kubala2021diffusion} the authors studied the effects of a crowded environment on the motion of random walkers subject to a constant velocity vector field. They showed that the effect of constant drift increases the chances of trapping and that the trajectories show a high degree of variability. Therefore, the interplay between drift and environment can be responsible for the emergence of both superdiffusive and subdiffusive behavior. For this reason, it's important to separate the parts of the movement that are not influenced by the wound attractant field and study them separately. A possible future direction that could contribute to the understanding of macrophage movement is the segmentation of the surrounding environment. The information about the environment in which the cells move could be used both to determine a model for the random migration of macrophages and to reconstruct the velocity vector field taking into account more precise boundary conditions within the domain.

Another interesting aspect is the study of macrophage movement under normal conditions, i.e. without the presence of a wound. In this case, the overall trajectories of the macrophages could be considered as random walks, and, with a sufficiently large data set, a model of migration under normal conditions could perhaps be determined. Finally, the smoothed trajectories and extracted random parts can be used in the classification of different states of macrophages. While many studies have reported differences in the morphology and dynamics of macrophage subtypes, such as M1, non-activated M1, and M2 macrophages, comprehensive quantitative analyses remain limited. By combining the segmented shapes of macrophages with properties of their smoothed trajectories and random motion components, a high-dimensional morpho-kinetic feature space can be defined and used for the graph-PDE supervised classification to differentiate between M1 and M2 macrophages \cite{Mikula2024}.

\begin{acknowledgement}
  This work has received funding from the European Union's Horizon 2020 research and innovation programme under the Marie Sk\l{}odowska-Curie grant agreement No 955576 and by the grants APVV-23-0186 and VEGA 1/0249/24.   
\end{acknowledgement}

\bibliographystyle{plain}
\bibliography{Bibliography}

\begin{thebibliography}{10}

\bibitem{ambroz2019numerical}
M.~Ambroz, M.~Bala{\v{z}}ovjech, M.~Medl’a, and K.~Mikula.
\newblock Numerical modeling of wildland surface fire propagation by evolving surface curves.
\newblock {\em Advances in Computational Mathematics}, 45:1067--1103, 2019.

\bibitem{ambroz2020semi}
M.~Ambroz, M.~Kollár, and K.~Mikula.
\newblock Semi-implicit scheme for semi-automatic segmentation in naturasat software.
\newblock In {\em Proceedings of the Conference Algoritmy}, pages 171--180, 2020.

\bibitem{Arcizet2008Temporal}
D.~Arcizet, B.~Meier, E.~Sackmann, J.O. R{\"a}dler, and D.~Heinrich.
\newblock Temporal analysis of active and passive transport in living cells.
\newblock {\em Physical review letters}, 101(24):248103, 2008.

\bibitem{barros2017live}
F.~Barros-Becker, P.Y. Lam, R.~Fisher, and A.~Huttenlocher.
\newblock Live imaging reveals distinct modes of neutrophil and macrophage migration within interstitial tissues.
\newblock {\em Journal of cell science}, 130(22):3801--3808, 2017.

\bibitem{doi:10.1137/080716396}
Z.~Belhachmi, D.~Bucur, B.~Burgeth, and J.~Weickert.
\newblock How to choose interpolation data in images.
\newblock {\em SIAM Journal on Applied Mathematics}, 70(1):333--352, 2009.

\bibitem{de2007cellular}
K.~De Bruin, N.~Ruthardt, K.~Von Gersdorff, R.~Bausinger, E.~Wagner, M.~Ogris, and C.~Br{\"a}uchle.
\newblock Cellular dynamics of egf receptor-targeted synthetic viruses.
\newblock {\em Molecular Therapy}, 15(7):1297--1305, 2007.

\bibitem{daumas2003confined}
F.~Daumas, N.~Destainville, C.~Millot, A.~Lopez, D.~Dean, and L.~Salom{\'e}.
\newblock Confined diffusion without fences of a g-protein-coupled receptor as revealed by single particle tracking.
\newblock {\em Biophysical Journal}, 84(1):356--366, 2003.

\bibitem{dieterich2008anomalous}
P.~Dieterich, R.~Klages, R.~Preuss, and A.~Schwab.
\newblock Anomalous dynamics of cell migration.
\newblock {\em Proceedings of the National Academy of Sciences}, 105(2):459--463, 2008.

\bibitem{Ellett}
F.~Ellett, L.~Pase, J.~Hayman, A.~Andrianopoulos, and G.~Lieschke.
\newblock Mpeg1 promoter transgenes direct macrophage-lineage expression in zebrafish.
\newblock {\em Blood}, 117(4):e49--e56, 2011.

\bibitem{DesignTangentVF}
M.~Fisher, P.~Schr\"{o}der, M.~Desbrun, and H.~Hoppe.
\newblock Design of tangent vector fields.
\newblock {\em ACM Transactions on Graphics (TOG)}, 26(3):56–es, 2007.

\bibitem{friedl2008interstitial}
P.~Friedl and B.~Weigelin.
\newblock Interstitial leukocyte migration and immune function.
\newblock {\em Nature immunology}, 9(9):960--969, 2008.

\bibitem{modes2010matrix}
E.~Van Goethem, R.~Poincloux, F.~Gauffre, I.~Maridonneau-Parini, and V.~Le Cabec.
\newblock Matrix architecture dictates three-dimensional migration modes of human macrophages: differential involvement of proteases and podosome-like structures.
\newblock {\em Journal of immunology}, 184(2):1049--1061, 2010.

\bibitem{harris2012generalized}
T.H. Harris, E.J. Banigan, D.A. Christian, et~al.
\newblock Generalized l{\'e}vy walks and the role of chemokines in migration of effector cd8+ t cells.
\newblock {\em Nature}, 486(7404):545--548, 2012.

\bibitem{hou1994removing}
T.Y. Hou, J.S. Lowengrub, and M.J.Shelley.
\newblock Removing the stiffness from interfacial flows with surface tension.
\newblock {\em Journal of Computational Physics}, 114(2):312--338, 1994.

\bibitem{huda2018levy}
S.~Huda, B.~Weigelin, K.~Wolf, et~al.
\newblock L{\'e}vy-like movement patterns of metastatic cancer cells revealed in microfabricated systems and implicated in vivo.
\newblock {\em Nature communications}, 9(1):4539, 2018.

\bibitem{kadirkamanathan2012neutrophil}
V.~Kadirkamanathan, S.R. Anderson, S.A. Billings, X.~Zhang, G.R. Holmes, et~al.
\newblock The neutrophil's eye-view: Inference and visualisation of the chemoattractant field driving cell chemotaxis in vivo.
\newblock {\em PloS one}, 7(4):e35182, 2012.

\bibitem{kimura1997numerical}
M.~Kimura.
\newblock Numerical analysis of moving boundary problems using the boundary tracking method.
\newblock {\em Japan journal of industrial and applied mathematics}, 14:373--398, 1997.

\bibitem{kubala2021diffusion}
P.~Kubala, M.~Cie{\'s}la, and B.~Dybiec.
\newblock Diffusion in crowded environments: Trapped by the drift.
\newblock {\em Physical Review E}, 104:044127, 2021.

\bibitem{kusumi1993confined}
A.~Kusumi, Y.~Sako, and M.~Yamamoto.
\newblock Confined lateral diffusion of membrane receptors as studied by single particle tracking (nanovid microscopy). effects of calcium-induced differentiation in cultured epithelial cells.
\newblock {\em Biophysical journal}, 65(5):2021--2040, 1993.

\bibitem{VectorfieldReconstruction}
M.~Lage, F.~Petronetto, A.~Paiva, H.~Lope, T.~Lewiner, and G.~Tavares.
\newblock Vector field reconstruction from sparse samples with applications.
\newblock In {\em 19th Brazilian Symposium on Computer Graphics and Image Processing}, pages 297--306, 2006.

\bibitem{lee1991direct}
G.M. Lee, A.~Ishihara, and K.A. Jacobson.
\newblock Direct observation of brownian motion of lipids in a membrane.
\newblock {\em Proceedings of the National Academy of Sciences}, 88(14):6274--6278, 1991.

\bibitem{li2012live}
L.~Li, B.~Yan, Y.Q. Shi, W.Q. Zhang, and Z.L. Wen.
\newblock Live imaging reveals differing roles of macrophages and neutrophils during zebrafish tail fin regeneration.
\newblock {\em Journal of Biological Chemistry}, 287(30):25353--25360, 2012.

\bibitem{LupiVF}
G.~Lupi and K.~Mikula.
\newblock Vector field reconstruction from sparse samples by triple-laplacian.
\newblock In {\em Proceedings of the Conference Algoritmy}, pages 85--98, 2024.

\bibitem{Lupi2022MacrophagesTS}
G.~Lupi, K.~Mikula, and S.A. Park.
\newblock Macrophages trajectories smoothing by evolving curves.
\newblock {\em Tatra Mountains Mathematical Publications}, 86(1), 2023.

\bibitem{article}
K.~Mikula, M.~Ambroz, and R.~Moko\v{s}ov\'a.
\newblock What was the river ister in the time of strabo? a mathematical approach.
\newblock {\em Tatra Mountains Mathematical Publications}, 80:71--118, 2021.

\bibitem{mikula2011inflow}
K.~Mikula and M.~Ohlberger.
\newblock Inflow-implicit/outflow-explicit scheme for solving advection equations.
\newblock In {\em Finite Volumes in Complex Applications VI, Problems \& Perspectives}, pages 683--692, 2011.

\bibitem{mikula2021automated}
K.~Mikula, J.~Urb\'an, M.~Koll\'ar, M.~Ambroz, I.~Jarol\'imek, J.~{\v{S}}ib{\'\i}k, and M.~{\v{S}}ib{\'\i}kov{\'a}.
\newblock An automated segmentation of natura 2000 habitats from sentinel-2 optical data.
\newblock {\em Discrete and Continuous Dynamical Systems}, 14(3):1017--1032, 2021.

\bibitem{sevcovic2001evolution}
K.~Mikula and D.~\v{S}ev\v{c}ovi\v{c}.
\newblock Evolution of plane curves driven by a nonlinear function of curvature and anisotropy.
\newblock {\em SIAM Journal on Applied Mathematics}, 61(5):1473--1501, 2001.

\bibitem{mikula2004direct}
K.~Mikula and D.~\v{S}ev\v{c}ovi\v{c}.
\newblock A direct method for solving an anisotropic mean curvature flow of plane curves with an external force.
\newblock {\em Mathematical Methods in the Applied Sciences}, 27(13):1545--1565, 2004.

\bibitem{mikula2008simple}
K.~Mikula, D.~\v{S}ev\v{c}ovi\v{c}, and M.~Balazovjech.
\newblock A simple, fast and stabilized flowing finite volume method for solving general curve evolution equations.
\newblock {\em Communications in Computational Physics}, 7(1):195--211, 2010.

\bibitem{ZHAO2009674}
Z.~Min.
\newblock Electrical fields in wound healing— an overriding signal that directs cell migration.
\newblock {\em Seminars in Cell and Developmental Biology}, 20(6):674--682, 2009.

\bibitem{Mussa-Ivaldi}
F.A. Mussa-Ivaldi.
\newblock From basis functions to basis fields: vector field approximation from sparse data.
\newblock {\em Biological Cybernetics}, 67(6):479--489, 1992.

\bibitem{Necas2012DirectMI}
J.~Ne\v{c}as.
\newblock {\em Direct methods in the theory of elliptic equations}.
\newblock Springer, Heidelberg, 2012.

\bibitem{Mikula2024}
S.A. Park, , G.~Lupi, R.~Ozbilgic, Mai~Nguyen Chi, and Karol Mikula.
\newblock Morpho-kinetic profiling of different modes of macrophages during wound healing in zebrafish.
\newblock {\em Work in Progress}, 2024.
\newblock Manuscript in preparation.

\bibitem{Sora}
S.A. Park, T.~Sipka, Z.~Kriv\'a, M.~Nguyen-Chi, G.~Lutfalla, and K.~Mikula.
\newblock Segmentation-based tracking of macrophages in 2d+time microscopy movies inside a living animal.
\newblock {\em Computers in Biology and Medicine}, 153:106499, 2023.

\bibitem{regner2013anomalous}
B.M. Regner, D.~Vu{\v{c}}ini{\'c}, C.~Domnisoru, T.M. Bartol, M.W. Hetzer, D.M. Tartakovsky, and T.J. Sejnowski.
\newblock Anomalous diffusion of single particles in cytoplasm.
\newblock {\em Biophysical journal}, 104(8):1652--1660, 2013.

\bibitem{rektorys2012variational}
K.~Rektorys.
\newblock {\em Variational Methods in Mathematics, Science and Engineering}.
\newblock Springer, 2012.

\bibitem{ruthardt2011single}
N.~Ruthardt, D.C. Lamb, and C.~Br{\"a}uchle.
\newblock Single-particle tracking as a quantitative microscopy-based approach to unravel cell entry mechanisms of viruses and pharmaceutical nanoparticles.
\newblock {\em Molecular therapy: the journal of the American Society of Gene Therapy}, 19(7):1199--1211, 2011.

\bibitem{Schnlieb2015PartialDE}
C-B. Sch{\"o}nlieb.
\newblock {\em Partial Differential Equation Methods for Image Inpainting}.
\newblock Cambridge University Press, 2015.

\bibitem{vsevvcovivc2011evolution}
D.~{\v{S}}ev{\v{c}}ovi{\v{c}} and S.~Yazaki.
\newblock Evolution of plane curves with a curvature adjusted tangential velocity.
\newblock {\em Japan journal of industrial and applied mathematics}, 28:413--442, 2011.

\bibitem{Tamara}
T.~Sipka, R.~Peroceschi, M.~Groß, F.~Ellett, C.~Pescia, C.~Gonzalez, G.~Lutfalla, and M.~Nguyen-Chi.
\newblock Damage-induced calcium signaling and reactive oxygen species mediate macrophage activation in zebrafish.
\newblock {\em Frontiers in Immunology}, 12, 2021.

\bibitem{sokolov2012models}
I.M. Sokolov.
\newblock Models of anomalous diffusion in crowded environments.
\newblock {\em Soft Matter}, 8(35):9043--9052, 2012.

\bibitem{viswanathan2011physics}
G.M. Viswanathan, M.G.E.~Da Luz, E.P. Raposo, and H.E. Stanley.
\newblock {\em The physics of foraging: an introduction to random searches and biological encounters}.
\newblock Cambridge University Press, 2011.

\end{thebibliography}
\end{document}